\documentclass[11pt]{amsart}
\usepackage{amsmath,amsfonts,amssymb,amsxtra,latexsym,amscd,enumerate,amsthm}

\usepackage{latexsym}
\usepackage{epsfig}
\usepackage{verbatim}

\def\g{\gamma}

\def\a{\alpha}

\def\b{\beta}

\def\p{\Phi}
\def\v{\varphi}
\def\ep{\epsilon}

\def\ll{\lambda}
\def\l{\Lambda}
\def\o{\omega}

\def\n{\nu}
\def\R{\mathbb{R}}

\def\c{\mathcal{C}}
\def\C{\mathbb{C}}
\def\Z{\mathbb{Z}}

\def\D{\mathcal{D}}
\def\P{\mathbb{P}}
\def\p{\mathcal{P}}
\def\U{\mathcal{U}}
\def\V{\mathcal{V}}
\def\F{\mathcal{F}}
\def\J{\mathcal{J}}
\def\I{\mathcal{I}}
\def\M{\mathcal{M}}
\def\S{\mathcal{S}}
\def\H{\mathcal{H}}

\def\A{\mathcal{A}}
\def\B{\mathcal{B}}
\def\L{\mathcal{L}}
\def\TT{\mathbb{T}}
\def\N{\mathbb{N}}
\def\cN{\mathcal{N}}
\def\n{\mathcal{N}}

\def\W{\mathcal{W}}
\def\Z{\mathbb{Z}}
\def\beq{\begin{equation}}
\def\eeq{\end{equation}}

\def\beq{\begin{equation}}
\def\eeq{\end{equation}}

\newtheorem{d0}{Definition}
\newtheorem{o0}{Observation}
\newtheorem{c0}{Corollary}
\newtheorem{t1}{Theorem}
\newtheorem{l1}{Lemma}
\newtheorem{p1}{Proposition}

\begin{document}
\title[]{The pointwise convergence of Fourier Series (I)\\
On a conjecture of Konyagin}

\author{Victor Lie}

\date{\today}
\address{Department of Mathematics, Purdue, IN 46907 USA}

\email{vlie@math.purdue.edu, vlie@math.princeton.edu}

\address{Institute of Mathematics of the
Romanian Academy, Bucharest, RO 70700, P.O. Box 1-764, Romania.}

\thanks{The author was supported by NSF grant DMS-1200932.}


\maketitle

\begin{abstract}
We provide a near-complete classification of the Lorentz spaces $\l_{\v}$
for which the sequence $\{S_{n}\}_{n\in \N}$ of partial Fourier
sums is almost everywhere convergent along lacunary subsequences.
Moreover, under mild assumptions on the fundamental function $\v$, we identify
$\l_{\v}:= L\log\log L\log\log\log\log L$ as the \emph{largest} Lorentz space
on which the lacunary Carleson operator is bounded as a map to $L^{1,\infty}$.
In particular, we disprove a conjecture stated by Konyagin in his 2006 ICM
address. Our proof relies on a newly introduced concept of a
``Cantor Multi-tower Embedding," a special geometric configuration of tiles
that can arise within the time-frequency tile decomposition of the Carleson
operator. This geometric structure plays an important role in the behavior
of Fourier series near $L^1$, being responsible for the unboundedness
of the weak-$L^1$ norm of a ``grand maximal counting function" associated
with the mass levels.
\end{abstract}
$\newline$

\section{\bf Introduction}

\subsection{Historical Background.} In this paper we address the problem of the pointwise convergence of the Fourier series along lacunary subsequences. Regarded
in a broader context, the problem of the pointwise convergence of the Fourier series has a rich history, tracing back to the cornerstone set by Fourier in his study on the heat propagation (\cite{Fou}). Since then, there has been a series of major advancements, from which we will mention only those closest to our topic of investigation: in 1873 du Bois-Reymond (\cite{Bois}) offers an example of a continuous function whose Fourier series diverges on the set of rational points. This surprising result stimulated the search for new grounds upon which one could reformulate the question of the pointwise convergence
for larger classes of functions by focusing only on the ``almost everywhere" behavior of the series and thus allowing for pathologies on ``negligible" sets (of ``measure $0$"). The appropriate setting was developed by Lebsesgue in his theory of integration (\cite{Leb}). Within this new framework,  Luzin (\cite{Lu}) conjectured that if a function $f$ is square integrable then its Fourier series converges to $f$ Lebesgue-almost everywhere. In 1922, Kolmogorov (\cite{Kol1}) constructed an example of an $L^1$-integrable function whose Fourier series diverges almost everywhere, suggesting that Luzin's conjecture may be false.
 However, after several decades of misconceptions, in 1966, the breakthrough of L. Carleson (\cite{c1}) provides the positive answer to this conjecture. Shortly after that, Hunt (\cite{hu}), extended the techniques in \cite{c1}, showing that Carleson's result holds for any function $f\in L^p(\TT)$ as long as $1<p<\infty$.

 Once at this point, we should state that though not providing a new result, the second proof of the almost everywhere convergence of the Fourier series for $L^2$ functions offered by C. Fefferman, (\cite{f}), marked a fundamental advancement in understanding the topic described here. A third proof of Luzin's conjecture was given in 2000 by Lacey and Thiele (\cite{lt3}) using the tools they developed
 for addressing the boundedness of the bilinear Hilbert transform (\cite{lt1}, \cite{lt2}).

\subsection{Formulation of the main problem(s); Context.}  We start this section by formulating (at first in a looser language), one of the main open questions in the area of Fourier series:
$\newline$

\noindent\textbf{Main Question.} What is the behavior of the (almost everywhere) pointwise convergence of the
Fourier Series between the two known cases for the Lebesgue-scale spaces $L^p(\TT)$:

\begin{itemize}
\item $p=1$, divergence of the Fourier Series (Kolmogorov)
\item $p>1$, convergence of the Fourier Series (Carleson-Hunt) ?

\end{itemize}

In order to make this Main Question precise, let us first introduce the following:

\begin{d0}\label{Cprop0} Let $Y$ be a r.i. (quasi-)Banach space.\footnote{For basic definitions and concepts referring to the theory of rearrangement invariant Banach spaces, including Lorentz spaces, see the Appendix.} We say that
$Y$ is a $\mathcal{C}-$space iff $\:\exists\:\:C_0=C_0(Y)>0$ such that the Carleson operator defined\footnote{Here, depending on the context, we either identify the torus $\TT:\R/\Z$ with $[-\frac{1}{2},\,\frac{1}{2}]$ or with $[0,1]$.} by $$C:\,C^{\infty}(\TT)\,\mapsto\,L^{\infty}(\TT)$$ with
\beq\label{carl}
Cf(x):=\sup_{N\in\N}\left|\int_{\TT}e^{i\,2\pi\,N\,(x-y)}\,\cot(\pi\,(x-y))\,f(y)\,dy\right|\:,
\eeq
obeys the relation\footnote{Recall that the weak-$L^1$ quasinorm is given by $\|f\|_{1,\infty}:
=\sup_{\lambda>0} \lambda\,|\{x\,|\,|f(x)|>\lambda\}|$.}
\beq\label{carlb}
\|Cf\|_{1,\infty}\leq C_0\,\|f\|_{Y}\:\:\:\:\:\forall\:\:f\in Y\:.
\eeq
\end{d0}
$\newline$
 With this definition, the Main Question can be reformulated as follows:
$\newline$

\noindent\textbf{Open Problem A.} 1) \textit{Give a satisfactory description of the Lorentz spaces $Y\subseteq L^1(\TT)$ that are also $\mathcal{C}-$spaces. If such exists, describe the maximal Lorentz $\mathcal{C}-$space $Y_0$.}

2) \textit{Let $Y$ be a r.i. (quasi-)Banach space. Provide necessary and sufficient conditions for $Y$ to be a $\mathcal{C}-$space.}
$\newline$

 The best known results in investigating the above problem are
 \begin{itemize}
 \item on the \textit{negative} side: Konyagin (\cite{koDivf}, \cite{koDivff}) proved that if $\phi(u)=o(u\sqrt{\frac{\log u}{\log\log u}})$ as $u\rightarrow \infty$ then the space $\phi(L)=\l_{\bar{\phi}}$ is not a $\mathcal{C}-$space, where here $\bar{\phi}(t):=\int_{0}^{t} s\,\phi(\frac{1}{s})\,ds$. Thus, there
     exists $f\in\phi(L)$ with
\beq\label{negl1}
\lim \sup_{m\rightarrow\infty} S_{m} f(x)=\infty\:\:\:\:\:\:\:\:\:\textrm{for all}\:\:\:x\in\TT\:,
\eeq
where here $S_{m}f$ stands for the $m^{th}$ partial Fourier sum associated with $f$.

\item on the \textit{positive} side: Antonov (\cite{An}) showed that \eqref{carlb} holds for the Lorentz space
$Y=L \log L \log\log\log L$; later Arias-de-Reyna (\cite{Ar}) showed that $Y$ can be enlarged to a rearrangement invariant quasi-Banach space, named $QA$, and strictly containing\footnote{For an interesting study on the properties of $QA$ and on the relationship(s) between Antonov and Arias-de-Reyna spaces, see \cite{MMR}. In the same paper, the authors prove that under suitable conditions on the function $\v$ the space $\l_{\v}=L\log L\log\log\log L$ is the largest Lorentz space contained in $QA$.} $L \log L \log\log\log L$.
\end{itemize}

We add here that the first results along these lines were obtained on the negative side by Chen (\cite{ch}), Prohorenko (\cite{P}) and K\"orner (\cite{Kor}), and on the positive side by Sj\"olin (\cite{sj3}) and F. Soria (\cite{So1}, \cite{So2}).

It is worth noting that all the progress mentioned above on the positive side involved tools from extrapolation theory. Recently, using methods that rely entirely on time-frequency arguments, the author was able to reprove all the positive results by a unified approach (\cite{lvCarl1}).

Now recall that both the maximal Hardy-Littlewood operator and the (maximal) Hilbert transform are operators that are bounded
from $L \log L$ to $L^1$. At a heuristic level, the Carleson operator may be thought as a superposition
between the maximal Hardy-Littlewood operator and modulated copies of the (maximal) Hilbert transform. Thus,
one is naturally led to the following
$\newline$

\noindent\textbf{Conjecture 1.} \textit{The Lorentz space $Y=L\log L$ is a $\mathcal{C}-$space.}
$\newline$

As a simplified model for better understanding the difficulties of the Open Problem A (and of the derived conjecture above)  one can formulate its \textit{lacunary} version. Recall here that a sequence $\{n_j\}_{j\in\N}\subset\N$ is called \textit{lacunary}
iff
\beq\label{lacseq}
\underline{\lim}_{j\:\rightarrow\:\infty}\: \frac{n_{j+1}}{n_j}>1\:.
\eeq

Now, by analogy with the previous situation, we first introduce the following

\begin{d0}\label{Cprop} Let $Z$ be a r.i. (quasi-)Banach space.

\noindent i)  Assume $\{n_j\}_{j\in\N}\subset\N$ is a lacunary sequence. We say that
$Z$ is a $\mathcal{C}_{L}^{\{n_j\}_j}-$space iff $\:\exists\:C_1=C_1(Z,\,\{n_j\}_j)>0$ such that the $\{n_j\}_{j\in\N}$ - lacunary Carleson operator defined by
$$C_{lac}^{\{n_j\}_j}:\,C^{\infty}(\TT)\,\mapsto\,L^{\infty}(\TT)$$ with
\beq\label{carlac}
C_{lac}^{\{n_j\}_j}f(x):=\sup_{j\in\N}\left|\int_{\TT}e^{i\,2\pi\,n_j\,(x-y)}\,\cot(\pi\,(x-y))\,f(y)\,dy\right|\:,
\eeq
obeys the relation
\beq\label{carlblac}
\|C_{lac}^{\{n_j\}_j}f\|_{1,\infty}\leq C_1\,\|f\|_{Z}\:\:\:\:\:\forall\:\:f\in Z\:.
\eeq

\noindent ii) We say that $Z$ is a $\mathcal{C}_L-$space iff $Z$ is a $\mathcal{C}_{L}^{\{n_j\}_j}-$space \textbf{for any} lacunary sequence $\{n_j\}_{j\in\N}$.

Moreover, trough out the paper, if $Z$ is a $\mathcal{C}_L-$space, we will (often) express this as\footnote{Given $A,\,B>0$, throughout the paper,  we will use the notations $A\lesssim B$ and $B\gtrsim A$ to specify that there exists $C>0$ such that $A\leq C\,B$ and $B\leq C\, A$ respectively.}

\beq\label{carlblacuni}
\|C_{lac}f\|_{1,\infty}\lesssim \|f\|_{Z}\:\:\:\:\:\forall\:\:f\in Z\:,
\eeq
where here $C_{lac}$ stands for ``the generic" lacunary Carleson operator. \footnote{In \eqref{carlblacuni}, the implicit constant is allowed to depend on the specific choice of the lacunary sequence and on the space $Z$ but not on the function $f\in Z$.}
\end{d0}

We can now formulate the analogue of Open Problem A:
$\newline$

\noindent\textbf{Open Problem B.}
1) \textit{Give a satisfactory description of the Lorentz spaces $Z$ that are also $\mathcal{C}_L-$spaces. If such exists, describe the maximal Lorentz $\mathcal{C}_L-$space $Z_0$. \footnote{One can formulate a more specific question by prescribing a lacunary sequence $\{n_j\}_{j\in\N}$ and asking for a satisfactory description of the Lorentz spaces $Z$ that are also $\mathcal{C}_L^{\{n_j\}_j}-$spaces.}}

2) \textit{Let $Z$ be a r.i. (quasi-)Banach space. Provide necessary and sufficient conditions for $Z$ to be a $\mathcal{C}_L-$space.}
$\newline$

Initial progress on this problem (in a more general context, \cite{Zyg}) has been made by Zygmund who showed that $Z=L \log L$ is a Lorentz $\mathcal{C}_L-$space. On the negative side, Konyagin (\cite{ko2}) proved that if $\phi:\,\R_{+}\,\rightarrow\,\R_{+}$ is an increasing
function with $\phi(0)=0$ and $\phi(u)=o(u\log\log u)$ as $u\,\rightarrow\,\infty$ then $\phi(L)=\l_{\bar{\phi}}$ is not a $\mathcal{C}_L-$space. This last result was reproved later in a slightly modified context by Antonov (\cite{An1}). \footnote{As an immediate application of the concepts developed in this paper, one can obtain a simplified proof of the results in \cite{ko2} and  \cite{An1} - see Remarks section.}

In his invited talk at the 2006 International Congress of Mathematicians in Madrid, Konyagin stated the following

$\newline$

\noindent\textbf{Conjecture 2.}(Konyagin, \cite{ko1}) \textit{The Lorentz space $L\log \log L$ is a $\mathcal{C}_L-$space.}
$\newline$

Once at this point, we should say that one can phrase an analogue of the above conjecture for the
Walsh-Fourier series (\textit{i.e.} Is it true that relation \eqref{carlblacuni} holds for $Z=L\log \log L$ and $C_{lac}$ replaced by the lacunary Walsh-Carleson operator ?). In this latter context, Do and Lacey (\cite{LaDo}) were the first to make progress
by showing that if one takes $Z=L\log\log L\log\log\log L$ then \eqref{carlblacuni} holds for the Walsh form of the lacunary Carleson operator. Their proof relies on a projection argument which is not transferable to the Fourier series case.
\footnote{In the same paper, using previous results from extrapolation theory, the authors proved that the Walsh form of \eqref{carlblacuni} holds for a slightly larger quasi-Banach space $Z=Q_D$. This last space turns out to be isomorphic with the space $\W$ introduced in \cite{lvKony1}, though we have designed $\W$ by different means independent of extrapolation theory.}

In \cite{lvKony1}, we were able to prove the following
$\newline$

\begin{t1}\label{w} (\cite{lvKony1}) Let $\W$ be the quasi-Banach space defined by\footnote{Throughout the paper we will use the following convention: $\log k$ stands for $\log_{2} k$.}:
$$\W:=\{f:\:\TT\mapsto C\,|\,f\:\textrm{measurable},\:\|f\|_{\W}<\infty \}\:,$$
where
$$\|f\|_{\W}:=\inf\left\{\sum_{j=1}^{\infty}(1+\log j)\|f_j\|_1\,\log\log\frac{4\,\|f_j\|_{\infty}}{\|f_j\|_1}\:\:\bigg|\:\:
\begin{array}{cl}
f=\sum_{j=1}^{\infty}f_j,\:\\
\sum_{j=1}^{\infty}|f_j|<\infty\:\textrm{a.e.}\\
f_j\in L^{\infty}(\TT)
\end{array}
\right\}\;.$$
Then \beq\label{K1}
\|\,C_{lac}(f)\,\|_{1,\infty}\lesssim \|f\|_{\W}\:.
\eeq
We thus have that $Z=\W$ is a $\mathcal{C}_L-$space.  Moreover, $\W$ contains the space
$L\log\log L\log\log\log L$.
\end{t1}

Taking as a black-box Theorem \ref{w} above, Di Plinio, (\cite{dp}), proved that $L\log\log L\log\log\log\log L$ is a $\mathcal{C}_L-$space. Indeed, relying entirely\footnote{The difficult combinatorial and time-frequency techniques are nedeed precisely in order to show that $\W$ is a $\mathcal{C}_L-$space. Our present paper shows that Theorem \ref{w} can not be essentially improved.} on standard extrapolation techniques, he showed that
$$Z':=L\log\log L\log\log\log\log L\subseteq \W\,,$$
which applying now \eqref{K1} immediately implies that $Z'$ is a $\mathcal{C}_L-$space.

\subsection{Main results.} In this section we present the main results of this paper. These results are based on a newly introduced concept called  \textsf{``Cantor multi-tower embedding"} (\textbf{CME}) whose nature will be detailed in the next subsection. With these being said, we state the following

$\newline$
\noindent\textbf{Main Theorem 1.}
$\newline$

\textsf{There exists} \textit{a lacunary sequence $\{n_j\}_j$ and there exists a sequence of (positive) functions $\{f_k\}_{k\in \N}$ such that the following hold:
\begin{itemize}
\item each $f_k\in L^\infty(\TT)$ with
\beq\label{normb}
\|f_k\|_{L\,\log\log L}\approx 1\,;
\eeq
\item there exists
 \beq\label{normub}
\lim_{k\rightarrow \infty}\|f_k\|_{L\,\log\log L\,\log\log\log\log L}=\infty\,;
\eeq
\item there exists $C>0$ absolute constant such that for any $k\in\N$
\beq\label{blowup}
\|C_{lac}^{\{n_j\}_j}(f_k)\|_{L^{1,\infty}}\geq C\,\|f_k\|_{L\,\log\log L\,\log\log\log\log L}\,.
\eeq
\end{itemize}}
$\newline$
Next result states that the conclusion of above theorem remains true for \textit{any} lacunary sequence $\{n_j\}_j$. More precisely one has

\begin{t1}\label{ww}
Given \textsf{any} lacunary sequence $\{n_j\}_j$ there exists a sequence of (positive) functions $\{f_k\}_{k\in \N}$ such that
\eqref{normb},  \eqref{normub} and \eqref{blowup} continue to hold.
\end{t1}
$\newline$

\noindent\textbf{Main Theorem 2.}
$\newline$

\noindent 1) \textit{Define $\v_0:\:[0,1]\,\rightarrow\,\R_{+}$ as $\v_0(s):=s\,\log\log\frac{17}{s}\,\log\log\log\log\frac{17}{s}$.}

\textit{Let now $\v:\:[0,1]\,\rightarrow\,\R_{+}$ be a non-decreasing concave function with $\v(0)=0$. Then we have:}

i) \textit{If $\:\:\:\underline{\lim}_{{s\rightarrow 0}\atop{s>0}}\:\frac{\v(s)}{\v_0(s)}>0\:\:\:$ then the Lorentz space $\l_{\v}$ is
a $\mathcal{C}_L-$space.}

ii) \textit{If $\:\:\:\overline{\lim}_{{s\rightarrow 0}\atop{s>0}}\:\frac{\v(s)}{\v_0(s)}=0\:\:\:$ then the Lorentz space $\l_{\v}$ is
not a $\mathcal{C}_L-$space.}

iii) \textit{If $\:\:\:\underline{\lim}_{{s\rightarrow 0}\atop{s>0}}\:\frac{\v(s)}{\v_0(s)}=0<\overline{\lim}_{{s\rightarrow 0}\atop{s>0}}\:\frac{\v(s)}{\v_0(s)}\:\:\:$ then both scenarios are possible. More precisely, one can choose a $\v$ such that
$\l_{\v}$ is a $\mathcal{C}_L-$space while for another proper choice of $\v$ one has that $\l_{\v}$ is not a $\mathcal{C}_L-$space.}
$\newline$

\noindent2) \textit{Let $\v:\:[0,1]\,\rightarrow\,\R_{+}$ be a quasi-concave function. Consider the following statements:}

(A) \textit{The function $\v$ obeys}
\beq\label{growth}
\int_{0}^{1}\frac{-s\,\v_0''(s)}{\v(s)}\,ds\approx\int_{0}^{1}\frac{\v_0(s)}{\v(s)}\,\frac{ds}{s\,\log\frac{4}{s}\,\log\log\frac{4}{s}} <\infty\:.
\eeq

(B)  \textit{Any r.i. Banach space $X$ with fundamental function $\v_{X}=\v$ is a $\mathcal{C}_L-$space.}

\textit{Then the following are true:}

i) (A)  \textit{implies} (B);

ii) (B) \textit{implies} $\underline{\lim}_{{s\rightarrow 0}\atop{s>0}}\:\frac{\v_0(s)}{\v(s)}=0$;

iii) \textit{if there exists $\ep>0$ such that $s\:\mapsto\:\frac{\v_0(s)}{\v(s)}$ is increasing on $(0,\ep)$ then}
$$(A)\:\textit{equivalent to}\:(B)\:.$$
$\newline$
\noindent\textbf{Remark.} In fact, one can derive Main Theorem 2 from Main Theorem 1 and Theorem \ref{GenTh}
displayed immediately below. However, we prefer to give special attentions to Main Theorem 1 and 2 since these statements
include the more conceptual nature of our results.
$\newline$

\begin{t1}\label{GenTh} Let $k\in\N$ and  $\{r_j\}_{1\leq j\leq k}$ be positive real numbers. For $1\leq j\leq k$ define $y_j=2^{-\log j\,2^{2^j}}$. Then, for any $x_j\in(y_{j+1},\,y_j]$, one can construct measurable sets $F_j\subset\TT$ such that:
\begin{itemize}
\item the sets $\{F_j\}_{j\leq k}$ are pairwise disjoint;
\item  $|F_j|=x_j$;
\item  there exist $C_1\geq C_2>0$ absolute constants such that the function
\beq\label{functk}
f_k:=\sum_{j=1}^k r_j\,\chi_{F_j}\,,
\eeq
obeys the estimate
\beq\label{W}
C_2\,\|f_k\|_{\mathcal{W}}\leq \|C_{lac}^{\{2^j\}_j} f_k\|_{1,\infty}\leq C_1\,\|f_k\|_{\mathcal{W}}\:.
\eeq
\end{itemize}
\end{t1}

$\newline$

\noindent\textbf{Consequences of Main Theorems 1 and 2.} From point 1) iii) in Main Theorem 2, (see also the corresponding proof) one can deduce that there exist (infinitely many) Lorentz $\mathcal{C}_L-$spaces $\l_{\v}$ such that
$$L\log\log L\log\log\log\log L\subsetneq \l_{\v}\subsetneq\W\:.$$
  While these $\l_{\v}$ spaces are non-canonical, their fundamental functions $\v$ still share at infinitely many space locations the same structure as that of $\v_0$. Thus, under suitable, mild conditions on $\v$, $\l_{\v_0}$ becomes the \textit{largest}  Lorentz $\mathcal{C}_L-$space, this being simply the content of the following:
$\newline$

\begin{c0}\label{Orlicz}[\textsf{Maximal characterization}]

Let $\v:\:[0,1]\,\rightarrow\,\R_{+}$ be a non-decreasing concave function with $\v(0)=0$.
Assume that there exists
\beq\label{las}
\lim_{{s\rightarrow 0}\atop{s>0}}\:\frac{\v(s)}{\v_0(s)}\in [0,\,\infty]\:.
\eeq
Then the largest Lorentz $\mathcal{C}_L$-space $\l_{\v}$ for which $\v$ obeys \eqref{las} is given by
$$Z_0=L\,\log\log L\,\log\log\log\log L\:.$$
\end{c0}

Taking in \eqref{las} the function $\v(s)=s\,\log\log \frac{4}{s}$, one further deduces

\begin{c0}\label{conjK}[\textsf{Resolution of Konyagin's conjecture}]

Conjecture 2 is false.
\end{c0}

Once at this point, we record this surprising at first glance

\begin{o0}\label{W} Define the following operators\footnote{For more on the definitions, notations and properties of the discrete Carleson and Walsh models see Section 12.}
\begin{itemize}
\item the $\{n_j\}_j$ - (lacunary) Lacey-Thiele \textbf{discretized Carleson periodic model} as
\beq\label{CMp1}
 \tilde{C}^{\{n_j\}_j}f(x):=\sup_{j\in\N}\left|\sum_{P\in\P^{0,0,0,+}} <f,\,\phi_{P_l}>\,\phi_{P_l}(x)\,\chi_{\o_{P_{u}}}(n_j)\,\right|\:;
\eeq

\item the $\{n_j\}_j$ - (lacunary) \textbf{discretized Walsh-Carleson operator} defined as
\beq\label{WalshC1}
\tilde{C}_{W}^{\{n_j\}_j}f(x)=\sup_{j\in\N}\left|\sum_{R\in\mathcal{R}} <f,w_{R_l}> \,w_{R_l}(x)\,\chi_{\o_{R_u}}(n_j)\right|\;.
\eeq

\item the $\{n_j\}_j$ - (lacunary) \textbf{Walsh-Carleson operator} defined as
\beq\label{wcarlw}
C_{W}^{\{n_j\}_j}f(x):=\sup_{j\in\N}\left|\sum_{k=0}^{n_j} <f,\,w_k>\,w_k(x)\right|\:,
\eeq
 where here $w_N(x)$ stands for the $N^{th}$ Walsh mode regarded as periodic function on $\R$.

\item the $\{n_j\}_j$ - (lacunary) \textbf{averaged Walsh-Carleson model} by
\beq\label{carlw}
C_{AW}^{\{n_j\}_j}f(x):=\sup_{j\in\N}\left|\int_{\TT} w_{n_j}(x)\,w_{n_j}(-y)\,\cot(\pi\,(x-y))\,f(y)\,dy\right|\:.
\eeq

\end{itemize}
(It is worth mentioning here that in \cite{T}, the author proved that, unlike the Fourier case, there is no distinction between the discretized and the standard (non-discretized) Walsh-Carleson operator, that is $C_{W}^{\{n_j\}_j}f=\tilde{C}_{W}^{\{n_j\}_j}f$.)
$\newline$

Now the following are true:

\begin{itemize}
\item Theorem \ref{ww} \textsf{holds} for the operator $C_{AW}^{\{n_j\}_j}$ (and obviously for $C^{\{n_j\}_j}$);

\item Theorem \ref{ww} \textsf{does not hold} for the operators $\tilde{C}^{\{n_j\}_j}$ and $C_{W}^{\{n_j\}_j}$.
\end{itemize}

All these facts will be discussed in great detail in Section 12. Notice that this is the first time when we are witnessing a sharp distinction between the behavior of the Carleson operator and that of the corresponding Lacey-Thiele discretizated Carleson model. This also provides a first instance when Fefferman's type discretization - \textit{which leaves the Carleson operator unchanged} - is a necessity and not a choice.
\end{o0}

\begin{o0}
The next corollary answers an \textbf{open question} related to the so called \textbf{Halo conjecture}\footnote{For more details on the connections between the Halo conjecture and the subject of the pointwise convergence of the Fourier Series the interested reader is referred to \cite{SS}, \cite{So1} and \cite{PH}.} (see \cite{PH}), regarding whether or not, given a sublinear, \textit{translation invariant} operator $T$, we must always have that the following are equivalent:
\begin{itemize}
\item  $T$ is of restricted weak type $(\l_{\v},\,L^1)$;
\item  $T$ is of weak type $(\l_{\v},\,L^1)$.
\end{itemize}
Here $\l_{\v}$ is some generic Lorentz space.
\end{o0}

\begin{c0}\label{Halo}[\textsf{Restricted weak-type does not imply weak-type}]

The $\{2^j\}_{j\in\N}$ - lacunary Carleson operator obeys the following:
\begin{itemize}
\item $C_{lac}^{\{2^j\}_{j}}$ is a sublinear, translation invariant operator;

\item (Theorem \ref{w} (\cite{lvKony1})) $C_{lac}^{\{2^j\}_{j}}$ is of restricted weak type $(L\log\log L, \,L^1)$; thus there exists an absolute constant $C>0$ such that
\beq\label{restrwtype}
\|C_{lac}^{\{2^j\}_{j}}(\chi_{E})\|_{1,\infty}\leq C\,|E|\,\log\log\frac{4}{|E|}\,,
\eeq
for any measurable set $E\subseteq\TT$;

\item (Main Theorem 1) $C_{lac}^{\{2^j\}_{j}}$ is not of weak type $(L\log\log L, \,L^1)$.

\end{itemize}
 \end{c0}

If one attempts to regard $L^{1,\infty}$ as a limiting space of the scale $\{L^{p,\infty}\}_{p>1}$ one obtains
the following:

\begin{c0}\label{Limextr}[\textsf{Limitations of extrapolation theory}]

The $\{2^j\}_{j\in\N}$ - lacunary Carleson operator obeys:
\begin{itemize}
\item  There exists $c>0$ such that for any $1<p\leq 2$
\beq\label{weaklp}
\|C_{lac}^{\{2^j\}_{j}} f\|_{p,\infty}\leq c \log (2+\frac{1}{p-1})\,\|f\|_{p}\:.
\eeq

\item There exists no $C>0$ such that
\beq\label{weakl1}
\|C_{lac}^{\{2^j\}_{j}} f\|_{1,\infty}\leq C \|f\|_{L\log\log L}\:\:\:\:\:\:\forall\:f\in L\log\log L\;.
\eeq
\end{itemize}
\end{c0}

The fact that \eqref{weaklp} holds can be easily derived from the proof of Theorem 1 (\cite{lvKony1}) as noticed in \cite{dpo}.

Now standard interpolation/extrapolation\footnote{See Section 1 in \cite{dpo} for details.} results show that if \eqref{weakl1} were to be true then \eqref{weaklp} would immediately follow. However
in the light of our Main Theorems 1 and 2 we notice that the reverse implication is false.

\begin{o0}\label{inter}
As a consequence of our last two corollaries we have the following important conclusion: \textsf{No general principle can be established\footnote{For large classes of operators that include the family of (maximal) operators associated with partial Fourier sums.} in terms of the \textit{equivalence} between the weak-$L^1$ type bounds and either the corresponding restricted weak-type $L^1$ bounds or weak-$L^p$ bounds ($p>1$). Moreover, extrapolation theory by itself is not suitable to provide sharp answers to the endpoint questions on the pointwise convergence of Fourier Series near $L^1$. In order to do so, one needs to take advantage of the special structure of the Carleson operator and hence to exploit the time-frequency analysis methods.}
\end{o0}

Finally, addressing a problem presented by the author in a previous work, we have:

\begin{c0}\label{revista}[\textsf{Lack of uniform control for Calderon-Zygmund tile partition}]

The open question raised by the author in \cite{lvCarl1} has a negative answer. More precisely, preserving the notations therein, there is no absolute constant $C>0$ such that, for $\a\in \N$, one has
$$\|T^{\P^{\a}}f\|_1\leq C\,\|f\|_1\:\:\:\:\:\:\:\:\:\forall\:f\in L^1(\TT)\,.$$
Moreover, there exists $f\in L^1(\TT)$ such that if one partitions $\P^{\a}:=\bigcup_{n=1}^N \P^{\a}_{n}$ with each $\P^{\a}_{n}$ having constant mass (i.e. $A(P)\approx 2^{-n}$ for any $P\in \P^{\a}_{n}$) one has
$$\|T^{\P^{\a}}f\|_{1,\infty}\gtrsim\, N\,\|f\|_1\:.$$
\end{c0}

\subsection{The fundamental idea.}  In this subsection we will describe at the philosophical level, a key aspect of the present work - that of introducing the concept of a \textsf{``Cantor multi-tower embedding"} (\textbf{CME}) \footnote{We warn the reader that these explanations can be truly understood only by experts within the time-frequency area since this newly defined concept goes very deeply at the heart of the time-frequency methods involved in analyzing the boundedness properties of the Carleson operator. Given this, a reader who is unfamiliar with these techniques might choose to skip for the moment being this subsection  and return to it only after being exposed gradually to the construction of our counter-example.}:
\begin{itemize}
\item \underline{\textit{What is it?}} It refers to a special geometric configuration of a set of tiles that, potentially, could be part of the time-frequency decomposition of the (lacunary) Carleson operator. The essence of this geometric configuration is that it is an \textit{extremizer} for the $L^{1,\infty}$-norm of a ``\textsf{grand maximal counting function}" (see \eqref{grandmaxcount} and \eqref{countj}), a newly defined object, that, as it turns out, plays a fundamental role in the behavior of the pointwise convergence of Fourier series near $L^1$.

    The existence of such tile configurations is a manifestation of the ``mass transference" phenomenon from the ``heavy" tiles ($P\in\P_n$ with $n\in\N$ close to $1$) towards ``light" tiles (i.e. those tiles $P\in\P_n$ for large $n\in\N$) that is capable of realizing a Cantor set structure for each of the sets $E(P)$ corresponding to a $P$ within the given tile configuration.

    Thus, in constructing a \textbf{CME}, a key role is played by the \textit{structure} of the sets $E(P)$ and not only by
    their relative \textit{size}.

\item \underline{\textit{Context within the literature.}} This particular configuration of tiles and the central role played by the corresponding grand maximal counting function are novel facts, which, after author's knowledge, do not have a direct correspondent within the previous time-frequency literature. However, the more elementary concept of a counting function has been used in many time-frequency papers and in the framework of the Carleson operator was first considered in \cite{f}.

    Regarded in a broader context, the idea of studying extreme geometric configurations along with their potential key role in deciding the answer to a (harmonic analysis) problem has been successfully applied in many instances. Two such classical examples are given by

    - the (un)boundedness properties of certain (sub-)linear operators,  e.g. Besicovitch/Kakeya sets relating the ball multiplier or Bochner-Riesz problems;

    - special topological/additive structure properties of sets,  e.g. Cantor sets.

\item \underline{\textit{What is its purpose?}} Based on the geometric location of the tiles within a \textbf{CME}, and in particular on the lacunary structure of the frequencies, we will first split the mass parameter $n$ into dyadic blocks. Then, for each such block, call it $B_j$, we will construct a corresponding set $F_j\subset\TT$ that realizes the alignment of the signum of all the components $\{T_P^{*} 1(x)|_{x\in F_j}\}_{P\in\F_j}$ where here $\F_j$ is the collection of all tiles $P$ inside the above mentioned \textbf{CME} that have the mass parameter $n\in B_j$. In this process we will make essential use of the fact that the exponentials $\{e^{i\,2^j\,2\,\pi\,\cdot}\}_{j\in\N}$ oscillate \textit{independently} in $[0,1]$ behaving similarly to a sequence of i.i.d. random variables. As a consequence of this property, we will thus be able to ``erase" the signum of the operators associated to various trees (of tiles), transforming the adjoint lacunary Carleson operator restricted to this tile configuration into a positive operator. Once at this point, taking the input function $f=\sum_j r_j \chi_{F_j}$ with $\{r_j\}$ arbitrary positive coefficients and $\{F_j\}$ as above, one concludes that $$\|C_{lac} f\|_{1,\infty}\approx\|f\|_{\mathcal{W}}\;.$$

\end{itemize}

The precise definition of \textbf{CME} is somewhat intricate at the technical and notational level and thus we will defer it for later (see Definition \ref{cme} below). For the time being, sacrificing a bit in the way of rigor, we state the following

\begin{t1}\label{Tcme} (\textsf{Existence})

The \textbf{CME} structures are compatible with the tile decomposition of the (lacunary) Carleson operator.
\end{t1}

More precisely, these structures can arise within the process of the time-frequency decomposition of the (lacunary) Carleson operator and are composed by tiles that, on the one hand, contain some prescribed ``amount" of the graph of the measurable function $N$ appearing in the linearization of the operator and, on the other hand, have some specific relative position one to the other. That is why the existence of such structures within the time-frequency decomposition process is non-trivial.

$\newline$
\subsection{Remarks.} 1) The present paper sheds new light on the topic of the pointwise convergence of Fourier series near $L^1$:
\begin{itemize}

\item Our Main Theorem 2 establishes near optimal necessary and sufficient conditions for a Lorentz space $\l_{\v}$ to be a $\mathcal{C}_L$-space. It also provides a very good description of when, given a quasi-convex function $\v$, \textit{any} r.i. Banach space $X$ with $\v_X=\v$ is a $\mathcal{C}_L$-space. Moreover, it \textit{essentially} states (see \eqref{W}) that the largest r.i. quasi-Banach space on which $C_{lac}$ is $L^{1,\infty}$ bounded is the space $\mathcal{W}$ introduced in \cite{lvKony1}.

In the literature regarding the (almost everywhere) pointwise convergence of Fourier Series, almost sharp results of this type constitute a novelty.

In a forthcoming paper (\cite{lvPCFS2}), we will unravel the subtle underlying principles that lie at the foundation of the present result and that lead naturally to the idea of considering the \textbf{CME} structure.

\item This second item, focuses on the \textit{method} of approaching the difficult problem of pointwise convergence of the Fourier series near $L^1$. Until very recently, all the progress made on the positive side of this topic was based on extrapolation theory.

    Using a different perspective - relying only on time-frequency reasonings - the author reproved (\cite{lvCarl1}) the best current positive results. The present paper goes beyond offering an alternate approach, as our results (Main Theorem 1 and 2, and Corollary 4) cannot be attained by pure extrapolation methods. Thus, in the context of the $L^1-$methods, this constitutes a first instance when the efficiency of time-frequency techniques is overtaking the canonical extrapolation approach used until now and serves for the idea advocated by the author that in order to make substantial progress on the problem of the convergence of Fourier series near $L^1$ one needs to leave the general extrapolation theory framework and make \textit{essential} use of the special \textit{structure} of the (lacunary) Carleson operator.

\item Thirdly, the spaces $Y=L\log L\log\log\log L$ and $Y=QA$ (i.e. the best known positive results for Open Problem A) viewed for a long time as mere byproducts of extrapolation techniques, are now revealed to be direct manifestations of the ``positive behavior" of the operators associated with generic \textbf{CME}. Indeed, as a consequence of the ideas introduced here\footnote{We will skip here detailed explanations as this goes beyond the purpose of the current paper; however, we will clarify in forthcoming work, including \cite{lvPCFS2}, many of the considerations discussed in the informal principle.} we have the following
$\newline$

\noindent\textbf{Informal principle}: \textit{The oscillation and mass transference from heavy to light tiles encapsulated in a generic \textbf{CME} structure represent the real enemy in advancing on Open Problem A. If one can reduce the behavior of the adjoint Carleson operator restricted to a general \textbf{CME} to the corresponding positive operator
(i.e. the absolute sum of the operators associated with the maximal, weight-uniform trees within it) then}, essentially, \textit{the largest Lorentz $\mathcal{C}$-space is precisely Antonov's space $L \log L\log\log\log L$. If on the contrary, considering the phase oscillation, one can remove (due to the extra-cancelation) the threat represented by the operators associated with these geometric structures, then Conjecture 1 can be answered affirmatively.}
$\newline$

In view of this heuristic, we can now summarize as follows:

In the \textit{lacunary} situation, based on the \textit{existence} of the \textbf{CME} configurations and on the oscillatory independence of the lacunary trigonometric system $\{e^{i\,2^j\,2\,\pi\,\cdot}\}_{j\in\N}$, one will be able to reduce the adjoint lacunary Carleson operator restricted to this special geometric configuration of tiles (and applied to a special input function) to the corresponding positive operator. Thus, as mentioned earlier, one concludes that $\W$ is essentially the largest $C_L-$space.

In the \textit{general} situation (i.e. that of the \textit{full} sequence of partial Fourier sums), while one can easily adapt the  extremal tile configuration to the new context, there is no analogue of the oscillatory independence of the lacunary trigonometric system. Thus the frequency locations of the tiles play now a determinative role in the boundedness properties of the Carleson operator: if one forms the positive counterpart\footnote{See the informal principle for the way in which this positive operator is defined.} of the adjoint of $C$ associated with a generic \textbf{CME} - call it $C^{*}_{+}$ - then one discovers that, essentially, the largest r.i. quasi-Banach space for which $C_{+}$ is $L^1-$weak bounded is given by $QA$ (notice the analogy with the lacunary case). Conversely, improving on Antonov's and Arias-de-Reyna's results would require precisely showing that there is some cancelation inside the operators associated with a \textbf{CME}. This is the key point where techniques from additive combinatorics might play an important role.

As a last remark, one may notice the following interesting analogy: in both the case of the Carleson operator and that of the Bilinear Hilbert Transform the current technology (producing the best up to date results) stops at the point where one needs to consider the signum/oscillation of some terms associated to particular structures:
\begin{itemize}
\item in the Carleson operator case - the current methods can't do better than providing bounds for the \textsf{positive} adjoint Carleson operator $C^{*}_{+}$ restricted to a generic \textbf{CME};

\item in the Bilinear Hilbert transform case (see e.g. \cite{lt2}) - the current methods can only deal with estimating
the absolute values of the elementary building blocks in the Gabor decomposition of the model operator.
\end{itemize}
This is the point where we believe that further progress (in either of the directions) require supplementary methods - very likely connected with the additive combinatoric structure of the frequencies of the trees in the time-frequency decomposition of the corresponding operators. In fact as a confirmation of the usefulness of additive combinatoric techniques in closely related time-frequency problems one can check paper \cite{Tt}.
\end{itemize}

2) Passing now to the negative results (i.e. finding ``the smallest" r.i. Banach space $X$ with $L\log L\subseteq X\subset L^1$ on which we have divergence of the Fourier series), it is likely that part of the ideas in our present paper will help improving the result(s) in \cite{koDivf} and \cite{koDivff}.

3) Finally, as briefly mentioned earlier, the geometry of the tile configuration introduced here is in fact the expression of the behavior of a specific ``grand maximal counting function" near $L^1$. Though we will not detail this subject here, it is worth saying that this function controls each of the counting functions of order $n,\: (n\in\N)$, i.e. those functions  that count the number of top maximal trees of mass $2^{-n}$ above each point $x\in[0,1]$. One should also add that the $BMO$ behavior of each of these counting functions of order $n$ played a fundamental role in removing the exceptional sets in the discretization of the Carleson operator. This last fact generated a first direct proof (i.e. without using interpolation), see \cite{lvPolynCar}, of the \textit{strong} $L^2$ boundedness of the (polynomial) Carleson operator (for an earlier approach on weak-$L^2$ bounds and strong $L^p$ bounds with $1<p<2$ see also \cite{lvQC}).
$\newline$

We plan to detail many of the above considerations regarding the pointwise convergence of the \textit{full} sequence of partial Fourier sums near $L^1$ in a subsequent paper.
$\newline$

\begin{o0}\label{He} In what follows we will build the sequence of steps that one needs to account for in constructing a sharp counterexample to the conjecture of Konyagin. In order to do so, we proceed as follows:
\begin{itemize}
\item Section 2 presents a very brief overview of the nature of the counterexample.

\item Section 3 reviews the discretization of the Carleson operator following Fefferman's approach (\cite{f}). As it turns out it is very important that this is an \textit{exact} discretization of the Carleson operator unlike the one provided by Lacey and Thiele in \cite{lt3}.

\item Section 4 introduces the main definitions required for our further reasonings; it is a technically involved section.

\item in Section 5 we present the main heuristic for our approach and also test the efficiency of the definitions and concepts introduced in Section 4 in order to treat a toy model of our problem that already strengthens the best results known to date; it is aimed to prepare the reader for the very technical sections to follow (especially the ones from 7 to 10);

\item Section 6 discusses a key concept introduced in the present paper - the grand maximal counting function.

\item Section 7 presents the generic construction of a Cantor Multi-tower Embedding (\textbf{CME}).

\item Section 8 explains in detail the construction of the input functions corresponding to the sets $\{F_j\}_{j=\frac{k}{2}+1}^{k}$.

\item Section 9 is meant to eliminate ``the back-ground noise" arising from the error terms; it can be skipped at the first reading of this paper.

\item Section 10 contains the proof of Main Theorem 1.

\item Section 11 presents the proof of Main Theorem 2 and can be read independently of all the other sections; it relies on extrapolations techniques.

\item Section 12 (roughly 13 pages) was added at the request of the referee and does not contribute ``per se" to either Main Theorems 1 or 2; it can be completely skipped at the first hand reading. It explains the sharp contrast between the behavior of the (lacunary) Carleson operator (using Fefferman's discretization) and the corresponding behavior of Lacey-Thiele discretized Carleson model and (lacunary) Walsh-Carleson operator respectively.

\item Section 13 addresses several final remarks.

\item finally in the Appendix we recall several standard facts about the theory of rearrangement invariant Banach spaces.
\end{itemize}

We encourage the reader to be patient with the sequence of technical definitions that will soon follow and, at the first glance through the paper, focus more on the main heuristics and ``big picture" information provided in Sections 2, 5 and 6.
\end{o0}

\section{Construction of the counterexample - an overview}
$\newline$

In this section we present the general strategy for proving Corollary 2.

We will show that for each $k= 2^{2^K}$ with $K\in\N$ large, there exists a function $f_k\in L^{\infty}$ with the following properties:
\begin{itemize}
\item $f_k$ is given by an expression of the form
\beq\label{fgenstruct}
f_k=\sum_{j=\frac{k}{2}+1}^{k} 2^{\log k\,2^{2^{j}}}\,\chi_{F_j}\,,
\eeq
where here $\chi_{F_j}$ designates the characteristic function of $F_j$ and each set $F_j$ has some prescribed properties that will be detailed shortly;

\item the $L\,\log\log L$ norm is under control:

$$\|f_k\|_{L\,\log\log L}\approx 1\,.$$

\item the $L^{1}-$weak norm of $C_{lac}^{\{2^j\}_{j}}(f_k)$ is large:

$$\|C_{lac}^{\{2^j\}_{j}}(f_k)\|_{L^{1,\infty}}\gtrsim \log k\,.$$

\end{itemize}

The construction of each $F_j$ requires some degree of technicalities and will be detailed later. As of now, we limit ourselves to only reveal the following properties:
\beq\label{propertFj}
\eeq
\begin{itemize}
\item $F_j\subseteq [0,1]$ measureable set;

\item $|F_j|\approx 2^{-\log k\,2^{2^{j}}}\,\cdot\,2^{-j}\,\cdot\frac{1}{k}$;

\item $F_j$ has a (finite) Cantor type structure.
\end{itemize}

We end this section by mentioning that the construction of the sets $\{F_j\}_j$ will be strongly dependent on the choice of the measurable function $N$ in the linearization of the $\{2^j\}_{j}$ - lacunary Carleson operator (see next section). Consequently, understanding/designing the structure of the set of tiles $\P$ appearing in the decomposition of the $\{2^j\}_{j}$ - lacunary Carleson operator $C_{lac}^{\{2^j\}_{j}}$ is a precondition for assigning the precise properties of each $F_j$.

\section{Discretization of our operator}

Let us first recall the main object of our study

\beq\label{carlac1}
C_{lac}^{\{2^j\}_{j}}f(x)=\sup_{j\in\N}\,|\int_{\TT} \frac{1}{x-y}\,e^{i\,2\pi\,2^j\,(x-y)}\,f(y)\,dy| \:.
\eeq

Applying Fefferman's discretization (\cite{f}), we follow the same steps as in \cite{lvKony1}:
\begin{itemize}
\item linearize our operator and write
$$Tf(x):=\int_{\TT} \frac{1}{x-y}\,e^{-i\,2\pi\,N(x)\,y}\,f(y)\,dy\,,$$
where here $N:\:\TT\rightarrow\:\{2^j\}_{j\in\N}\;$ measurable function. (Here, for technical reasons, we erase the term $N(x)\,x$ in the phase of the exponential, as later
in the proof this will simplify the structure of the adjoint operators $T^{*}$.)

\item use the dilation symmetry of the kernel and express
$$\frac{1}{y}=\sum_{k\geq 0} \psi_k(y)\:\:\:\:\:\:\:\:\:\forall\:\:0<|y|<1\:,$$
where  $\psi_k(y):=2^{k}\psi(2^{k}y)$ (with $k\in \N$) and $\psi$ an odd
$C^{\infty}$ function such that
\beq\label{psi}
\operatorname{supp}\:\psi\subseteq\left\{y\in \R\:|\:2<|y|<8\right\}\:.
\eeq

\item write\footnote{Throughout the paper we use the convention $0\in\N$. Thus $\N=\{0,\,1,\,2\ldots\}$.} $$Tf(x)=\sum_{k\in\N}\int_{\TT}e^{-i\,2\pi\,N(x)\,y}\,\psi_{k}(x-y)\,f(y)\,dy\:.$$

\item for each $k\in\N$, we partition the time-frequency plane in tiles (rectangles of area one) of the form $P=[\o,I]$ with $\o,\,I$ dyadic intervals (with respect to the canonical dyadic grids on $\R$ and respectively $[0,1]$) such that $|\o|=|I|^{-1}=2^{k}$. The set of all such tiles will be denoted by $\bar{\P}^{k}$. Further on, we let $\bar{\P}=\bigcup_{k\in\N} \bar{\P}^{k}$.

\item to each  $P=[\o,I]\in\bar{\P}$ we assign the set
$$E(P):=\{x\in I\,|\,N(x)\in\o\}\,,$$
that is responsible for the mass (or ``weight") of the tile - $\frac{|E(P)|}{|I|}$.
This mass concept will play a key role in partitioning the set $\bar{\P}$.

\item for $P=[\o,I]\in\bar{\P}^{k}$ with $k\in\N$ we define the operators
\beq\label{TP}
T_{P}f(x)=\left\{\int_{\TT}e^{-i\,2\pi\,N(x)\,y}\,\psi_{k}(x-y)\,f(y)\,dy\right\}\chi_{E(P)}(x)\:,
\eeq
and conclude that
\beq\label{discret}
Tf(x)=\sum_{P\in\bar{\P}} T_P f(x)\,.
\eeq
\end{itemize}

Notice that if we think of $N:\:\TT\rightarrow\:\{2^j\}_j\;$ as a predefined measurable function then the above decomposition is
\textit{independent} on the function $f$.
$\newline$

 \begin{o0}\label{LocTT} For $P=[\o_P,I_P]\in\bar{\P}$ let $c(I_P)$ be the center of the interval $I_P$ and define $I_{P^*}=[c(I_P)-\frac{17}{2}|I_P|,\,c(I_P)-\frac{3}{2}|I_P|]\cup[c(I_P)+\frac{3}{2}|I_P|,\,c(I_P)+\frac{17}{2}|I_P|]$.
 Based on \eqref{psi}, we then deduce that:
\beq\label{supportin}
\textrm{supp}\,T_{P}\subseteq I_P\,,
\eeq
while the adjoint operator of $T_P$ denoted with $T_{P}^{*}$ obeys\footnote{This is a direct consequence of \eqref{supportin} and of the fact that $\psi_{k}$ is compactly supported.}
\beq\label{support}
\textrm{supp}\,T_{P}^{*}\subseteq I_{P^*}\:.
\eeq

As a consequence, if $P_1,\,P_2\subset \P$ such that $I_{P_1}\subset I_{P_2}$ and $|I_{P_1}|\leq 2^{-10}\,|I_{P_2}|$, then
we have
\beq\label{supportint}
\textrm{supp}\,T_{P_1}^{*}\cap\textrm{supp}\,T_{P_2}^{*}=\emptyset\:.
\eeq
By a standard reasoning  we will be able to arrange that the following holds: if $P_1,\,P_2\in\bar{\P}$ and $|I_{P_1}|\not=|I_{P_2}|$ then $|I_{P_1}|\leq 2^{-10}\,|I_{P_2}|$ or $|I_{P_2}|\leq 2^{-10}\,|I_{P_1}|$.
Thus \eqref{supportint} is automatically guaranteed if $I_{P_1}\subsetneq I_{P_2}$.
\end{o0}

 We will make repeated use of this observation in our construction process.
$\newline$

\noindent\textbf{Notation.} Throughout the paper, if $I$ is a (dyadic) interval of center $c(I)$ and $d>0$ a positive constant, then $d\,I$ designates an interval having the same center as $I$ - namely $c(I)$ - and having the length $|d\,I|:=d\,|I|\:.$
Also, if $P=[\o,\, I]$ and $a>0$ then we define the tile-dilation $a\,P:=[a\,\o,\,I]$.

\section{Main Definitions and Preparatives.}

In this section we will introduce several of the basic concepts which will be used later in the proof. The first
three definitions were introduced in \cite{f}, while definitions \ref{stree} and \ref{infforest}  were first developed in \cite{lvPolynCar}.

\begin{d0}\label{mass}(\textsf{weighting the tiles})

We define the \textbf{mass} of $P=[\o,I]\in\bar{\P}$ as

\beq\label{v1} A(P):=\sup_{{P'=[\o',I']\in\:\P}\atop{I\subseteq
I'}}\frac{|E(P')|}{|I'|}\:\frac{1}{(1+\frac{\textrm{distance}(10\o,\,10\o')}{|\o|})^{N_0}} \eeq where here $N_0$ is a fixed large natural number and if $A,\,B\subseteq\R$ then we refer at $\textrm{distance}(A,\,B)=\inf_{{a\in A}\atop{b\in B}} |a-b|$.

We also refer to the \textbf{restricted mass} or r-mass of $P=[\o,I]\in\bar{\P}$ as
$$A_0(P):=\frac{|E(P)|}{|I_P|}\,.$$
\end{d0}

\begin{d0}\label{ord} (\textsf{ordering the tiles})

   Let $P_j=[\o_j,I_j]\in\bar{\P}$ with $j\in\left\{1,2\right\}$. We say that
 $P_1\leq P_2$       iff       $\:I_1\subseteq I_2$
and  $\o_1\supseteq\o_2$. We write $P_1< P_2$ if $P_1\leq P_2$ and $|I_1|<|I_2|$.
\end{d0}

Notice that $\leq$ defines a partial order relation on the set $\bar{\P}$.
$\newline$

\begin{o0}\label{Partfam} In what follows we will define various families of tiles with some prescribed analytic and geometric properties (relative to the mass of a tile and to the order relation $\leq$, respectively.) In order to do so,
we will introduce several refinements of the set $\bar{\P}$ keeping always in mind that the analytic and geometric properties that we will describe are strongly influenced by the key fact
 \beq\label{imageN}
 \textsf{Image} (N)\subseteq \{2^j\}_{j\in\N}\;.
 \eeq

Let us first define
$$\P(0):=\{P\in\bar{\P}\,|\,0\in 100\,\o_P\}\,,$$
and
$$\bar{\P}_0:=\{P\in\bar{\P}\setminus \P(0)\,|\,A_0(P)=0\}\;.$$

Next set $$\P:=\bar{\P}\setminus \left(\P(0)\cup\bar{\P}_0\right)\;.$$

For $n\in\N$, we further define
\beq\label{pn}
\P_n:=\{P\in\P\,|\,A(P)\in(2^{-n-1},\,2^{-n}]\}\:.
\eeq
From now on, we will say that a tile $P$ has weight $n$ iff $P\in\P_n$.

Later on, in our construction of a \textbf{CME}, will be useful to impose the following restriction on the measurable function $N$:
\beq\label{imN}
  \textsf{Image}(N)\subseteq\{2^{2^{2^{2^k}}+100 m}\}_{m\in \{0,\,\log k\,2^{2^k-1}\}}\:,
\eeq
where we recall that here $k\in\N$ is a fixed large parameter.

Also we will ask that each tile $P\in\P$ obeys
\beq\label{ap}
A(P)\geq 2^{-2^{2k}}\;.
\eeq
Consequently, we deduce that
\beq\label{famtil}
\P=\bigcup_{n\leq 2^{2k}} \P_n\,.
\eeq
In particular we will only work under the assumption that $\P$ is finite.
\end{o0}

\begin{d0}\label{tree} (\textsf{tree})

 We say that a set of tiles $\p\subset\P$ is a \textbf{tree} with \emph{top} $P_0\in\P$ if
the following conditions hold:
$\newline 1)\:\:\:\:\:\forall\:\:P\in\p\:\:\:\Rightarrow\:\:\:\:P\leq P_0 $
$\newline 2)\:\:\:\:\:$if $P_1,\:P_2\: \in\p$ and $P_1\leq P \leq P_2$ then $P\in\p\:.$
\end{d0}

\begin{d0}\label{stree} (\textsf{sparse tree})

 We say that a set of tiles $\p\subset\P$ is a \textbf{sparse tree} if
 $\p$ is a tree and for any $P\in\p$ we have
 \beq\label{cmstree}
 \sum_{{P'\in\p}\atop{I_{P'}\subseteq I_P}}|I_{P'}|\leq C\,|I_P|\:,
 \eeq
 where here $C>0$ is an absolute constant.
 \end{d0}

\begin{d0}\label{infforest} (\textsf{forest} - $L^{\infty}$ control over union of trees)

Fix $n\in\N$. We say that $\p\subseteq\P_n$ is an $L^\infty$-\textbf{forest} (of $n^{th}$-generation) if

i) $\p$ is a collection of separated trees, {\it i.e.}
$$\p=\bigcup_{j\in\N}\p_j\,$$
with each $\p_j$ a tree with top $P_j=[\o_j,I_j]$ and such that
 \beq\label{wellsep}
 \forall\:\:j'\not=j\:\:\&\:\:\forall\:\:P\in\p_j\:\:\:\:\:\:\:\:2P\nleq 10P_{j'}\:.
 \eeq

ii) the counting function
\beq\label{cfunct}
 \cN_{\p}(x):=\sum_{j}\chi_{I_j}(x)
 \eeq
 obeys the estimate $\|\cN_{\p}\|_{L^{\infty}}\lesssim 2^n$.

 Further, if $\p\subseteq\P_n$ only consists of sparse separated trees then
we refer at $\p$ as a \textbf{sparse $L^\infty$-forest}.
 \end{d0}

 \begin{o0}\label{Forest} In this paper we will choose to focuss on the decomposition of our set of tiles $\P$ into
$L^\infty$-forests with some prescribed properties (see below). For this reason, unlike the precedent tile
decompositions in \cite{lvPolynCar} and  \cite{lvCarl1} (there we have introduced the concept of a BMO-forest), we will refer to an $L^{\infty}$ forest as simply a \textit{forest}.
\end{o0}

\begin{d0}\label{genf} (\textsf{generalized forest} - \textbf{GF})

Let $r,\,n\in\N$ with $r\leq n$. We say that $\p\subseteq\P$ is a \textbf{generalized forest} of $(r,n)$-generation  if we can decompose
\beq\label{genfordec}
 \p=\bigcup_{j=r}^{n} \p[j]\,,
\eeq
such that the following hold:
\begin{itemize}

\item each $\p[j]$ is an ($L^\infty$-) forest of $j^{th}$ generation;

\item if
$$\p[j]=\bigcup_l \p_l[j]\,,$$
is the decomposition of $\p[j]$ into maximal separated trees then for any pair $(j,\,j')$ with $r\leq j<j'\leq n$, any $l$ and any $P\in\p_l[j]$ there exists $l'$ and $P'\in \p_{l'}[j']$ such that
$$P\leq P'\:.$$
\end{itemize}

 \end{d0}

\begin{d0}\label{sgenf} (\textsf{saturated generalized forest} - \textbf{SGF})

Let $r,\,n\in\N$ with $r\leq n$. We say that $\p\subseteq\P$ is a \textbf{saturated generalized forest} of $(r,n)$-generation  if the following are true
\begin{itemize}
\item $\p$ is a GF of $(r,n)$-generation;

\item if $\exists\:P\in\P_j$ s.t. $\exists\:P'\in\p[j]$ with $P>P'$ then $P\in\p[j]$;

\item if $\exists\:P\in\P_j$ s.t. $\exists\:P'\in\p[j]$ with $P<P'$ then $P\in\p[j]$.

\end{itemize}
 \end{d0}

 \begin{o0}\label{LinN} All the previous definitions make perfect sense in a general context, with no particular restriction imposed on the linearization function $N$. The structures introduced in the next two definitions though, are not present (in a non-trivial form) for an arbitrary choice of $N$. However, in the context of this paper, we will have the liberty of \textit{choosing} $N$, and thus allow their existence. The precise form of the definitions below is chosen for simplifying as much as possible the general tile-configuration of the counterexample. For other, more general purposes, the requirements in these definitions can be significantly relaxed.
 \end{o0}

\begin{d0}\label{usgenf} (\textsf{uniform saturated generalized forest} - \textbf{USGF})

Let $r,\,n\in\N$ with $r\leq n$. We say that $\p\subseteq\P$ is a \textbf{uniform saturated generalized forest} of $(r,n)$-generation if the following are true
\begin{itemize}
\item $\p$ is a SGF of $(r,n)$-generation;

\item for each $j\in\{r,\ldots,n\}$ there exists $C_j\in (0,1]$ such that
$$|I_{P_l^j}|=C_j\:\:\:\:\:\:\:\forall\:l\:,$$
where, with the notations above, $P_l^j$ stands for the top of $\p_l[j]$.
\end{itemize}
 \end{d0}

\begin{d0}\label{usgentopf} (\textsf{uniform saturated generalized top-forest} - \textbf{USGTF})

We say that $\p\subseteq\P$ is a \textbf{uniform saturated generalized top-forest} of $(r,n)$-generation if
$\p$ is a USGF of $(r,n)$-generation and each tree $\p_l[j]$ in the definition above consists of just a single tile - its top.
\end{d0}

In our paper we will only work with some \textit{special} type of USGTF's. We will describe their properties in what follows, but first we need some more notations.

If $J\subseteq[0,1]$ dyadic interval and $m\in\N$, then we define
\beq\label{defimj}
\I_{m}(J):=\{I\subseteq J\,|\,I\:\textrm{dyadic},\:|I|=|J|\,2^{-m}\,\}\:.
\eeq
From now on, we will always apply the following convention: if $\{I_s\}_s$ is a collection of disjoint dyadic space intervals, then the indexing $s$  reflects the relative position of these intervals in $[0,1]$ from left to right, i.e., if $s_1<s_2$ then $c(I_{s_1})<c(I_{s_2})$.

\begin{o0}\label{co}
In what follows we will use (for awhile) an alternative description of a tile $P=[\o, I]$, namely $P:=I\times \a$, where here we write $\o=[l(\o),\,r(\o))$ and set $\a:=l(\o)$. This is justified since knowing $I$ and $l(\o)$ completely determines $P$ (recall that the area of the rectangle determined by $P$ is always assumed to be one). Notice that due to the definition of $C_{lac}^{\{2^j\}_j}$ we can always assume that $l(\o)\in\{2^j\}_{j\in\N}$.
\end{o0}

\begin{o0}\label{NUSGT} From now on, whenever we refer to a family $\F\subset\P$ as a USGTF (of $(r,n)$-generation) we will specify three sets of parameters, that, following the algorithm described below, will completely determine $\F$:

$\newline$
\noindent \textsf{The three sets of parameters}
$\newline$
\beq\label{param}
\eeq
\begin{itemize}
\item the collection of disjoint dyadic same-length space intervals
$$ITop(\F):=\{I_j\}_{j}\,.$$

\item the collection of \textit{distinct} frequencies (arranged in an increasing order)
$$\a(\F):=\{\a_u\}_{u=1}^{2^{n-1}}\,.$$

\item the collection of disjoint dyadic (same-length) space intervals
$$IBtm(\F):=\bigcup_{j} \I_{n-r} (I_j)\,.$$
\end{itemize}

$\newline$
\noindent \textsf{The algorithm (that completely  determines $\F$)}
$\newline$

We \textit{define} $\F$ the USGTF (of $(r,n)$-generation) obeying:
\begin{itemize}
\item the collection of tiles of weight $n$ inside $\F$ denoted with $\F[n]$ is given by:
 $$\F[n]=\{I\times \a\,|\,I\in ITop(\F)\:\:\&\:\:\a\in \a(\F)\}\:.$$
Also each $P\in \F[n]$ has the property that
\beq\label{weightpn}
A(P)=A_0(P)=2^{-n}\;.
\eeq
\item the collection of tiles of weight $r$ inside $\F$, denoted with $\F[r]$, is given by
$$\F[r]:=\bigcup_{j}\,\bigcup_{l=1}^{2^{n-r}}\,\{I_l\times \a\,|\,I_l\in\I_{n-r}(I_j) \:\:\&\:\:\a\in \{\a_{(l-1)\,2^{r-1}+1},\ldots,\,\a_{l\,2^{r-1}}\}\,\}\:.$$
As before, we require that $P\in\F[r]$ implies
\beq\label{weightpr}
A(P)=A_0(P)=2^{-r}\;.
\eeq
Thus, if $P\in \F[n]$, we impose that
\beq\label{keycompress0}
\textrm{there exists a unique}\:\:P'\in \F[r]\:\:\textrm{such that}\:\:P'<P\,,
\eeq
and on top of this we require
\beq\label{keycompress}
E(P)=E(P')\;.
\eeq
\end{itemize}

Notice that due to the above requirements, with a particular emphasis on \eqref{keycompress0} and \eqref{keycompress}, we have
that all the tiles in the intermediate families $\{\F[l]\}_{r<l<n}$ are now completely determined.

Thus, we have indeed that $\F$ is completely determined by the three sets of parameters defined in \eqref{param} once we make the convention that we will always run the above algorithm.
\end{o0}

\noindent{\textbf{Notation.}} Let $\A=\{\A_j\}_j,\,\B=\{\B_k\}_k$ be two collections of disjoint dyadic intervals. We write
$$\A\prec\B\,$$
iff
$$\forall\:j\:\:\:\exists\:k\:\:\textrm{s.t.}\:\:\:\A_j\subset \B_k\:.$$
Also we refer at $\tilde{\A}$ as
$$\tilde{\A}:=\bigcup_{j} \A_j\:.$$
If $\A_j^{lt}$ designates the left child of the interval $\A_j$ then we set
$$\A^{lt}:=\{\A_j^{lt}\}_j\:.$$
The same convention applies to $\A^{rt}$ (i.e. the collection of the right children of the intervals in $\A$.)

Recalling the definition of $\I_m(J)$ in \eqref{defimj}, for $\A$ as before we set
$$\I_m(\A)=\bigcup_{j}\I_m(\A_j)\:.$$

Also, for $J$ a dyadic interval, if $\I_m(J)=\{I_l\}_l$ then $\I_{m}^{lt}(J):=\{I_l^{lt}\}_l$ and $\I_{m}^{rt}(J):=\{I_l^{rt}\}_l$. Further on, $\I_m^{lt}(\A)=\bigcup_{j}\I_m^{lt}(\A_j)$ and similarly for $\I_m^{rt}(\A)$.

If $\A=IBtm(\F)$ for $\F$ a USGTF, then we set $\tilde{I}Btm(\F)=\tilde{\A}$, $I^{rt}Btm(\F)=\A^{rt}$ and
$I^{lt}Btm(\F)=\A^{lt}$. With the obvious changes, same applies to the case  $\A=ITop(\F)$.
 $\newline$

\begin{d0}\label{tow} (\textsf{tower})

We say that $\p\subseteq\P$ is a \textbf{tower} of $(r,n)$-generation if there exists $m\in\N,\,m\geq 1$ such that
$$\p=\bigcup_{l=1}^{m}\p_l\,,$$
and
\begin{itemize}
\item each $\p_l$ is a USGTF of $(r,n)$-generation;

\item $ITop(\p_{l+1})\prec IBtm(\p_{l})$ for any $l\in\{1,\ldots,\,m-1\}$;

\item from the two items above, we infer $\forall\:l\not=l'$ and $\forall\:P\in \p_{l}$ and $\forall\:P'\in \p_{l'}$ one has
\beq\label{incomp}
P\nleq P'\:\:\:\:\:\:\:\textrm{and}\:\:\:\:\:\:\:\:P'\nleq P\:.
\eeq
\end{itemize}

In particular we have\footnote{If \eqref{nocomfreq} were not true then condition
$\exists\:a\in \a(\p_l)\cap\a(\p_{l'})$ would imply that there exists $P\in \p_{l}$ and $P'\in \p_{l'}$ with $\a(P)=\a(P')=a$;
 this together with the second item in Definition \ref{tow} would imply $P\leq P'$ or $P'\leq P$ thus contradicting \eqref{incomp}.}
\beq\label{nocomfreq}
\forall\:l\not=l'\:\:\:\Rightarrow\:\:\:\a(\p_l)\cap\a(\p_{l'})=\emptyset\,.
\eeq

The number of USGTF's is called the \textit{height} of the tower while  $\tilde{I}Top(\p_{1})$ stands for its \textit{basis}. Thus, for $\p$ as above, we have
 $$Height(\p)=m\:\:\:\:\:\textrm{and}\:\:\:\:\:Basis(\p)=\tilde{I}Top(\p_{1})\;.$$
\end{d0}

\begin{d0}\label{tow} (\textsf{multi-tower})

We say that $\M\subseteq\P$ is a \textbf{multi-tower} of $(r,n)$-generation if one can decompose it as
$$\M=\bigcup_{l}\M_l\,,$$
such that
\begin{itemize}
\item each $M_l$ is a tower of $(r,n)$-generation;

\item $Basis(M_l)\cap Basis(M_{l'})=\emptyset$ for any $l\not=l'$.
\end{itemize}
\end{d0}

\begin{d0}\label{inclusion} (\textsf{multi-tower embedding})
If $\F^1,\,\F^2$ are two (multi)towers with $\F^j$ of generation $(r_j,\,n_j)$, we say that $\F^1$ \textbf{embeds} into $\F_2$ and write
$$\F^1\sqsubset\F^2\,,$$
iff
$$n_1\leq r_2\,,$$
and
$$\forall\:\:P\in \F^1[n_1]\:\:\:\exists\:\:P'\in \F^2[r_2]\:\:\:\textrm{such that}\:\:\:P\leq P'\:.$$
In particular, if $\F^1,\,\F^2$ USGTF's, we must have
$$ITop(\F_1)\prec IBtm(\F_2)\:\:\:\:\textrm{and}\:\:\:\a(\F_1)\subset\a(\F_2)\:.$$
\end{d0}

This finishes the preparatives required for presenting the main components of our proof.
$\newline$

\section{Heuristics and a warm-up example}
$\newline$

In order to smoothen the transition between two important technical sections of our paper and to help clarifying the ``big picture" in our reasonings, we start with the following
$\newline$

\noindent\textbf{Main Heuristic}. \textit{Our aim is to design some special function $N$ that will give rise to a family of embedded multi-towers, i.e., a \textbf{multi-tower chain} with respect to embedding relation ``$\sqsubset$". This chain will loosely have the form}
\beq\label{chain}
\F=\bigcup_{\frac{k}{2}<j\leq k} \F_j\,,
\eeq
such that
\begin{itemize}
\item each $\F_j$ a multi-tower of generation $(2^{j-1},\,2^j)$

\item $\F_j\sqsubset\F_{j+1}$.
\end{itemize}

 \textit{At the informal level, a (lacunary) \textbf{CME} will be a chain that maximizes the $L^{1,\infty}$-norm of a \textbf{grand maximal counting function} - a notion that will be our main focus in the section to follow.
\newline
As a consequence of this requirement, for $\F$ a \textbf{CME}, we have that for a generic tile $P\in\F$, the set $E(P)$ has a Cantor-type distribution inside $I_P$.}
$\newline$

Next, we would like to motivate the necessity of considering the \textbf{CME} concept and why we were required to develop the notions of multi-tower and chain of multi-towers in the previous section. For this, we will first discuss a simpler toy model, naturally developing from our introductory discussion in Section 1:

Indeed, as mentioned in the Introduction, from \cite{ko2} (see also \cite{An1}) we know that if $\phi:\,\R_{+}\,\rightarrow\,\R_{+}$ is increasing with $\phi(0)=0$ and $\phi(u)=o(u\log\log u)$ as $u\,\rightarrow\,\infty$ then $\phi(L)=\l_{\bar{\phi}}$ is not a $\mathcal{C}_L-$space. This result can now be easily deduced from the the following stronger claim:
$\newline$

\begin{p1}\label{S}. There exists $C>0$ absolute constant and a sequence of measurable sets $\{F_k\}_{k\in\N}$ with the following properties:
\begin{itemize}
\item $F_k\subseteq \TT$ with $|F_k|\rightarrow^{k\rightarrow\infty} 0$;

\item for any $k\in\N$ one has
\beq\label{claccc}
\|C_{lac}^{\{2^j\}_j}(\chi_{F_k})\|_{1,\infty}\geq C\,|F_k|\,\log\log \frac{4}{|F_k|}\:.\eeq
\end{itemize}
\end{p1}
$\newline$

The \textit{idea} of the proof relies on the newly introduced concept of \textsf{tower}: we let $\F_k$ be a tower of height $1$ that is a single USGTF! More precisely, using the same language as the one introduced in Observation \ref{NUSGT}, we define the collection of tiles \textsf{$\F_k$ be the USGTF of $(0, 2^{k})$-generation} given by the following characteristics:

\begin{itemize}
\item the collection of disjoint dyadic same-length space intervals
$$ITop(\F_k):=\{[0,1]\}\,.$$

\item the collection of \textit{distinct} frequencies\footnote{One can notice that here we made a minor modification of the USGTF model described in Observation \ref{NUSGT} by requiring that the stack of tiles having the mass $2^k$ has the height $2^{2^k}$ instead of $2^{2^k-1}$.}
$$\a(\F_k):=\{2^{2^{2^{2^k}}+100 m}\}_{m\in \{1,\ldots, 2^{2^k}\}}\,.$$

\item the collection of disjoint dyadic (same-length) space intervals
$$IBtm(\F):=\bigcup_{j} \I_{2^k} ([0,1])\,.$$
\end{itemize}

Next, one can construct a set $F_k$ with the following properties:
\begin{itemize}
\item the \textit{size} given by
$$|F_k|\approx 2^{-2^{2^k}}\,,$$
\item the \textit{structure} of $F_k$ is chosen such that\footnote{The mechanism of realizing the size and structure conditions of $F_k$ is described in the full generality (i.e. multi-tower case) in Section 8.1. and thus it will not be detailed here.}
$$\int\,\textrm{Re}\left(\chi_{F_k}\,T_{P}^{*}(1)\right)(\cdot)\approx\int\,|\chi_{F_k}\,T_{P}^{*}(1)(\cdot)|\:,$$
holds for all $P\in \F_{k}$.
\end{itemize}

Then, recalling that $\F_k[n]$ stands for the tiles in $\F_k$ of weight $n$, we conclude that

\beq\label{clac}
\|C_{lac}^{\{2^j\}_j}(\chi_{F_k})\|_{1,\infty}\gtrsim \sum_{n=1}^{2^k} \sum_{P\in \F_k[n]} 2^{-n}\,\frac{|I_P\cap F_k|}{|I_P|}\,|I_P|\geq 2^{k}\,|F_k|\:,
\eeq
thus proving our proposition.
$\newline$

\begin{o0} We stress here that for the proof of the above proposition there was no need to consider a chain of towers since we were not aiming for the extra $\log\log\log\log$ term in \eqref{claccc}; equivalently saying, our reasoning involved a single characteristic function of a set instead of an input function $f_k$ expressed as in \eqref{fgenstruct}, as a linear combination of $\approx\,k$ characteristic functions of sets. It is thus natural, that once we move our attention towards proving our Main Theorem 1, we need to consider the more involved concept of multi-tower and finally that of \textbf{CME}.
\end{o0}
$\newline$

\section{The grand maximal counting function}

As announced, in this section we will elaborate and motivate on the introduction of the new concept of \textbf{grand maximal counting function}, which is defined as

\beq\label{grandmaxcount}
\n:=\sup_{j}\n_j\,,
\eeq
where, recalling \eqref{pn}, we set
\beq\label{countj}
\n_j:=\frac{1}{2^{j-1}}\sum_{n=2^{j-1}+1}^{2^j} \frac{1}{2^{n-1}}\sum_{P\in\P_n^{max}}\chi_{I_P}\;,
\eeq
with $\P_n^{max}$ here designating the maximal tiles\footnote{Relative to $``\leq"$.} $P$ inside $\P_n$ such that $A(P)>2^{-n-1}$.

\begin{o0}\label{Countfunctmotiv}
The motivation for defining the counting functions $\n_j$ and $\n$ originates in \cite{lvPolynCar} where the author used
a complex greedy algorithm (involving more basic counting functions) in order to remove the exceptional sets arising in the
time-frequency discretization and provide direct strong $(2,2)$ bounds for the (standard) Carleson operator. Notice that $\n_j$
roughly controls the average spacial density location\footnote{That is the number of maximal trees siting above a point $x$ as $x$ runs through the interval $[0,1]$.} of the maximal trees of mass parameter $n\approx2^j$, while $\n$ pics the worst (largest) such density location among all the possible dyadic mass scales. The normalization factor
$\frac{1}{2^j}$ in \eqref{countj} preserves a uniform control of the $BMO$ norm of $\n_j$, that is one has $\|\n_j\|_{BMO(\TT)}\leq 10$ for any $j\in\N$.
\end{o0}

We move now to a further elaboration of the Main Heuristic whose key message is:  ``a \textbf{CME} is a chain that maximizes the $L^{1,\infty}$-norm of the grand maximal counting function".

\begin{o0}\label{Countfunct} Define for the beginning the $k-$truncated grand maximal counting function as
\beq\label{grandmaxcount}
\n^{[k]}:=\sup_{1\leq j\leq k}\n_j\,
\eeq
and notice that under the assumptions \eqref{imN}-\eqref{famtil} we trivially have that
\beq\label{equa}
\n^{[2k]}=\n\,.
\eeq
Next, as a consequence of the work in \cite{lvPolynCar}, we have that for each $j\in\N$, $(j\leq 2k)$
\beq\label{fa1}
\|\n_j\|_{BMO}\lesssim 1\:.
\eeq
This, together with standard John-Nirenberg inequality gives
\beq\label{fa2}
\|\n^{[2k]}\|_{1,\infty}\leq \|\n^{[2k]}\|_1\lesssim \log k\:.
\eeq
The crux of our main results in this paper is based on the fact that one can construct special configurations inside the family of tiles $\P$ corresponding to chains as described in \eqref{chain}, such that the inequality \eqref{fa2} can be reversed. In these instances, one can thus show that the following holds
\beq\label{fa3}
\|\n^{[2k]}\|_{1,\infty}\approx \log k\:.
\eeq
\end{o0}

We start now a more detailed analysis of the properties of the grand maximal counting function $\n$ by first presenting a short proof of \eqref{fa2}.

For this we start by defining
$$\overline{\n}_{n}:=\frac{1}{2^{n-1}}\sum_{P\in\P_n^{max}}\chi_{I_P}\,,$$
and notice that
\beq\label{functv}
\n_j:=\frac{1}{2^{j-1}}\sum_{n=2^{j-1}+1}^{2^j} \overline{\n}_{n}\:.
\eeq
Now, following the reasonings described in \cite{lvPolynCar}, we have that for every $n\in\N$ the function $\bar{\n}_{n}$ belongs to $BMO(\TT)$ with
$$\|\overline{\n}_{n}\|_{BMO(\TT)}\leq 10\:.$$
Now, since $\n_j$ is an arithmetic mean of functions of the type $\overline{\n}_{n}$ we further deduce that
\beq\label{functv}
\|\n_j\|_{BMO(\TT)}\leq 10\,,
\eeq
and hence applying John-Nirenberg inequality one has that there exists a universal constant $c>0$ such that for any $\g>0$ and any $j\in \N$
\beq\label{JN}
|\{x\in \TT\,|\,\n_j(x)>\g\}|\leq e^{-c\,\g}\;.
\eeq
Next, for any $C>0$, we have that
\beq\label{use}
\|\n^{[k]}\|_1\leq C\,\log k \,+\,\sum_{j=1}^k \|\n_j\big|_{\n_j>C\log k}\|_1\:.
\eeq
Choosing now in \eqref{JN} $\g=C\,\log k$ with $C=\frac{1}{c}$ and replacing it in \eqref{use} we conclude
\beq\label{c1}
\|\n^{[k]}\|_1\leq (C+1)\,\log k \:,
\eeq
thus proving \eqref{fa2}.

We pass now to the proof of \eqref{fa3}. Recall that we want to show that, if $\P$ contains a family of tiles $\F$ which is a suitable \textbf{CME} (for the moment being, the reader is invited to think at $\F$ in the more vague terms described in the Main Heuristic corresponding to \eqref{chain}; later on, if desired, one can consult the precise version displayed in Definition \ref{cme}), then there exists $\bar{C}>0$ such that for any
$k\in\N$ one has
\beq\label{fa5}
\|\n^{[k]}\|_{1,\infty}\geq \bar{C}\,\log k \:.
\eeq
First we present the heuristic for why would one believe such a statement. This is based on the following list of loosely stated observations:
\beq\label{obs}
\eeq
\begin{itemize}
 \item for $\g>>C\,\log k$ the level sets $\{\{\n_j>\g\}\}_{j=1}^{k}$ do not significantly contribute to the norm $\|\n^{[k]}\|_{1,\infty}$;

 \item similarly, for $\g<< C\,\log k$ the level sets $\{\{\n_j<\g\}\}_{j=1}^{k}$ can not provide an estimate of type
 \eqref{fa5};

 \item there exists a conformation of tiles $\P$ such that for suitable $C_1,\, C_2$ absolute positive constants
 \beq\label{obs1}
 |\{\n_j> C_1 \log k\}|\geq \frac{C_2}{k}\:\:\:\:\:\:\:\:\:\:\forall\:\frac{k}{2}<j\leq k\;.
\eeq
 \end{itemize}
 The first two items are just simple consequences of \eqref{JN} and \eqref{use}. The third item will be a direct byproduct of the construction of the \textbf{CME}
 $$\F=\bigcup_{j=\frac{k}{2}+1}^{k}\F_j$$
 presented in the next section.

 Now, if one would assume that it is possible that on top of property \eqref{obs1} one could arrange that the functions
$\{\n_j\}_{j=\frac{k}{2}+1}^{k}$ behave morally as if they were \textbf{independent} random variables then we would immediately conclude that
\beq\label{obs2}
 |\{\n^{[k]}> C_1 \log k\}|\gtrsim\sum_{j=\frac{k}{2}+1}^{k}|\{\n_j> C_1 \log k\}|\geq_{see\:\eqref{obs1}} \frac{C_2}{2}\;,
\eeq
 thus proving \eqref{fa5}.

 The main point of the construction in Section 6 is that the \textbf{CME} as given by Definition \ref{cme} provides exactly
 the mutual independence behavior of $\{\n_j\}_{j=\frac{k}{2}+1}^{k}$ mentioned above.

 Properties \eqref{obs1} and \eqref{obs2} will be a byproduct of the construction of $\F$ presented in the next section (see Definition \ref{cme}). Indeed, writing with the usual notations $\F=\bigcup_{j=\frac{k}{2}+1}^{k}\F_j$, one will be able to  decompose each multi-tower $\F_j$ into a controlled number of towers, that is $\F_j=\bigcup_{l=1}^{\log k}\F_j^l$, and then deduce the following key properties\footnote{In what follows we refer to $Basis (\F_j^{\log k})$  as $\bigcup_{\F_{j}^{\log k,r}} \tilde{I}Top(\F_{j}^{\log k,r})$ where here $\F_{j}^{\log k,r}$ ranges through the decomposition of $\F_j^{\log k}$ into maximal USGTF's. For more details one is invited to consult Section 7.}
\beq\label{obs3}
|Basis (\F_j^{\log k})|\geq \bar{c}\,\frac{1}{ k}\:\:\:\forall \;j\in\{\frac{k}{2}+1,\ldots,\, k\}\;,
\eeq
and
\beq\label{obs4}
Basis (\F_{j_1}^{\log k})\cap Basis (\F_{j_2}^{\log k})=\emptyset\:\:\:\forall \;\:j_1\not=j_2\in\{\frac{k}{2}+1,\ldots,\, k\}\;,
\eeq
where here $\bar{c}>0$ is an absolute constant.

Relations \eqref{obs3} and \eqref{obs4} will then imply
\begin{itemize}
\item $\forall \;j\in\{\frac{k}{2}+1,\ldots,\, k\}$ and $\forall \;m\in\{2^{j-1}+\log\log k,\ldots,\, 2^{j}\}$
\beq\label{obs5}
|\{\overline{\n}_{m}\geq \log k\}|\geq \bar{c}\,\frac{1}{k}\;,
\eeq
\item $\forall \;m_1\in\{2^{j_1-1}+\log\log k,\ldots,\, 2^{j_1}\}$, $\forall \;m_2\in\{2^{j_2-1}+\log\log k,\ldots,\, 2^{j_2}\}$ and
$\forall \;j_1\not=\,j_2\in\{\frac{k}{2}+1,\ldots,\, k\}$ one has
\beq\label{obs6}
\{\overline{\n}_{m_1}\geq \log k\}\cap \{\overline{\n}_{m_2}\geq \log k\}=\emptyset\;.
\eeq
\end{itemize}

This ends our description on the motivation and main properties of the grand maximal counting function.

\section{The structure of the set $E$ via the collection of tiles $\P$. Definition of \textbf{CME}.}
$\newline$

In this section we will make a certain choice for $N$. This will not be done directly but through the structure that we impose on the set of tiles $\P$. More precisely, as described above, we will run an algorithm of constructing a chain of multi-towers with some prescribed properties, this way giving rise to the concept of \textbf{CME}.
$\newline$

 We start with several general observations/heuristics:
\begin{itemize}
\item our construction of the tile configurations will focus on the set\footnote{In reality the set $\F$ will have a more complicated structure; for example $\F$ will also contain tiles from the families $\{\P_j\}_{j< 2^{\frac{k}{2}}}$. Indeed, during the construction algorithm we will express $\F=\bigcup_{j=\frac{k}{2}}^{k}\F_j$ with each tower $\F_j$ further decomposed as $\F_j=\bigcup_{l=1}^{\log k}\F_{j}^l$. While each of the families $\{\F_{j}^l\}_{l=1}^{\log k-2}$ will only have tiles within the set $\bigcup_{j\geq 2^{\frac{k}{2}}}\P_j$, the remaining families $\F_{j}^{\log k-1}$ and $\F_{j}^{\log k}$ would contain also tiles from the set $\bigcup_{j<2^{\frac{k}{2}}}\P_j$. However, in the spirit of the reasons in Section 9, the tile set $\bigcup_{j<2^{\frac{k}{2}}}\P_j$ will play a secondary role in the behavior of $\|T(f_k)\|_{1,\infty}$.}
$$\F\approx\bigcup_{j=2^{\frac{k}{2}}}^{2^{k}}\P_j\,.$$
 Later we will show that, as a consequence of our choice for the tile structure, the contribution of the tiles in $\bar{\P}\setminus\F$ to the $L^{1,\infty}-$``norm"  of our operator $T$ is small in an appropriate sense.

 \item depending on the mass parameter, we will partition the set $\F$ into $\frac{k}{2}$ (dyadic) levels
(preparing thus the future generations):
 $$\F=\bigcup_{l=\frac{k}{2}+1}^{k}\F_{l}\,$$
 with each
 $$\F_l\approx\bigcup_{j=2^{l-1}+1}^{2^l} \P_{j}\:.$$

\item our main task will be to design our set of tiles such that each $\F_l$ is a multi-tower of generation (roughly) $(2^{l-1}+1,\,2^l)$ and of height $\log k$. This construction will be realized through an \textit{inductive process} that will move downwards from the  largest $l$, i.e. $l=k$, towards the smallest one, i.e. $l=\frac{k}{2}+1$.
\end{itemize}
$\newline$

We now present an \textbf{outline} of this process:

\begin{itemize}
\item At the first stage, we will design
        $$\F_k=\bigcup_{l=1}^{\log k}\F_k^{l}\,,$$ such that $\F_k$ is a tower of generation (roughly) $(2^{k-1}+1,\,2^k)$ and of height $\log k$.

        For this we require that
\begin{itemize}
\item each $\F_k^l$ is a USGTF of generation\footnote{In the actual construction process, for technical reasons, we will require that $\F_k^l$ is a USGTF of generation $(2^{k-1}+\log \log k, 2^{k}).$ Same observation applies to the other multi-towers at level $j$, i.e. the actual generation will be $(2^{j-1}+\log \log k, 2^{j})$. The appearance of the factor $\log\log k$ relies on the following loose statement: within the structure formed by the tiles at the bottom scale of each USGTF of generation $(2^{j-1}+\log \log k, 2^{j})$ we can embed towers of generation $(2^{j-2}+\log \log k, 2^{j-1})$ and height precisely $\log k$, this being the hight threshold that plays an important role in our proof. Another way of saying this, is that the bottom structure - mass and number of the tiles at the bottom - of a USGTF determines the height of a tower of a given (smaller) generation that can be embedded within it.}
    $\approx(2^{k-1}+1, 2^{k})$.

\item the set of frequencies $\a(\F_k^l)$ sits entirely below and largely separated from that corresponding to $\a(\F_k^{l+1})$.

\item $ITop(\F_k^{l+1})\prec IBtm(\F_k^{l})$.

\end{itemize}

\item In general, having constructed the (multi-)tower $\F_{j+1}$ of generation $(2^{j}+1,\,2^{j+1})$ and height $\log k$ we will divide it in maximal USGTF's $\{\F_{j+1,r}\}_r$ and within each USGTF $\F_{j+1,r}$ we will embed a specially designed multi-tower $\F_{j}[\F_{j+1,r}]$ of generation $(2^{j-1}+1,\,2^{j})$ and height $\log k$.

\item We will repeat this algorithm until we exhaust the family $\F$ by reaching the level $j=\frac{k}{2}+1$.

\end{itemize}

Now, let us make the above description precise.
$\newline$

\noindent\textbf{$1^{st}$ stage.} \textbf{Constructing the tower} $\F_{k}\;.$
$\newline$

As mentioned before, we will split our  family $\F_{k}$ into $\log k$ sets
$$\F_{k}=\bigcup_{l=1}^{\log k}\F_{k}^{l}\,,$$
with each $\F_{k}^{l}$ being a USGTF of generation $(2^{k-1}+\log \log k, 2^{k})$ (except for $\F_{k}^{\log k-1}$ and $\F_{k}^{\log k}$ which are USGTF's of generation $(2, 2^{k})$). Based on the description made in the previous section it will be enough to specify three parameters: the top, the bottom and the frequency set of each $\F_{k}^{l}$:
We will proceed by induction.
$\newline$

\noindent\textsf{Step $1$.} \textbf{Defining} $\F_{k}^{1}$.
$\newline$

  The key parameters of $\F_{k}^{1}$ are:
\begin{itemize}
\item the top
$$ITop(\F_k^1):=\{[0,1]\}\,.$$

\item the frequency set
$$\a(\F_k^1):=\{2^{2^{2^{2^k}}+100m}\}_{m\in \{0,\,2^{2^k-1}-1\}}\,.$$

\item the bottom
$$IBtm(\F_k^1):=\I_{2^{k-1}-\log\log k}([0,1])\,.$$
\end{itemize}
$\newline$

\noindent\textsf{Step $2$.} \textbf{Defining} $\F_{k}^{2}$.
$\newline$

 The parameters of $\F_{k}^{2}$ are:
\begin{itemize}
\item the top
$$ITop(\F_k^2):=\I^{rt}_{2^{k-1}-\log\log k}([0,1])\,.$$

\item the frequency set
$$\a(\F_k^2):=\{2^{2^{2^{2^k}}+100m}\}_{m\in \{2^{2^k-1},\,2^{2^k}-1\}}\,.$$

\item the bottom
$$IBtm(\F_k^2):=\I_{2^{k-1}-\log\log k}[\I^{rt}_{2^{k-1}-\log\log k}([0,1])]\,.$$
\end{itemize}
$\newline$

\noindent\textsf{Step $l$.} \textbf{Defining $\F_{k}^{l}$ assuming $\F_{k}^{l-1}$} (here $2\leq l\leq \log k-2$)
$\newline$

 Assume we have given $\F_{k}^{l-1}$:
\begin{itemize}
\item the top
$$ITop(\F_k^{l-1}):=\I\,.$$

\item the frequency set
$$\a(\F_k^{l-1}):=\{2^{2^{2^{2^k}}+100m}\}_{m\in \{(l-2)\,2^{2^k-1},\,(l-1)\,2^{2^k-1}-1\}}\,.$$

\item the bottom
$$IBtm(\F_k^{l-1}):=\I_{2^{k-1}-\log\log k}[\I]\,.$$
\end{itemize}
$\newline$

Then $\F_{k}^{l}$ is given by
$\newline$

\begin{itemize}
\item the top
$$ITop(\F_k^l):=\I^{rt}_{2^{k-1}-\log\log k}[\I]\,.$$

\item the frequency set
$$\a(\F_k^l):=\{2^{2^{2^{2^k}}+100m}\}_{m\in \{(l-1)\,2^{2^k-1},\,l\,2^{2^k-1}-1\}}\,.$$

\item the bottom
$$IBtm(\F_k^l):=\I_{2^{k-1}-\log\log k}[\I^{rt}_{2^{k-1}-\log\log k}[\I]]\,.$$
\end{itemize}
$\newline$

\textsf{Step $\log k-1$ and $\log k$.} \textbf{Defining $\F_{k}^{\log k-1}$ and $\F_{k}^{\log k}$.}
$\newline$

For the last two USGTF's we make some minor changes. We will require that both $\F_{k}^{\log k-1}$ and $\F_{k}^{\log k-1}$
be of generation $(1, 2^k)$ and assuming that we are given $\F_{k}^{\log k-2}$ we have:

 $\F_{k}^{\log k-1}$ is given by

\begin{itemize}
\item the top
$$ITop(\F_{k}^{\log k -1}):=\I^{rt}_{0}[IBtm(\F_{k}^{\log k-2})]\,.$$

\item the frequency set
$$\a(\F_{k}^{\log k-1}):=\{2^{2^{2^{2^k}}+100m}\}_{m\in \{(\log k-2)\,2^{2^k-1},\,(\log k-1)\,2^{2^k-1}-1\}}\,.$$

\item the bottom
$$IBtm(\F_{k}^{\log k-1}):=\I_{2^{k}-1}[ITop(\F_{k}^{\log k -1})]\,.$$
\end{itemize}
$\newline$

Finally, the set $\F_{k}^{\log k}$ is given by
$\newline$

\begin{itemize}
\item the top
$$ITop(\F_{k}^{\log k}):=ITop(\F_{k}^{\log k-1})\,.$$

\item the frequency set
$$\a(\F_{k}^{\log k}):=\{2^{2^{2^{2^k}}+100m}\}_{m\in \{(\log k-1)\,2^{2^k-1},\,\log k\,2^{2^k-1}-1\}}\,.$$

\item the bottom
$$IBtm(\F_{k}^{\log k}):=IBtm(\F_{k}^{\log k-1})\,.$$
\end{itemize}

This ends the process of defining the set $\F_{k}$.
$\newline$

\noindent\textbf{$2^{nd}$ stage.} \textbf{Constructing the family $\F_{k-1}\;.$}
$\newline$

We start with the following observation: the family $\F_{k-1}$ has $\log k$ disjoint components according
to the information carried by the graph of $N$ within each of the previously constructed family  $\{\F_{k}^{l}\}_{l\in\{1,\ldots, \log k\}}$. Thus we actually have
$$\F_{k-1}=\bigcup_{l=1}^{\log k} \F_{k-1}[\F_{k}^{l}]\:.$$
For the particular case $l=\log k-1$ and $l=\log k$ we already have determined $\F_{k-1}[\F_{k}^{l}]$ since the sets $\F_{k}^{l}$ are themselves completely determined (up to tiles of mass one) by the requirement that they are USGTF's of generation $(1, 2^k)$. (This is in contrast with the case $l<\log k-1$ where we only require that $\F_{k}^{l}$ is a USGTF of generation $(2^{k-1}+\log \log k, 2^{k})$.)

Thus it only remains to discuss the construction of the families $\{\F_{k-1}[\F_{k}^{l}]\}_{l=1}^{\log k-2}$. In our algorithm we will demand that each $\F_{k-1}[\F_{k}^{l}]$ be a multi-tower of generation $(2^{k-2}\,+\,\log\log k,\,2^{k-1})$ and height $\log k$ embedded into $\F_{k}^{l}$.

In what follows, we will only detail the construction of $\F_{k-1}[\F_{k}^{1}]$ since the remaining multi-towers are constructed in the same way adapting our reasonings inside $\F_{k}^{1}$ to the corresponding $\F_{k}^{l}$.

Recall now the properties of $\F_{k}^{1}$:

\begin{itemize}
\item the top
$$ITop(\F_k^1):=\{[0,1]\}\,.$$

\item the frequency set
$$\a(\F_k^1):=\{2^{2^{2^{2^k}}+100m}\}_{m\in \{0,\,2^{2^k-1}-1\}}\,.$$

\item the bottom
$$IBtm(\F_k^1):=\I_{2^{k-1}-\log\log k}([0,1])\,.$$
\end{itemize}

Write now
$$IBtm(\F_k^1)= I^{lt}Btm(\F_k^1)\cup I^{rt}Btm(\F_k^1)\;,$$
and further express
$$I^{lt}Btm(\F_k^1)=\{J_s\}_{s=1}^{2^{2^{k-1}-\log\log k}}\:.$$
(Recall here the index convention $s_1<s_2$ implies $c(J_{s_1})<c(J_{s_2})$.)

Fix such an interval $J_s$ and consider the set
$$\I_{\log\log k}(J_s)=\{I_r^s\}_{r=1}^{\log k}\:.$$
We then define the family
$$\F_{k-1}[\F_{k}^{1}][J_s]\;,$$
consisting of  $\log k$ towers
$$\F_{k-1}[\F_{k}^{1}][J_s]=\bigcup_{r=1}^{\log k} \F_{k-1}[\F_{k}^{1}][I_r^s]\,.$$
Each tower can be decomposed as
$$ \F_{k-1}[\F_{k}^{1}][I_r^s]=\bigcup_{l_1=1}^{\log k} \F_{k-1}^{l_1}[\F_{k}^{1}][I_r^s]\,,$$
with $\F_{k-1}^{l_1}[\F_{k}^{1}][I_r^s]$ a USGTF.

To see this, we first describe the maximal USGTF within each of the towers ($l_1=1$) :

\begin{itemize}
\item $\F_{k-1}^{1}[\F_{k}^{1}][I_r^s]$ is a USGTF of generation $(2^{k-2}+\log\log k, 2^{k-1})$;

\item the top
$$ITop(\F_{k-1}^{1}[\F_{k}^{1}][I_r^s]):=\{I_r^s\}\,.$$

\item the frequency set
$$\a(\F_{k-1}^{1}[\F_{k}^{1}][I_r^s]):=\{2^{2^{2^{2^k}}+100m}\}_
{m\in \{(s-1)\,(\log k\,2^{2^{k-1}-1})+(r-1)\,(2^{2^{k-1}-1}) ,\,(s-1)\,(\log k\,2^{2^{k-1}-1})+r\,2^{2^{k-1}-1}-1\}}\,.$$

\item the bottom
$$IBtm(\F_{k-1}^{1}[\F_{k}^{1}][I_r^s]):=\I_{2^{k-2}-\log\log k}(I_r^s)\,.$$
\end{itemize}

Now, once we have established the base of each tower the rest of procedure should follow the lines of the tower construction described at stage $1$. For the sake of concreteness we will specify the following:

 Assume we have constructed $\F_{k-1}^{l_1-1}[\F_{k}^{1}][I_r^s]$. Then for $l_1\leq \log k-2$ we have
\begin{itemize}
\item $\F_{k-1}^{l_1}[\F_{k}^{1}][I_r^s]$ is a USGTF of generation $(2^{k-2}+\log\log k, 2^{k-1})$;

\item the top
$$ITop(\F_{k-1}^{l_1}[\F_{k}^{1}][I_r^s]):= I^{rt}Btm(\F_{k-1}^{l_1-1}[\F_{k}^{1}][I_r^s])\,.$$

\item the frequency set
$$\a(\F_{k-1}^{l_1}[\F_{k}^{1}][I_r^s]):=\{2^{2^{2^{2^k}}+100m}\}_
{m\in \A_{k,\,r,\,s,\,l_1}}$$
where
\beq\label{AK1}
\A_{k,\,r,\,s,\,l_1}:=
\bigcup_{m=(s-1)\,(\log k\,2^{2^{k-1}-1})+(r-1+l_1-1)\,(2^{2^{k-1}-1})}^{(s-1)\,(\log k\,2^{2^{k-1}-1})+(r+l_1-1)\,2^{2^{k-1}-1}-1} \{m\}
\:\:\:\:\:\:\:\:\:\:\:\:\textrm{if}\:\:\:r+l_1-1\leq \log k\,.
\eeq
and respectively
\beq\label{AK2}
\A_{k,\,r,\,s,\,l_1}:= \bigcup_{m=(s-1)\,(\log k\,2^{2^{k-1}-1})+(r-1+l_1-1-\log k)\,(2^{2^{k-1}-1})}
^{(s-1)\,(\log k\,2^{2^{k-1}-1})
+(r+l_1-1-\log k)\,2^{2^{k-1}-1}-1}\{m\}\:\:\:\:\:\:\:\:\:\:\:\:\textrm{otherwise}\,.
\eeq
\item the bottom
$$IBtm(\F_{k-1}^{l_1}[\F_{k}^{1}][I_r^s]):=\I_{2^{k-2}-\log\log k}(ITop(\F_{k-1}^{l_1}[\F_{k}^{1}][I_r^s]))\,.$$
\end{itemize}

The construction of $\F_{k-1}^{l_1}[\F_{k}^{1}][I_r^s]$ with $l_1\in \{\log k-1,\,\log k\}$ follows a similar
pattern with the following changes (see also the corresponding changes at the $1^{st}$ Stage):
$\newline$

For $\F_{k-1}^{\log k-1}[\F_{k}^{1}][I_r^s]$ we have
$\newline$
\begin{itemize}
\item $\F_{k-1}^{\log k-1}[\F_{k}^{1}][I_r^s]$ is a USGTF of generation $(1, 2^{k-1})$;

\item the top
$$ITop(\F_{k-1}^{\log k-1}[\F_{k}^{1}][I_r^s]):= I^{rt}Btm(\F_{k-1}^{\log k-2}[\F_{k}^{1}][I_r^s])\,.$$

\item the frequency set
$$\a(\F_{k-1}^{\log k-1}[\F_{k}^{1}][I_r^s]):=\{2^{2^{2^{2^k}}+100m}\}_
{m\in \A_{k,\,r,\,s,\,\log k-1}}$$
where we preserve the definitions \eqref{AK1} and \eqref{AK2}.
\item the bottom
$$IBtm(\F_{k-1}^{\log k-1}[\F_{k}^{1}][I_r^s]):=\I_{2^{k-1}-1}(ITop(\F_{k-1}^{\log k-1}[\F_{k}^{1}][I_r^s]))\,.$$
\end{itemize}

$\newline$
For $\F_{k-1}^{\log k}[\F_{k}^{1}][I_r^s]$ we have
$\newline$

\begin{itemize}
\item $\F_{k-1}^{\log k}[\F_{k}^{1}][I_r^s]$ is a USGTF of generation $(1, 2^{k-1})$;

\item the top
$$ITop(\F_{k-1}^{\log k}[\F_{k}^{1}][I_r^s]):= ITop(\F_{k-1}^{\log k-1}[\F_{k}^{1}][I_r^s])\,.$$

\item the frequency set
$$\a(\F_{k-1}^{\log k}[\F_{k}^{1}][I_r^s]):=\{2^{2^{2^{2^k}}+100m}\}_
{m\in \A_{k,\,r,\,s,\,\log k}}\:.$$

\item the bottom
$$IBtm(\F_{k-1}^{\log k}[\F_{k}^{1}][I_r^s]):=IBtm(\F_{k-1}^{\log k-1}[\F_{k}^{1}][I_r^s])\,.$$
\end{itemize}
$\newline$

\begin{o0}\label{Prec} \textit{Notice the following key property of our construction:
$$\a(\F_{k-1}^{l_1}[\F_{k}^{1}][I_r^s])\cap \a(\F_{k-1}^{l_1}[\F_{k}^{1}][I_{r'}^s])=\emptyset\:\:\:\:
\:\:\forall\:r\not=r'\:\:\textrm{and}\:\:\forall\:l_1\:.$$
Moreover, for each $r,\,s$, the sets $\{\a(\F_{k-1}^{l_1}[\F_{k}^{1}][I_r^s])\}_{l_1=1}^{\log k}$ form a partition
of the frequency set $\{2^{2^{2^{2^k}}+100m}\}_
{m\in \{(s-1)\,(\log k\,2^{2^{k-1}-1}),\,s\,\log k\,2^{2^{k-1}-1}-1\}}\,.$}
\end{o0}
$\newline$

This ends the process of defining $\F_{k}$ and $\F_{k-1}$.

We now repeat this algorithm and further construct by induction $\F_{k-2},\ldots, \F_{\frac{k}{2}+1}$.

$\newline$

\noindent\textbf{$3^{rd}$ stage.} \textbf{Constructing a generic tower} $\F_{j}$, $\frac{k}{2}+1\leq j\leq k-2$.
$\newline$

Assume we have constructed the multi-tower $\F_{j+1}$ of height $\log k$. We first write as before it's layer decomposition
$$\F_{j+1}=\bigcup_{l=1}^{\log k} \F_{j+1}^l\,.$$
For the sake of clarity, we mention here the process of obtaining $\{\F_{j+1}^l\}_l$. Thus, $\F_{j+1}^1$ consists of the union of maximal USGTF's of generation $(2^{j}+\log\log k,\,2^{j+1})$
$$\F_{j+1}^1=\bigcup_{m}\F_{j+1}^{1,m}\:,$$
such that

- $\tilde{I}Top(\F_{j+1}^{1,m})$ is maximal with respect to inclusion among all the sets $\tilde{I}Top(\A)$ with $\A$ maximal
USGTF inside $\F_{j+1}$.

- the sets $\{\tilde{I}Top(\F_{j+1}^{1,m})\}_m$ are pairwise disjoint.

Erase now $\F_{j+1}^1$ from $\F_{j+1}$ and repeat the above algorithm to obtain $\F_{j+1}^2$. Continue this process inductively. (Notice that once we reach $l=\log k-1$, the generation of maximal USGTF's in the decomposition of  $\F_{j+1}^l$ changes to $(1,\,2^{j+1})$.) From our construction this process will end in precisely $\log k$ steps.

With this done, fix a family $\F_{j+1}^l$ (here we assume $l\leq \log k-2$ otherwise trivial considerations), and, with the previous notations, write
$$\F_{j+1}^l=\bigcup_{m}\F_{j+1}^{l,m}\:.$$

Now taking $I^{lt}Btm(\F_{j+1}^{l,m})=\{J_s\}_s$, and then $$\I_{\log\log k}(J_s)=\{I_r^s\}_{r=1}^{\log k}\:,$$
we can initiate the same algorithm as in the case of the construction of $\F_{k-1}$ (see above). Adapting the description made at the $2^{nd}$ stage we have: given $J_s$, we design precisely $\log k$ towers with each such
tower being:

- of generation $(2^{j-1}+\log\log k, 2^{j})$ and height $\log k$;

- embedded in $\F_{j+1}^{l,m}$.

- with basis equal with the corresponding $I_r^s$ in the partition of $J_s$.

To make things clear, for each $J_s$, we will construct the family
$$\F_{j}[\F_{j+1}^{l,m}][J_s]\;,$$
consisting of  $\log k$ towers
$$\F_{j}[\F_{j+1}^{l,m}][J_s]=\bigcup_{r=1}^{\log k} \F_{j}[\F_{j+1}^{l,m}][I_r^s]\,.$$

Again as in the case of $\F_{k-1}$, we have
$$ \F_{j}[\F_{j+1}^{l,m}][I_r^s]=\bigcup_{l_1=1}^{\log k} \F_{j}^{l_1}[\F_{j+1}^{l,m}][I_r^s]\,,$$
with  $\F_{j}^{l_1}[\F_{j+1}^{l,m}][I_r^s]$ a USGTF of generation $(2^{j-1}+\log\log k,\,2^{j})$ (excepting the
cases $l_1=\log k-1$ and $l_1=\log k$.)

For expository reasons
we will no longer give further details on this construction since one follows the same steps (with the obvious changes) displayed in the $2^{nd}$ stage of our construction.

In this way, repeating the above algorithms, our construction process ends by specifying the multi-towers
$$\{\F_{j}\}_{j\in\{\frac{k}{2}+1,\ldots,k\}}\,.$$

Lastly, just for having specified the set of tiles $P\in\P_n$ with $n\leq 2^{\frac{k}{2}}$ (though this information is not in any way \textit{essential} for our later reasonings), we slightly modify the structure of the last constructed family $\F_{\frac{k}{2}+1}$ by requiring that all the maximal USGTF's that sit inside are of generation $(1,2^{\frac{k}{2}+1})$.
$\newline$

\begin{d0}\label{cme} We say that $\F\subset\P$ is a (lacunary) Cantor Multi-tower Embedding (\textbf{CME}) if
\beq\label{defcme}
\F:=\bigcup_{l=\frac{k}{2}+1}^{k} \F_{j}\,,
\eeq
with $\F_j$ constructed as above.
\end{d0}
$\newline$

Notice that as a byproduct of the above construction process we have that Theorem \ref{Tcme} holds.
$\newline$

\begin{o0}\label{Cme} One could relax many of the requirements in Definition \ref{cme} above. For example,
the lacunary structure of the frequencies of each $\F_{j}$ as well as the dyadic splitting of the mass parameter (when forming the generations) are just an expression of the particular operator considered in this paper, i.e. the lacunary Carleson operator. Also, the precise number of multi-towers (in this case $\frac{k}{2}$) is irrelevant in general and should be adapted to the nature of the problem under discussion. Thus, the CME structure can easily be adapted to the study of the pointwise convergence of the full sequence of the partial Fourier sums near $L^1$.

However the key property that should be present in every variation on the theme generated by Definition \ref{cme} is that a CME structure is required to have a tile configuration that maximizes the $L^{1,\infty}-$norm of a grand maximal counting function
similar to \eqref{grandmaxcount}.
\end{o0}

\section{Construction of the set(s) $F_j$}

In this section we will focus on defining the sets $\{F_j\}_{j=\frac{k}{2}+1}^{j=k}$ appearing in the definition of $f$. Each set $F_j$ will be constructed independently. Its structure will be completely determined solely by the \textit{normal} component of the family $\F_{j}$ - as given below:

\begin{d0}\label{nmbdc} Let $j\in\{\frac{k}{2}+1,\,k\}$ and $\F_j$ constructed as before. Partition
$$\F_j=\bigcup_{r}\F_{j,r}$$ into maximal USGTF's. Set now
\beq\label{nmc}
\F_{j,r}^{nm}:=\{P\in\F_{j,r}\,|\,I_{P^*}\cap (\R\setminus\tilde{I} Top(\F_{j,r}))=\emptyset\}\:.
\eeq
 We then define the \textbf{normal component} of $\F_j$ as
\beq\label{nmctot}
\F_{j}^{nm}:=\bigcup_{r}\F_{j,r}^{nm}\,,
\eeq
with the \textbf{boundary component} of $\F_j$ given by
\beq\label{bdctot}
\F_{j}^{bd}:=\F_{j}\setminus\F_{j}^{nm}\,.
\eeq
\end{d0}
$\newline$

\subsection{Construction of $F_k$}

As mentioned previously, we will only work with collection $\F_{k}$.

Recall now the following key property:
\beq\label{obsf}
ITop(\F_{k}^{l})=I^{rt}Btm(\F_{k}^{l-1})\:\:\:\:\:\:\:\:\:\forall\:2\leq l\leq \log k-1 \,.
\eeq
As a consequence we further have
\beq\label{meas}
|\tilde{I}Top(\F_{k}^{l})|=\frac{1}{2}|\tilde{I}Btm(\F_{k}^{l-1})|\:\:\:\:\:\:\:\:\:\forall\:2\leq l\leq \log k-1 \,.
\eeq
We next specify the \textbf{size} of $F_k$ and its \textbf{approximate location}.

Thus we will impose on $F_k$ the following conditions:
\begin{itemize}

\item the total measure of $F_k$ is taken $|F_k|\approx 2^{-\log k\,2^{2^{k}}}\,\cdot\,2^{-k}\,\cdot\frac{1}{k}$;

\item  $F_k\subset \tilde{I}Top(\F_{k}^{\log k-1})= \tilde{I}Top(\F_{k}^{\log k})$;

\item for any $P\in \F_{k}^{\log k-1}\cup \F_{k}^{\log k}$ with $I_P\in IBtm(\F_{k}^{\log k-1})\cup IBtm(\F_{k}^{\log k})$ one has
\beq\label{fconc}
\frac{|I_P\cap F_k|}{|I_P|}\approx 2^{\log k}\,|F_k|\:.
\eeq
\end{itemize}

It remains now to specify the \textbf{concrete location} of $F_k$ inside each of the intervals $I_P$ mentioned in
\eqref{fconc}.

First we notice that from our construction $IBtm(\F_{k}^{\log k-1})= IBtm(\F_{k}^{\log k})$ and that
$IBtm(\F_{k}^{\log k})$ consists of pairwise disjoint intervals. Given this,
it will be enough to specify the structure of
$F_k$ inside a given $I\in IBtm(\F_{k}^{\log k})$.

Now, fixing $I\in IBtm(\F_{k}^{\log k})$, our set $F_k$ inside $I$ will be determined by
\begin{itemize}

\item the frequencies of tiles $P\in \F_{k}^{nm}$;

\item the signum on $I$ of the functions $e^{-2\,\pi\,i\,l(\o_P)\,\cdot}\,T_P^{*}(\chi_{[0,1]})(\cdot)$.

\end{itemize}

Notice that given the choice of our frequencies for the tiles inside $\P$, i.e.
$$\{2^{2^{2^{2^k}}+100m}\}_{m\in \{0\,\ldots\,\log k\,2^{2^{k}-1}-1\}}\,$$
we have that for any $P\in \F_k$ one has that
$T_P^{*}$ is highly oscillatory on $I_P^{*}$, or equivalently, the period of the frequency oscillation $e^{2\,\pi\,i\,l(\o_P)\,\cdot}$ is much smaller than $|I_P|$.

Now the process of constructing $I\cap F_k$ will follow a fractal pattern: we will start with the tiles at the lowest frequency and as we move up in the frequency scales we will trim more and more from the possible location(s) of $I\cap F_k$ inside $I$. (A detailed construction is only done for the $j=k$ case;  with some small adaptations, everything can be repeated for the construction of a general $F_j$).

Here are three \textbf{key} observations derived from properties \eqref{support}, \eqref{supportint} and Observations \ref{LocTT} and \ref{NUSGT}:
 \beq\label{keyobs}
 \eeq
\begin{itemize}
 \item the frequencies $\{e^{2\,\pi\,i\,(2^{2^{2^{2^k}}+100m})\,\cdot}\}_{m\in \{0\,\ldots\,\log k\,2^{2^{k}-1}-1\}}$ are oscillating independently.

 \item $\forall\:\:P\in\F_{k}^{nm}$ we have $e^{-2\,\pi\,i\,l(\o_P)\,\cdot}\,T_P^{*}(\chi_{[0,1]})(\cdot)$ keeps the same
 signum\footnote{This is because in definition \eqref{TP} the function $\psi_k$ is smooth, odd with $x\,\psi_{k}(x)\geq 0$ for any $x\in\R$ and $\chi_{[0,1]}\geq 0$.} on $I$ and moreover due to the smoothing effect of the convolution is morally constant on $I$.

 \item if $P,\,P'\in\F_k^{nm}$ ($P\not=P'$) such that $l(\o_P)=l(\o_{P'})$ then
 $$\textrm{supp}\,T_{P}^{*}\cap\textrm{supp}\,T_{P'}^{*}=\emptyset\:.$$

 \end{itemize}

We start by focussing on the first level of our (multi-)tower - $\F_{k}^{1}$. In this USGTF we have $2^{2^k-1}$ frequencies.

 Making use of observations \eqref{keyobs} above, we have that the following function is well defined:
  $$\S[I,\F_{k}^{1}]:\:\:\a(\F_{k}^{1})\:\rightarrow\:\{-1,\,1\}\;,$$
\begin{itemize}
\item given by $$\S[I,\F_{k}^{1}](a)=1$$ if
 \beq\label{zero}
 \forall\:P\in \F_{k}^{1}\cap\F_{k}^{nm}\:\:\textrm{with}\:\:l(\o_P)=a\:\:\textrm{we have}\;\:T_{P}^{*}\equiv 0\:\:\textrm{on}\:\:I\:.
 \eeq
\item given by\footnote{Notice that based on \eqref{TP} and \eqref{nzero} the expression to which the signum function is applied in \eqref{nonzero} is always \textit{real}.}
\beq\label{nonzero}
\S[I,\F_{k}^{1}](a):=\textrm{sgn} \left(\int_{I} e^{-2\,\pi\,i\,a\,x}\,T_P^{*}(\chi_{[0,1]})(x)\,dx\right)\,,
\eeq
if there exists (see the observation below) $P\in\F_{k}^{1}\cap\F_{k}^{nm}$ tile such that
\beq\label{nzero}
 l(\o_P)=a\:\:\textrm{and}\:\:T_P^{*}(\chi_{[0,1]})\:\:\textrm{is not}\:\:0\:\: \textrm{a.e. on}\:\:I\:.
\eeq
\end{itemize}

\begin{o0}\label{SGN}
Under the assumptions made in Observation \ref{LocTT}, such a tile $P$, if exists, is unique; however this last fact is not actually essential in defining $\S[I,\F_{k}^{1}]$ since the signum of the expression defined in the RHS of \eqref{nonzero} remains the same for all the tiles $P\in\F_{k}^{1}\cap\F_{k}^{nm}$ obeying \eqref{nzero}.
\end{o0}

Next, for $a\in \a(\F_{k}^{1})$, we let
 \beq\label{u0}
\U_{\F_k^{1}}[I](a):=\{x\in I\,|\, \textrm{sgn}\,\left(Re\left(e^{2\,\pi\,i\,a\,x}\right)\right)=\S[I,\F_{k}^{1}](a)\}\:.
\eeq
Now remark that $\U_{\F_k^1}[I](a)$ is a union of disjoint, equidistant, same size dyadic intervals. Moreover, one has that
 \beq\label{monot}
\forall\:a,\,b\in \a(\F_{k}^{1})\:\:\:\textrm{with}\:\;\:a<b\:\:\:\Rightarrow\:\:\:|\U_{\F_k^1}[I](a)\cap \U_{\F_k^1}[I](b)|=\frac{1}{2}\,|\U_{\F_k^1}[I](a)|\:.
\eeq

We now simply impose the following requirement on $F_k$:
\beq\label{reqFk1}
I\cap F_k\subseteq \bigcap_{a\in \a(\F_{k}^{1})} \U_{\F_k^1}[I](a)\:.
\eeq

The ends the process of restricting the location of $F_k$ relative to the first level $\F_{k}^{1}$.

Naturally, the same idea is further extended inside each of the higher level USGTF's
 $\{\F_{k}^{l}\}_{l\leq \log k}$. At the end of this inductive process, putting together all the levels, one concludes
\beq\label{reqFk2}
I\cap F_k\subseteq \bigcap_{l=1}^{\log k}\:\bigcap_{a\in \a(\F_{k}^{l})} \U_{\F_k^l}[I](a)\:,
\eeq
where $\U_{\F_k^l}[I](a)$ for $a\in \a(\F_{k}^{l})$ designates the obvious notational extension from the case $l=1$.

Let us write
$$\bigcup_{m} U_m(I)\,,$$
as the decomposition of $\bigcap_{l=1}^{\log k}\:\bigcap_{a\in \a(\F_{k}^{l})} \U_{\F_k^l}[I](a)$ into
maximal (disjoint) dyadic intervals.

Notice that based on \eqref{monot} and  \eqref{reqFk2} one has
\beq\label{reqFk3}
\frac{|\bigcup_{m}U_m(I)|}{|I|}\approx 2^{-\log k \,(2^{2^{k}-1})}\:.
\eeq

 Also, if $a_0=a_0(I,\,k)\in\a(\F_{k})$ stands for the minimal $l(\o_P)$ with $P\in\F_k$ and $I_P\supseteq I$ then write
\beq\label{min}
\U_{F_k}[I](a_0)=\bigcup_{l} R_{l}\,,
\eeq
the decomposition of $\U_{F_k}[I](a_0)$ into maximal disjoint dyadic intervals (as one can easily notice that for the specific case $j=k$, we actually have $a_0=2^{2^{2^{2^k}}}$). Notice that the intervals $\{R_l\}_l$ have the same length and are equidistant inside $I$. Moreover
\beq\label{u}
|\U_{F_k}[I](a_0)|=\frac{|I|}{2}\,.
\eeq

Then, we can finally determine $F_k$ (up to small perturbations within each $U_m(I)$) by requiring that for any $I\in IBtm(\F_{k}^{\log k})$ the following hold:
\beq\label{REQ}
\eeq
\begin{itemize}
\item  the set $F_k$ obeys \eqref{reqFk2}, or equivalently
$$I\cap F_k\subset \bigcup_{m} U_m(I)\:.$$

\item the set $I\cap F_k$ is equidistributed inside $\{U_m(I)\}_m$, i.e.:
$$|U_m(I)\cap F_k|=|U_{m'}(I)\cap F_k|\:\:\:\:\:\:\:\:\forall\:m,\,m'\;.$$

\item for each $R_l\subset I$ one has
\beq\label{rl}
\frac{|R_l\cap F_k|}{|R_l|}\approx 2^{\log k}\,|F_k|\:.
\eeq
As a consequence, from \eqref{min}, \eqref{u} and \eqref{rl}, one immediately has that
\beq\label{imp}
\frac{|I\cap F_k|}{|I|}\approx 2^{-k}\cdot 2^{-\log k \,2^{2^{k}}}\approx 2^{\log k}\,|F_k|\:.
\eeq
\item inside each $U_m(I)$ the set $U_m(I)\cap F_k$ consists of a single dyadic interval called from now $U_m(I,\,F_k)$;
 its position is irrelevant for our purposes but for the sake of clarity pick this dyadic interval such that its left end-point is the center of $U_m(I)$.
\end{itemize}

This ends the construction of $F_k$.

\subsection{Construction of an arbitrary $F_j$}

In the general situation, we start by first decomposing the multi-tower $\F_j$ into maximal towers, i.e. :
\begin{itemize}
\item first we apply the layer-decomposition
$$\F_j=\bigcup_{l=1}^{\log k} \F_j^l\,;$$
\item we then decompose each $\F_j^l$ into maximal USGTF's:
$$\F_j^l=\bigcup_{m}\F_{j}^{l,m}\,;$$
\item  finally we form the maximal chains (of length $\log k$) which are precisely the maximal towers that we were looking for.
\end{itemize}

For each such maximal tower we repeat the main steps from the $j=k$ case. Since this process is pretty straight-forward we will not give further details here. Putting together the set specifications relative to each maximal tower, at the end of the day, we will have constructed the set $F_j$.

\subsection{Consequences of this construction}

In this section we analyze how the presence of a \textbf{CME} structure within $\P$ and the adapted construction of $\{F_j\}_j$ reflects on the properties of our operator $T$ (or $T^{*}$).

In the next section, we will prove the following
$\newline$

\noindent\textbf{Proposition} (\textsf{heuristic}). \textit{The main component of our operator $T$ relative to the $\|\cdot\|_{1,\infty}$-norm is given by}

\beq\label{maincomp}
T_M(f_k):=\sum_{\frac{k}{2}<j\leq k}\sum_{P\in \F_{j}^{nm}} \,T_{P}(2^{\log k\,2^{2^{j}}}\,\chi_{F_j})\:.
\eeq

In what follows we want to present a glimpse of the methods that we will employ in order to prove our main result. As expected,
the fundamental role in our reasonings will be played by
\begin{itemize}
\item the properties of our \textbf{CME} $\F=\bigcup_{j=\frac{k}{2}+1}^{k}\F_{j}$;
\item the properties of the sets $\{F_j\}_{j=\frac{k}{2}+1}^{k}$.
\end{itemize}
 As a consequence of the above items based on the dual formulation
 \beq
\int T_M(f_k)=\int f_k\,T_{M}^{*} (1)\:,
\eeq
 we make the  following fundamental observation:

 \textsf{the real part of the function $ f_k\,T_{M}^{*} (1)$  is a \textit{positive} function who's integral is bounded from
 below by the expression $\|f_k\|_{L\,\log\log L\,\log\log\log\log L}$. }

This is the main reason for proving the following\footnote{It is worth noticing that in general one can show a much stronger estimate for the \textsf{full} operator $C_{lac}$ since we know that in particular $\|H f\|_1\approx \|f\|_{L\log L}$ where here $H$ is the Hilbert transform. However, \textsf{relative to the $L^1$ norm}, our operator $T_M$ becomes an error term due to the specific \textsf{choice} of our \textbf{CME}. Indeed, one can trivially modify the proof of Lemma \ref{l1blow} and obtain that in fact  $\| T_M(f_k)\|_1\approx\|f_k\|_{L\,\log\log L\,\log\log\log\log L}\approx\log k\:.$ In this context, the main contribution is given by the operator representing the tiles at frequency $0$, that is the operator $T_O$ (see the notations/defintions from the next section). Notice that $T_O$ behaves as a variant of the maximal Hilbert transform.}

\begin{l1}\label{l1blow} (\textsf{$L^1$-blowup})
With the previous notations, the following holds:
\beq\label{l1bl}
 \| T_M(f_k)\|_1\gtrsim \log k\:.
\eeq

\end{l1}

\begin{proof}

Based on the construction of our sets $F_j$ and of our family of tiles, we have
$\newline$

\noindent\textit{Claim.} The following relation holds:
\begin{equation}\label{KEY}
\begin{split}
\int\,\textrm{Re}\left(\chi_{F_j}\,T_{P}^{*}(1)\right)(\cdot)&\geq \frac{1}{500}\,\int\chi_{F_j}(\cdot)\,\left|\,\int \psi_{P}(\cdot-y)\,\chi_{E(P)}(y)\,dy\,\right|\\
&=\frac{1}{500}\,\int\,|\chi_{F_j}\,T_{P}^{*}(1)(\cdot)|\:.
\end{split}
\end{equation}

holds for all\footnote{The condition $P\in \F_{j}^{nm}\setminus (\F_{j}^{\log k-1}\cup\F_{j}^{\log k})$ could be relaxed here at the expense of some extra technicalities.} $P\in \F_{j}^{nm}\setminus (\F_{j}^{\log k-1}\cup\F_{j}^{\log k})$.
$\newline$

\noindent\textit{Proof of our claim.} We will appeal repeatedly to the fundamental relations \eqref{zero}-\eqref{reqFk3}.

Assume wlog that
\beq\label{asump}
P\in \F_{j}^{nm}\cap\F_{j}^{l}\:\:\textrm{for some}\:\:l\in\{1,\ldots, \log k-2\}\,.
\eeq
Now, from the way in which we constructed the set $F_j$ we know that $F_j\subseteq\tilde{I}Btm(\F_{j}^{\log k})$. As a consequence we have
\beq\label{be1}
\int\,\textrm{Re}\left(\chi_{F_j}\,T_{P}^{*}(1)\right)(\cdot)=
\sum_{I\in \tilde{I}Btm(\F_{j}^{\log k})}\int_{I}\,\textrm{Re}\left(\chi_{F_j}\,T_{P}^{*}(1)\right)(\cdot)\:.
\eeq
Next, from our assumption \eqref{asump}, we have that $\tilde{I}Btm(\F_{j}^{\log k})\cap I_{P}^{*}\not=\emptyset$ and that if
$I\in IBtm(\F_{j}^{\log k})$ with $I\cap I_{P}^{*}\not=\emptyset$ then $I\subseteq I_{P}^{*}$ and $|I|<2^{-1000}\,|I_{P}|$ as long as $k$ is large enough.

Thus, it is enough to restrict our attention to a fixed interval $I\in IBtm(\F_{j}^{\log k})$ with $I\subseteq I_{P}^{*}$.

Define now
\beq\label{jp}
\J_{P}:=\{J\:\textrm{dyadic}\,|\,J\subset I_{P}^{*},\: |J|=\frac{1}{4}\,l(\o_P)^{-1}\}\:.
\eeq
Notice that from our construction of the \textbf{CME}  $\F$ we have that $J\in \J_{P}$ implies $|J|<<|I|$.

With the notations from Section 6 we define\footnote{Notice that definition \eqref{nonzero} can be extended from $I\in IBtm(\F_{j}^{l})$ to any (dyadic) interval $J\in \J_{P}$.}
 \beq\label{u1}
\U_{P}[I]:=\bigcup_{{J\in\J_{P}}\atop{J\subset I}}\{x\in J\,|\, \textrm{sgn}\,\left(Re\left(e^{2\,\pi\,i\,l(\o_P)\,x}\right)\right)=\S[J,\F_{j}^{l}](l(\o_P))\}\:.
\eeq
Observe that there exists a subcollection of intervals inside $I$ and belonging to $\J_{P}$, denoted with $\J_{P}^{+}[I]$, such that
\beq\label{u2}
\U_{P}[I]=\bigcup_{J\in\J_{P}^{+}[I]}J\:.
\eeq
Moreover, as a consequence of our $F_j$ construction, we also have
\beq\label{u3}
I\cap F_j=\U_{P}[I]\cap F_j\:.
\eeq
Now one important observation is that due to the choice for the distribution of our frequencies - recall condition \eqref{imN}, any two distinct frequencies $l(\o_1)>l(\o_2)$ in our CME must obey the condition $l(\o_1)\geq 2^{10} l(\o_2)$. Fix now $J\in\J_{P}^{+}[I]$. Then, from our construction of the set $F_j$ we deduce the following uniformity distribution condition:
 \begin{equation}\label{distrib}
\forall\:J_1,\,J_2\in\I_{5}(J)\,\:\:\Rightarrow\:\:
|J_1\cap F_j|=|J_2\cap F_j|=2^{-5}\,|J\cap F_j|\:.
\end{equation}
This implies that for any $J\in\J_{P}^{+}[I]$ one has
 \beq\label{unidistrib}
 \int_{J}\,\textrm{Re}\left(\chi_{F_j}\,T_{P}^{*}(1)\right)(\cdot)\geq
 \frac{1}{4}\,\int\chi_{F_j\cap \frac{1}{16}J}(\cdot)\,\left|\,\int \psi_{P}(\cdot-y)\,\chi_{E(P)}(y)\,dy\,\right|\,.
\eeq

 Thus, as a consequence of \eqref{unidistrib}, \eqref{u3}, \eqref{zero}, \eqref{nonzero} and \eqref{be1}, we have
  \begin{equation}\label{KEY1}
\begin{split}
&\int\,\textrm{Re}\left(\chi_{F_j}\,T_{P}^{*}(1)\right)(\cdot)=
\sum_{{I\in IBtm(\F_{j}^{\log k})}\atop{I\subseteq I_{P}^{*}}}\sum_{J\in\J_{P}^{+}[I]}\int_{J}  |\textrm{Re}\left(\chi_{F_j}\,T_{P}^{*}(1)\right)(\cdot)|\\
&\geq \sum_{{I\in IBtm(\F_{j}^{\log k})}\atop{I\subseteq I_{P}^{*}}}\sum_{J\in\J_{P}^{+}[I]}
\frac{1}{4}\,\int\chi_{F_j\cap \frac{1}{16}J}(\cdot)\,\left|\,\int \psi_{P}(\cdot-y)\,\chi_{E(P)}(y)\,dy\,\right|\\
&\geq\frac{1}{500}\,\int\,\chi_{F_j}(\cdot)\,\left|\,\int \psi_{P}(\cdot-y)\,\chi_{E(P)}(y)\,dy\,\right|\:.
 \end{split}
\end{equation}

This concludes the proof of our claim.
$\newline$

Now, based on \eqref{KEY}, we deduce that
\beq\label{maincomprev}
\int\textrm{Re}\left(\sum_{\frac{k}{2}<j\leq k}\sum_{P\in \F_{j}^{nm}}\,2^{\log k\,2^{2^{j}}}\,\chi_{F_j}\,T_{P}^{*}(1)\right)\approx\int \sum_{\frac{k}{2}<j\leq k}\sum_{P\in \F_{j}^{nm}}\,2^{\log k\,2^{2^{j}}}\,\chi_{F_j} \,|T_{P}^{*}(1)|\:.
\eeq

Thus, we have that
$$|\int T_M(f_k)|\gtrsim \sum_{\frac{k}{2}<j\leq k}\sum_{P\in \F_{j}^{nm}} 2^{\log k\,2^{2^{j}}}\, |F_j\cap I_{P}^{*}|\,\frac{|E(P)|}{|I_P|} $$
$$=\sum_{\frac{k}{2}<j\leq k}\sum_{1\leq l\leq \log k}\sum_{r=2^{j-1}+\log\log k}^{2^j} 2^{\log k\,2^{2^{j}}}\,\sum_{P\in \F_{j}^{nm,l}[r]}\frac{ |F_j\cap I_{P}^{*}|}{|I_P|}\,\frac{|E(P)|}{|I_P|}\,|I_P|$$
$$\gtrsim \sum_{\frac{k}{2}<j\leq k}\sum_{1\leq l\leq \log k}\,2^{\log k\,2^{2^{j}}}\, 2^{j}\,
\frac{|F_j\cap\tilde{I}Top(\F_{j}^{l})|}{|\tilde{I}Top(\F_{j}^{l})|}\,|\tilde{I}Top(\F_{j}^{l})|$$
$$\gtrsim  \log k\,\sum_{\frac{k}{2}<j\leq k} 2^{\log k\,2^{2^{j}}}\,2^j\,|F_j|\approx \|f_k\|_{L\,\log\log L\,\log\log\log\log L}\:,$$
where here for the third relation we used the fact that $\forall\:P\in \F_{j}^{nm,l}[r]$ one has $\frac{ |F_j\cap I_{P}^{*}|}{|I_P|}\approx \frac{|F_j\cap\tilde{I}Top(\F_{j}^{l})|}{|\tilde{I}Top(\F_{j}^{l})|}$.
\end{proof}

Our goal in this paper will be to show that the above Lemma remains true if one replaces in \eqref{l1bl} the $L^1$ norm with an $L^{1,\infty}$ estimate.

\section{Removing the small terms}

As already announced, in this section we intend to show that the contribution of the $(T-T_{M})(f_k)$ component of the operator
is small. More precisely, we will prove that
\beq\label{controlrest}
\|(T-T_{M})(f_k)\|_{1,\infty} \lesssim \|f_k\|_{L\,\log\log L}\:.
\eeq

For this, let us first introduce several notations. Let

$$T_O(f_k):=\sum_{P\in\P(0)} T_{P}(f_k)\,,$$
$$T_0(f_k):=\sum_{P\in\P_0} T_{P}(f_k)\,,$$
$$T_{R,<}(f_k):= \sum_{\frac{k}{2}<j\leq k}\sum_{n< 2^{j-1}+\log\log k}\sum_{P\in\P_n\setminus\F_j} T_P(2^{\log k\,2^{2^{j}}}\,\chi_{F_j})\,,$$
$$T_{R,>}(f_k):= \sum_{\frac{k}{2}<j\leq k}\sum_{n>2^{j}}\sum_{P\in\P_n} T_P(2^{\log k\,2^{2^{j}}}\,\chi_{F_j})\,,$$
and
$$T_R^{bd}(f_k):=\sum_{\frac{k}{2}<j\leq k}\sum_{P\in \F_{j}^{bd}} \,T_{P}(2^{\log k\,2^{2^{j}}}\,\chi_{F_j})\:.$$

Also recall the formula
$$T_M(f_k):=\sum_{\frac{k}{2}<j\leq k}\sum_{P\in \F_{j}^{nm}} \,T_{P}(2^{\log k\,2^{2^{j}}}\,\chi_{F_j})\:.$$

With this, we obviously have
\beq\label{decompT}
\eeq
$$T(f_k)=T_{O}(f_k)\,+T_{0}(f_k)\,+\,T_{R,<}(f_k)\,+\,T_R^{bd}(f_k)\,+\,T_M(f_k)\,+T_{R,>}(f_k)\:.$$

In what follows, we will show that

\beq\label{rests0}
\|T_{O}(f_k)\|_{1,\infty}\lesssim \|f_k\|_{L^1}\:,
\eeq

\beq\label{rests}
\|T_{0}(f_k)\|_{1,\infty}\lesssim \|f_k\|_{L^1}\:,
\eeq

\beq\label{rests1}
\|T_{R,>}(f_k)\|_{1}\lesssim \|f_k\|_{L\,\log\log L}\:,
\eeq

\beq\label{rests2}
\|T_{R,<}(f_k)\|_{1}\lesssim \|f_k\|_{L\,\log\log L}\:,
\eeq
and
\beq\label{rests3}
\|T_{R}^{bd}(f_k)\|_{1}\lesssim \|f_k\|_{L\,\log\log L}\:,
\eeq

First relation is just a consequence of the simple geometric observation that the tiles in $P(0)$ form a tree - where in this case the tree concept must be slightly modified, by replacing the first item appearing in Definition \eqref{tree} with the more relaxed requirement $100 P\leq P_0$.  As a consequence of this observation we have that $T_0$ behaves essentially like the maximal Hilbert transform which we know to be bounded from $L^1$ to $L^{1,\infty}$.\footnote{A stronger estimate that implies the $L^1\,\rightarrow\,L^{1,\infty}$ bound for a tree is the content of Lemma 6.2. in \cite{lvCarl1}.} Thus \eqref{rests0} holds.

Next, relation \eqref{rests} is a direct consequence of Theorem 1.1. part b) proved by the author in \cite{lvCarl1}.

For \eqref{rests2}, we start by decomposing
$$T_{R,<}(f_k)=T_{R,<}^{nm}(f_k)+T_{R,<}^{bd}(f_k)\,,$$f
where
$$T_{R,<}^{nm}(f_k):= \sum_{\frac{k}{2}<j\leq k}\sum_{l<j}\sum_{P\in\F_l^{nm}\setminus\F_j} T_P(2^{\log k\,2^{2^{j}}}\,\chi_{F_j})\,.$$
Next, notice that from the construction of  $\{F_j\}_j$ and $\F$ we have that
\beq\label{kkey}
T_{R,<}^{nm}(f_k)=0\;.
\eeq
For the remaining boundary term, we proceed as follows:

First, fixing $j$, for any $n<2^{j-1}+\log\log k$ we set
$$\F_{<j}^{bd}[n]:=\{P\in \P_n\,|\,P\notin\F_j\:\:\textrm{and}\:\:I_{P^*}\cap F_j\not=\emptyset\}\,,$$
$$\J_{<j}^{bd}[n]:=\{I_{P}\,|\,P\in \F_{<j}^{bd}[n]\}\,.$$
Further on, we define
$$\F_{<j}^{bd}:=\bigcup_{n<2^{j-1}+\log\log k}\F_{<j}^{bd}[n]\:.$$
The key observation here is the following Carleson packing property\footnote{The Carleson packing property can be viewed as a $BMO$-type condition arising naturally from the concept of a Carleson measure; for more details, see \cite{c2}.}: for any $n<2^{j-1}+\log\log k$ and $J\in\J_{<j}^{bd}[n]$
one has
\beq\label{kkey1}
\sum_{{I\in \bigcup_{r\leq n}\J_{<j}^{bd}[r]}\atop{I\subset 1000 J}} |I|\lesssim |J|\;.
\eeq
With this in our mind, we have
$$\|T_{R,<}^{bd}(f_k)\|_1\leq \|\sum_{\frac{k}{2}<j\leq k}\sum_{P\in\F_{<j}^{bd}} T_P(2^{\log k\,2^{2^{j}}}\,\chi_{F_j})\|_1$$
$$\lesssim \sum_{\frac{k}{2}<j\leq k}\sum_{P\in\F_{<j}^{bd}} 2^{\log k\,2^{2^{j}}}\, \frac{|I_{P^*}\cap F_j|}{|I_{P}|}
\frac{|E(P)|}{|I_{P}|}\,|I_{P}|\lesssim \sum_{\frac{k}{2}<j\leq k} 2^{\log k\,2^{2^{j}}}\,2^{\log k}\,|F_j|\lesssim k\,\|f\|_1\,.$$
 Thus  \eqref{rests2} holds.

The proof of \eqref{rests3} is similar in spirit with that of \eqref{rests2} and thus we will not provide further details.

We pass now to \eqref{rests1}. To prove this statement we need to essentially re-do one of the main components of the proof
in \cite{lvKony1}. Setting
\beq\label{defrests1}
T_{R,>}[j]:=\sum_{n>2^{j}}\sum_{P\in\P_n} T_P(2^{\log k\,2^{2^{j}}}\,\chi_{F_j})\:,
\eeq
it is enough to show that
\beq\label{show1}
\|T_{R,>}[j]\|_1\lesssim 2^{\log k\,2^{2^{j}}}\,|F_j|\, 2^{j}\:,
\eeq
since $$T_{R,>}(f_k)=\sum_{\frac{k}{2}<j\leq k}\,T_{R,>}[j]\;.$$

Proceeding as in \cite{lvKony1}, and further using the properties of the set $F_j$ (corresponding to the analogue of \eqref{REQ}), we follow the steps below:

\begin{itemize}
\item  first, for some $g\in L^{\infty}(\TT)$ with $\|g\|_{\infty}=1$, we write
$$\|T_{R,>}[j,2]\|_1=2^{\log k\,2^{2^{j}}}\,\int \chi_{F_j}\,T_{R,>}^{*}[j,2](g)$$
where by definition
$$T_{R,>}^{*}[j,2](g):= \sum_{n\geq 2^{j}}\:\:\sum_{P\in\P_n} T_P^{*}(g)\;.$$

\item next, for $f\in L^1(\TT)$ and $I\in IBtm(\F_j^{\log k})$, define the function
$$\L_{I}(f):=\frac{\int_I f}{|I|}\,\chi_{I}\,,$$
and
$$\L(f):=\sum_{I\in IBtm(\F_j^{\log k})} \L_{I}(f)\:.$$

\item set the operator
\beq\label{av}
T_{R,>}^{*}[av,j,2](g)(\cdot):= \sum_{n\geq 2^{j}}\:\:\sum_{P\in\P_n}
 e^{2\,\pi\,i\,l(\o_{P})\,\cdot}\,\:\L(T_{P^0}^{*}(g))(\cdot)\;,
\eeq
where here if $|I_P|=2^{-k}$ then $T_{P^0}^{*}(g)(x)=-\int \psi_{k}(x-y)\,\chi_{E(P)}(y)\,g(y)\,dy$ is the operator $T_P^{*}(g)$ shifted to the real axis.

\item we have now the key relation (see \cite{lvKony1})
\beq\label{keydif}
\|T_{R,>}^{*}[j,2](g)\|_{L^1(F_j)}\lesssim \|T_{R,>}^{*}[av,j,2](g)\|_{L^1(F_j)}\,+\,|F_j|\;.
\eeq

\item next, decompose $\P_n=\bigcup \p$ into maximal trees, and then, for each $I\in IBtm(\F_j^{\log k})$
and $n\geq 2^{j}$, we further deduce
\beq\label{orthog1}
\eeq

$$|\int_{I\cap F_j} \sum_{\p\in\P_n} e^{2\,\pi\,i\,l(\o_{\p})\,\cdot}\,\L_{I}(T_{{\p}^0}^{*}(g))(\cdot)|
\lesssim
\| \S^{*}_{\P_n}(g) \|_{L^2(I)}\,\left\{ \sum_{\p\in\P_n}
\frac{|<\chi_{F_j\cap I},\,e^{2\,\pi\,i\,l(\o_{\p})\,\cdot}>|^2}{|I|} \right\}^{\frac{1}{2}}\;,$$

where here $\S^{*}_{\P_n}$ is the square-function associated with $\P_n$ given by
$$\S^{*}_{\P_n}:=\left\{\sum_{\p\in\P_n} |{T^{\p}}^{*}|^2\right\}^\frac{1}{2}\;.$$

\item recall now Zygmund's inequality:
\beq \label{zyg}
\|\sum_{j} a_j\,e^{i\,m_j\,x}\|_{\exp(L^2(\TT))}\lesssim \{\sum_{j} |a_j|^2\}^{\frac{1}{2}}\,,
\eeq
where here $\{m_j\}_j\subset\N$ is any lacunary sequence. Applying this (see also \cite{LaDo}, \cite{lvKony1} for a similar treatment) one has
\beq\label{Zyg}
\left\{\sum_{\p\in\P_n}\frac{|<\chi_{F_j\cap I},\,e^{2\,\pi\,i\,l(\o_{\p})\,\cdot}>|^2}{|I|}\right\}^\frac{1}{2}
\lesssim \frac{|I\cap F_j|}{|I|}\,\left(\log \frac{|I|}{|I\cap F_j|}\right)^{\frac{1}{2}}\,|I|^{\frac{1}{2}}\;.
\eeq

\item Due to the \textit{global} control of the $L^2$ norm of $\S^{*}_{\P_n}$ in terms of the mass parameter $n$ we have
\beq\label{masscontrol}
\|\sum_{I\in IBtm(\F_j^{\log k})} \chi_{I}\,\S^{*}_{\P_n}\|_{2}\lesssim 2^{-\frac{n}{2}}\;.
\eeq

\item Now, putting together the $j$ analogue of \eqref{REQ} with \eqref{keydif}, \eqref{orthog1}, \eqref{Zyg} and \eqref{masscontrol}, one concludes
\beq\label{control}
\|T_{R,>}^{*}[j,2](g)\|_{L^1(F_j)}\lesssim |F_j|\,+\, |F_j|\,\sum_{n\geq 2^{j}}\,
 2^{-\frac{n}{2}}\,k \,(\log k\,2^{2^j})^{\frac{1}{2}}\lesssim k^2\,|F_j|\;.
\eeq
\end{itemize}

 Thus, we have actually improved on \eqref{rests1}, proving that in fact
\beq\label{rests11}
\|T_{R,>}(f_k)\|_{1}\lesssim \|f_k\|_{L\,(\log\log \log L)^{2}}\:.
\eeq

\section{Proof of main theorem 1}

In this section we will prove the key result of our paper, encoded in the following

\begin{t1}\label{mainterm} (\textsf{$L^{1,\infty}-$ blowup})

With the previous notations one has\footnote{Notice that as a consequence of Theorem 1, we actually have that $\|T_M(f_k)\|_{1,\infty}\approx\log k$.}
\beq\label{weakblow}
\|T_M(f_k)\|_{1,\infty}\gtrsim \log k\:.
\eeq
\end{t1}

\begin{proof}

Let us start with the following observation: for proving \eqref{weakblow} we need to show that $\exists\:\:G\subset[0,1]$ such that $\forall\:\:G'$ major set, (i.e. $G'\subset G $ and $|G'|\approx|G|$), one has
\beq\label{weakblow1}
\large|\int_{G'} T_M(f_k)\,\large|\gtrsim\log k\:.
\eeq

In order to match Theorem \ref{mainterm} with relation \eqref{blowup} we need to make a ``wise" choice for our set $G$. Recall for the beginning relations \eqref{rests0} and \eqref{rests}. In particular, choosing $G_0=[0,1]$ we know that we can find
$G\subseteq G_0$ major set such that $|G|\geq 1-100^{-1000}$ and
\beq\label{weakblow11}
\int_{G}(|T_{O}(f_k)| +|T_{0}(f_k)|)\lesssim \|f_k\|_1\:.
\eeq

\noindent\textbf{Main Proposition.} \textit{Taking the set $G$ above, for any $G'\subseteq G$ with $|G'|\geq 1-10^{-1000}$ we have}

\beq\label{weakblow2}
\large|\int f_k\,T_M^{*}(\chi_{G'})\,\large|\gtrsim\log k\:.
\eeq

Notice that if we assume for the moment that our claim holds, then, putting together  \eqref{decompT}-\eqref{rests3},
\eqref{weakblow11} and \eqref{weakblow2}, we conclude that Main Theorem 1, and hence Corollary \ref{conjK}, hold.

We start now the proof of our Proposition.

First let us notice that from the previous discussion we have
$$\large|\int f_k\,T_M^{*}(\chi_{G'})\,\large|\geq \int f_k\,\textrm{Re}\left(T_M^{*}(\chi_{G'})\right)$$
$$\gtrsim \sum_{\frac{k}{2}<j\leq k} 2^{\log k\,2^{2^{j}}}\,\sum_{P\in \F_{j}}
\int \chi_{F_j}\,|T_P^{*}(\chi_{G'})|\,.$$
(Strictly speaking in the definition of $T_M$ we only have to deal with the normal components of $\F$, namely $\{\F_j^{nm}\}_j$; however since in the previous section, we proved that the term involving the boundary component is an error term, we will allow in the following estimates terms arising from considering the full family $\F$.)

With this observation we begin the analysis of the structure of the multi-tower $\F_j$. As usual, we first apply the layer-decomposition:
$$\F_j=\bigcup_{l=1}^{\log k} \F_j^l\:.$$
Next, we decompose each $\F_j^l$ into maximal USGTF's:
$$\F_j^l=\bigcup_{m} \F_j^{l,m}\:.$$
We then, introduce the following notations:
$$E(\F_j^{l,m}):=\bigcup_{P\in \F_{j}^{l,m}} E(P)\,, $$
$$G'(\F_j^{l,m}):=E(\F_j^{l,m})\cap G'\,.$$
Also we set
$$T^{\F_{j}}:=\sum_{l=1}^{\log k} T^{\F_{j}^l},\:\:\:\:\:\:\:\:\:\:\:T^{\F_{j}^{l}}:=\sum_{m} T^{\F_{j}^{l,m}}\,,$$
where
$$T^{\F_{j}^{l,m}}:=\sum_{P\in \F_{j}^{l,m}}\,T_P\:.$$

We first claim that if
\beq\label{cut}
|G'(\F_{j}^{l,m})|\geq \frac{1}{2}\,|E(\F_{j}^{l,m})|\:,
\eeq
then there exists $c>10^{-3}$ absolute constant\footnote{The value of the constant $c$ here is not relevant. The condition $c>10^{-3}$ while not at all fundamental for the reasonings to follow can be easily verified and is used only for writing explicit quantitative bounds of the sort of \eqref{difset} and \eqref{conclusion}. Alternatively, one can choose to work with unspecified constants $c$ but then the statement of the Main Proposition above must be adjusted correspondingly.} such that
\beq\label{large}
|\int \chi_{F_j}\,{T^{\F_{j}^{l,m}}}^{*}(\chi_{G'})|\geq c\,2^j\,|F_j\cap \tilde{I}Top(\F_{j}^{l,m})|\:.
\eeq
This is a direct consequence of the special properties of the family $\F$ resulted from our construction.

Indeed, let us first notice that $\F_{j}^{l,m}[2^{j-1}+\log \log k]$ is precisely the bottom of the family $\F_{j}^{l,m}$ (i.e. the time intervals of these tiles form precisely the set $IBtm(\F_{j}^{l,m})$ ).

Set now $$A_{0}^{G'}(P):= \frac{|E(P)\cap G'|}{|I_P|}\,,$$
and let the \textit{heavy} component of $\F_{j}^{l,m}[2^{j-1}+\log \log k]$ be
$$\F_{j}^{l,m}[2^{j-1}+\log \log k](H):=\{P\in\F_{j}^{l,m}[2^{j-1}+\log \log k]\,|\,A_{0}^{G'}(P)\geq \frac{1}{4} A_0(P)\}\:.$$ Now, from \eqref{cut}, we must have that
\beq\label{heavy}
\#\:\F_{j}^{l,m}[2^{j-1}+\log \log k](H)\geq \frac{1}{4}\,\# \F_{j}^{l,m}[2^{j-1}+\log \log k]\;.
\eeq
Observe now that from our construction, more precisely from the USGTF properties \eqref{weightpn}-\eqref{keycompress}, we have: given $I\in ITop(\F_{j}^{l,m})$, $a\in \a(\F_{j}^{l,m})$ and
$s\in[2^{j-1}+\log \log k,\,2^{j}]$ there exists a \textit{unique} $P=[\o_{P},\,I_P]\in\F_{j}^{l,m}$ such that
\beq\label{obss}
\eeq
\begin{itemize}
\item $I_P\subseteq I$;
\item $l(\o_P)=a$;
\item $A_0(P)\in [2^{-s},\,2^{-s+1})$.
\end{itemize}
Moreover, we also have that for any $P,\,P'\in \F_{j}^{l,m}$ with $I_P,\,I_{P'}\subseteq I$ and
$l(\o_P)=l(\o_{P'})=a$ the following is true: $E(P)=E(P')$.

Thus, if for each $I\in ITop(\F_{j}^{l,m})$ and $a\in \a(\F_{j}^{l,m})$ we specify the information $E(\bar{P})\cap G'$ carried by the unique tile $\bar{P}=[\o_{\bar{P}},\,I_{\bar{P}}]\in \F_{j}^{l,m}[2^{j-1}+\log \log k]$ with $l(\o_{\bar{P}})=a$ and
$I_{\bar{P}}\subseteq I$, we then determined the entire structure of the family $\F_{j}^{l,m}$.

Define now
$$\F_{j}^{l,m}(H):=
\{P\in \F_{j}^{l,m}\,|\,\exists\:P'\in \F_{j}^{l,m}[2^{j-1}+\log \log k](H)\:\:\textrm{with}\:\:P'\leq P\}\;,$$
and
$$\F_{j}^{l,m}[s](H):=\F_{j}^{l,m}(H)\cap\F_{j}^{l,m}[s]\;.$$

Then from our previous observations \eqref{obss} and the geometry of our tiles in
$\F_{j}^{l,m}$ we must have that

$$\tilde{\n}_j^{l,m}(H):=\frac{1}{2^{j-1}}\sum_{s=2^{j-1}+\log\log k}^{2^{j}}\sum_{P\in\F_{j}^{l,m}[s](H)}\frac{1}{2^{s-1}}\, \chi_{I_P} \;,$$
obeys
\beq\label{heavycount}
\|\tilde{\n}_j^{l,m}(H)\|_1\geq 10^{-2}\,|\tilde{I}Top(\F_{j}^{l,m})|\,.
\eeq
Since $$|\int \chi_{F_j}\,{T^{\F_{j}^{l,m}}}^{*}(\chi_{G'})|\gtrsim
\sum _{P\in \F_{j}^{l,m}(H)}\int\chi_{F_j}\,|T_{P}^{*}(\chi_{G'})|\,,$$
and
\beq\label{samemaxfunction}
\frac{|I_{P}^{*}\cap F_j|}{|I_{P}^{*}|}\approx
\frac{|\tilde{I}Top(\F_{j}^{l,m})\cap F_j|}{|\tilde{I}Top(\F_{j}^{l,m})|}\approx 2^{l}\,|F_j |\,\:\:\:\;\:\:\:\:\:\:\:\:\forall\:\:P\in\F_{j}^{l,m}\cap\F_{j}^{nm}\;,
\eeq
we deduce that \eqref{heavycount} implies the validity of \eqref{large}.

Moreover, combining \eqref{large} with \eqref{samemaxfunction} and \eqref{propertFj}, we have
\beq\label{boundbel}
2^{\log k\,2^{2^j}}\,|\int \chi_{F_j}\,{T^{\F_{j}^{l,m}}}^{*}(\chi_{G'})|\geq c\,\frac{1}{k}\,2^{l}\,|\tilde{I}Top(\F_{j}^{l,m})|\:.
\eeq

Define now
\beq\label{heavyjl}
\H_{j}^l:=\bigcup_{{m}\atop{\F_{j}^{l,m}\:\textrm{verifies}\:\eqref{cut}\,}}   \F_{j}^{l,m}\:,
\eeq
\beq\label{heavyj}
\H_{j}:=\bigcup_{l}\H_{j}^l\:,
\eeq
and
\beq\label{heavyj}
\H:=\bigcup_{j}\H_{j}\:.
\eeq

Then, based on \eqref{boundbel}, one deduces that
\beq\label{boundbelow}
\int f_k\,\textrm{Re}\left(T_M^{*}(\chi_{G'})\right)\gtrsim \frac{1}{k}\,\sum_{\F_{j}^{l,m}\in\H}\,2^{l}\,|\tilde{I}Top(\F_{j}^{l,m})| \,+\,\sum_{\F_{j}^{l,m}\notin\H} 2^{\log k\,2^{2^{j}}}\,
\int \chi_{F_j}\,|{T^{\F_{j}^{l,m}}}^{*}(\chi_{G'})|\:.
\eeq

Take now $C>0$ a very large constant ($C>2^{100}$ is enough).

Assume by contradiction that
\beq\label{small}
\sum_{\F_{j}^{l,m}\in\H}\,2^{l}\,|\tilde{I}Top(\F_{j}^{l,m})| <\frac{k\,\log k}{C}\,,
\eeq
since otherwise we are done.

Set now
\beq\label{complay}
\V_{j}^l:=\sum_{\F_{j}^{l,m}\in\F_{j}^l}\,2^{l}\,|\tilde{I}Top(\F_{j}^{l,m})|\:,
\eeq
\beq\label{heavylay}
\V_{j}^l(H):=\sum_{\F_{j}^{l,m}\in\H_{j}^l}\,2^{l}\,|\tilde{I}Top(\F_{j}^{l,m})|\:,
\eeq
and respectively
\beq\label{lightlay}
\V_{j}^l(L):=\sum_{\F_{j}^{l,m}\in\F_{j}^l\setminus\H_j^l}\,2^{l}\,|\tilde{I}Top(\F_{j}^{l,m})|\:.
\eeq
Notice that
$$\V_{j}^l=\V_{j}^l(H)\,+\,\V_{j}^l(L)\;.$$
On the one hand from the construction of $\F$ we know
\beq\label{complaybel}
\V_{j}^l\geq 2^{-50}\:\:\:\:\:\:\:\:\:\:\:\forall\:\:\:j,\,l\:.
\eeq
On the other hand, reformulating \eqref{small} we have
\beq\label{data}
\sum_{j=\frac{k}{2}+1}^{k}\sum_{l=1}^{\log k}\V_{j}^l(H) <\frac{k\,\log k}{C}\,.
\eeq
Let now
\beq\label{light}
\D:=\{j\in\{\frac{k}{2}+1,\ldots,\,k\}\,|\,\exists\:l\:\:\textrm{s.t.}\:\V_{j}^{l}(L)\geq \V_{j}^l(H)\}\,.
\eeq
Then, as a consequence of \eqref{complaybel} and \eqref{small}, we have
\beq\label{size}
\#\D\geq \frac{k}{100}\,.
\eeq

For $j\in\D$, let $l_j$ be the smallest value of $l$ that appears in \eqref{light}.

Set now
$$E(j,l_j)=\bigcup_{\F_{j}^{l_j,m}\in\F_{j}^{l_j}\setminus \H_{j}^{l_j}} E(\F_{j}^{l_j,m})\,,$$
$$G'(j,l_j)=\bigcup_{\F_{j}^{l_j,m}\in\F_{j}^{l_j}\setminus \H_{j}^{l_j}} G'(\F_{j}^{l_j,m})\,,$$
and
$$\A_j:=E(j,l_j)\setminus G'(j,l_j)\;.$$
Also define
\beq\label{countfuncj}
\n_j(x):=\frac{1}{2^{j-1}-\log\log k}\sum_{n=2^{j-1}+\log\log k+1}^{2^j} \frac{1}{2^{n-1}}\sum_{P\in\F_j[n]}\chi_{I_P}\;.
\eeq

Then, for any $j\in\D$, we have that
\beq\label{excize}
\eeq
\begin{itemize}
\item $|\A_j|\geq\frac{1}{2}\,|E(j,l_j)|\geq \frac{1}{1000} \sum_{\F_{j}^{l_j,m}\in\F_{j}^{l_j}}\,|\tilde{I}Top(\F_{j}^{l_j,m})|\:;$

\item $\A_j\subseteq\{x\,|\,\n_j(x)\geq l_j\}\;;$

\item $|\A_j|\geq \frac{1}{1000}\,|\{x\,|\,\n_j(x)\geq l_j\}|\;.$
\end{itemize}

Let us now turn our attention towards the level sets of the function(s) $\n_j$.

For $l\in\{1,\ldots,\log k\}$ we decompose the set
\beq\label{not1}
C_{j}^l:=\{x\,|\,\n_j(x)\geq l\}\,,
\eeq
into maximal disjoint (dyadic) intervals
\beq\label{not2}
C_{j}^l=\bigcup_{r} C_{j}^l(r)\:.
\eeq

Now, from our construction of the \textbf{CME} $\F$, we have the following key properties
\beq\label{propertcountf}
\eeq
\begin{itemize}
\item if $j_1\geq j_2$ then for any pairs $(l_1,\,r_1),\,(l_2,\,r_2)$ one has
$$\textrm{either}\:\:C_{j_2}^{l_2}(r_2)\subset C_{j_1}^{l_1}(r_1)\:\:\textrm{or}\:\:C_{j_2}^{l_2}(r_2)\cap C_{j_1}^{l_1}(r_1)=\emptyset\:.$$
\item if we denote with $C_{j_1}^{l_1}(r_1)[j_2,\,l_2]:=\cup_{C_{j_2}^{l_2}(r_2)\subset C_{j_1}^{l_1}(r_1)} C_{j_2}^{l_2}(r_2)$ then one has the \textit{John-Nirenberg} type condition
$$|C_{j_1}^{l_1}(r_1)[j_2,\,l_2]|< 2^{-l_2+10}\,|C_{j_1}^{l_1}(r_1)|\,.$$
\end{itemize}

Indeed, properties \eqref{propertcountf} rely on the key observation
\beq\label{keyobs}
\textrm{given any}\:C_{j}^l(r)\:\:\:\exists\:\F_{j}^{l,m}\:\:\textrm{and}\:\:\exists\:!\:I\in ITop(\F_{j}^{l,m})\:\:
\textrm{s.t.}\:\:I=C_{j}^l(r)\:.
\eeq
Now first item in \eqref{propertcountf} is a direct consequence of \eqref{keyobs} and of the construction of $\F$ which requires
\beq\label{keyobs1}
\begin{split}
if\: j_1\geq j_2\:\:&\textrm{then}\:\:\forall\:\:I_s\in ITop(\F_{j_s}^{l_{j_s},m_s})\:\:\textrm{with}\:\:s\in\{1,\,2\}\:\:\textrm{we have}\\
&\textrm{either}\:\:I_1\cap I_2=\emptyset\:\:\:\textrm{or}\:\:\:I_2\subset I_1\;.
\end{split}
\eeq
The second item in \eqref{propertcountf} follows from \eqref{keyobs}, \eqref{keyobs1} and
\beq\label{keyobs2}
\begin{split}
if\: l_1\leq l_2\:\:&\textrm{then}\:\:\forall\:\:\frac{k}{2}<j\leq k\:\:\textrm{and}\:\:\forall\:I_s\in ITop(\F_{j}^{l_s,m_s})\:\:\textrm{with}\:\:s\in\{1,\,2\}\\
&\textrm{either}\:\:I_1\cap I_2=\emptyset\:\:\:\textrm{or}\:\:\:I_2\subset I_1\:\:\:\textrm{with}\:\:
|I_2|\leq 2^{10}\,2^{l_1-l_2}\,|I_1|\;.
\end{split}
\eeq

Based on \eqref{excize} and \eqref{propertcountf}, one deduces the following fundamental behavior of the sets $\{\A_j\}_{j\in\D}$:
\beq\label{propertAj}
\eeq
\begin{itemize}
\item  $2^{100}\,2^{-l_j}> |\A_j|>2^{-100}\,2^{-l_j} \:.$

\item  $|\A_{j_1}\cap\A_{j_2}|\leq 2^{100}\,|\A_{j_1}|\,|\A_{j_2}|\,.$

\end{itemize}

To see this, we first notice that with notations \eqref{not1} and \eqref{not2} and from \eqref{excize} the first item in
\eqref{propertAj} trivially follows from
\beq\label{nott1}
\A_j\subseteq C_{j}^{l_{j}}\:\:\:\textrm{and}\:\:\:\:|\A_j|\geq \frac{1}{1000}\,|C_{j}^{l_{j}}|\;,
\eeq
and
\beq\label{nott2}
  2^{-55}  \leq 2^{l_j}\,|C_{j}^{l_{j}}|\leq 2^{10}\,,
\eeq
where for the last relation we used \eqref{complaybel} and the last item in \eqref{propertcountf} with $j_2=j$, $j_1=k$, $l_2=l_j$ and $l_1=1$.

For the second item of \eqref{propertAj} we use \eqref{nott1} and first notice that
\beq\label{nott3}
|\A_{j_1}\cap\A_{j_2}|\leq |C_{j_1}^{l_{j_1}}\cap C_{j_2}^{l_{j_2}}|\;.
\eeq
Now, assuming wlog that $j_1\geq j_2$ and making use of the last item in \eqref{propertcountf} we have
\beq\label{nott4}
|C_{j_1}^{l_{j_1}}\cap C_{j_2}^{l_{j_2}}|\leq 2^{-l_{j_2}+10}\,|C_{j_1}^{l_1}|\leq 2^{70}\,|C_{j_1}^{l_{j_1}}|\, |C_{j_2}^{l_{j_2}}|\leq 2^{100} |\A_{j_1}|\,|\A_{j_2}|\;,
\eeq
where here we used again \eqref{nott1} and \eqref{nott2}.

Once at this point, let us notice that
\beq\label{difset}
|G\setminus G'|\geq |\bigcup_{j\in \D} \A_j|\,-\,100^{-1000}\:.
\eeq

Using now \eqref{difset}, \eqref{propertAj}, \eqref{size}, \eqref{complaybel} and the inclusion-exclusion principle we conclude
\beq\label{conclusion}
|G\setminus G'|\geq 2^{-500}\:,
\eeq
which contradicts the requirement that $|G'|\geq 1-10^{-1000}$ thus proving our proposition.

This concludes the proof of Theorem \ref{mainterm} and of Main Theorem 1.
\end{proof}
$\newline$

\section{Proof of main theorem 2}

In this section we will give the proof of Main Theorem 2. The main ingredients that we will use are Theorem \ref{w} and Main Theorem 1. To these ones, we will need to add: for 1 i) the simple observation
\beq\label{simple}
L\log\log L\log\log\log\log L\subset \mathcal{W}\,,
\eeq
noticed in \cite{dp} and for 1 iii) reasonings in the spirit of the proof of Theorem 2.3. in \cite{MMR}.

\subsection{Proof of 1 i)} In this case we just simply remark that
$$\l_{\v}\subseteq\l_{\v_0}=L\log\log L\log\log\log\log L$$
and from \eqref{simple} and Theorem \ref{w} we deduce that $\l_{\v}$ is a $\mathcal{C}_{L}-$space.

\subsection{Proof of 1 ii)} Assume by contradiction that $\l_{\v}$ is a $\mathcal{C}_{L}-$space. Then, there exists $C>0$ such that
\beq\label{RA}
\|T f\|_{L^{1,\infty}}\leq C\,\|f\|_{\l_{\v}}\:\:\:\:\:\:\:\:\:\:\:\forall\:\:f\in\l_{\v}\:.
\eeq

Take now $f=f_k$ with $k,\,f_k$  as designed in the proof of Main Theorem 1. Further, let
\beq\label{dr}
f_k^{*}=\sum_{j=\frac{k}{2}+1}^{k} 2^{\log k\,2^{2^j}}\,\chi_{F_j^{*}}\,,
\eeq
be the decreasing rearrangement of $f_k$.

Recall the definition of the $\l_{\v}-$norm as
$$\|f_k\|_{\l_{\v}}=\int_{0}^{1} f^{*}(t)\,d (\v)(t)\,.$$
Using this and \eqref{dr} we have that
$$\|f_k\|_{\l_{\v}}=\sum_{j=\frac{k}{2}+1}^{k} 2^{\log k\,2^{2^j}}\,\int_{F_j^{*}}d(\v)
\lesssim\sum_{j=\frac{k}{2}+1}^{k} 2^{\log k\,2^{2^j}}\,\v(|F_j|)\:.$$
Next, we notice that $\v_{0}(|F_j|)\approx 2^{-\log k\,2^{2^j}}\,\frac{\log k}{k}\:.$

Thus, we deduce that
\beq\label{fknorm}
\|f_k\|_{\l_{\v}}\lesssim\frac{\log k}{k}\,\sum_{j=\frac{k}{2}+1}^{k} \frac{\v(|F_j|)}{\v_{0}(|F_j|)}\,.
\eeq

At this point, from Main Theorem 1, we know that there exists $C'>0$ such that
\beq\label{T1}
\|T f_k\|_{L^{1,\infty}}\geq C'\,\log k\:.
\eeq
Combining \eqref{fknorm} and \eqref{T1} with assumption \eqref{RA}, we deduce that
\beq\label{T11}
C_0\leq \frac{1}{k}\,\sum_{j=\frac{k}{2}+1}^{k} \frac{\v(|F_j|)}{\v_{0}(|F_j|)}\,,
\eeq
where here is $C_0>0$ is an absolute constant.

Letting now $k\rightarrow \infty$, we notice that the hypothesis $\overline{\lim}_{{s\rightarrow 0}\atop{s>0}}\:\frac{\v(s)}{\v_0(s)}=0$ contradicts \eqref{T1} thus proving that our assumption \eqref{RA} can not be true.

\subsection{Proof of 1 iii)} Assume now that
\beq\label{limit}
\underline{\lim}_{{s\rightarrow 0}\atop{s>0}}\:\frac{\v(s)}{\v_0(s)}=0<\overline{\lim}_{{s\rightarrow 0}\atop{s>0}}\:\frac{\v(s)}{\v_0(s)}\:.
\eeq

We will first show that there exists $\v$ obeying \eqref{limit} such that
\beq\label{llarge}
 \l_{\v_0}=L\log \log L\log\log\log\log L\subsetneq \l_{\v}\subset\W\,,
\eeq
and hence $\l_{\v}$ is a Lorentz $\mathcal{C}_{L}-$space strictly larger than  $\l_{\v_0}$.

Given the analogy between \eqref{llarge} and Theorem 2.3. in \cite{MMR}, we will only present here an outline of the main steps required for our proof (just a simple adaptation of the corresponding steps in \cite{MMR}):

\begin{itemize}
\item Define $\mu:\,[0,1]\,\rightarrow\,\R_{+}$ by $\mu(0)=0$ and $\mu(t):=t\,\log\log \frac{4}{t}$ with $t\in(0,1]$. Notice
that
$$\|f\|_{\W}:=\inf\left\{\sum_{j=1}^{\infty}(1+\log j)\|f_j\|_{\infty}\,\mu(\frac{\|f_j\|_{1}}{\|f_j\|_{\infty}})\:\:\left|\right.\:\:
\begin{array}{cl}
f=\sum_{j=1}^{\infty}f_j,\:\\
\sum_{j=1}^{\infty}|f_j|<\infty\:\textrm{a.e.}\\
f_j\in L^{\infty}(\TT)
\end{array}
\right\}\;.$$
\item Preserving the notations in \cite{MMR}, for $s=\{s_n\}_{n\in\N}$ with each $s_n\in[0,1]$, define $\l^{(s)}$ as the space
of all measurable functions $f:\,\TT\,\rightarrow\,\C$ such that there exists a sequence $\{f_n\}_{n\in\N}$ with $f_n\in L^{\infty}(\TT)$ satisfying  $f=\sum_{n\in\N} f_n$ (with convergence in $L^1(\TT)$) and
$$\sum_{n=1}^{\infty} \max\{\|f_n\|_1,\,s_n\,\|f_n\|_{\infty}\}\,\frac{\mu(s_n)}{s_n}\,(1+\log n) <\infty\;.$$
We endow $\l^{(s)}$ with the norm
\beq\label{normls}
 \|f\|_{\l^{(s)}}:=
 \inf\{\sum_{n=1}^{\infty} \max\{\|f_n\|_1,\,s_n\,\|f_n\|_{\infty}\}\,\frac{\mu(s_n)}{s_n}\,(1+\log n)\,|\,f=\sum_{n\in\N} f_n\}\,.
\eeq
\item Using the simple observation $\mu(\a)\leq \max\{1,\,\frac{\a}{\b}\}\,\mu(\b)$ for any $\a,\,\b\in(0,1)$ one has
that for any sequence $s=\{s_n\}_{n\in\N}\in (0,1]^{\N}$, $\l^{(s)}$ is a r.i. Banach space such that
\beq\label{emb}
\l^{(s)}\hookrightarrow \W\,,
\eeq
with the norm inclusion $\leq 1$.
\item For $s\in(0,1]^{\N}$ as before, let $\v^{(s)}$ be the quasi-concave function on $[0,1]$ defined by
$\v^{(s)}(0)=0$ and
$$\v^{(s)}(t)=\inf_{n\in\N} \max\{t,\,s_n\}\,\frac{\mu(s_n)}{s_n}\,(1+\log n)\,,\:\:\:\:\:\:\:\:\:\:\:\:\forall\:t\in(0,1]\:.$$
Then, for $\tilde{\v}^{(s)}$ the least concave majorant of $\v^{(s)}$, we have
\beq\label{emb1}
\l_{\tilde{\v}^{(s)}}\hookrightarrow\W\,,
\eeq
with inclusion norm smaller or equal than $1$.

To prove \eqref{emb1} one uses \eqref{emb}, the last observation in the Appendix - \eqref{continclusion} and
the fact that the fundamental function of $\l^{(s)}$ is precisely $\v^{(s)}$. This last fact is pretty straightforward and we leave it for the reader.
\item If $s=\{s_n\}_{n\in\N}\in (0,1]^{\N}$ is given by $s_n=2^{-2^{2^n}}$, then
$$\l_{\tilde{\v}^{(s)}}=L\log\log L\log\log\log\log L\:.$$
\end{itemize}
With this done we can now end the proof of \eqref{llarge} trough the following reasoning:

- let $s=\{s_n\}_{n\in\N}\in (0,1]^{\N}$ with $s_n=2^{-2^{2^{2^n}}}$ and define
\beq\label{vdefinition}
\v(t)=\min\{\v_{0}(t),\,\tilde{\v}^{(s)}\}\:.
\eeq
Then we clearly have
$$\l_{\v_{0}}\subseteq \l_{\v}\subseteq \l_{\tilde{\v}^{(s)}}+\l_{\v_{0}}\subset\W\:.$$
In order to now prove that
\beq\label{strict}
\l_{\v_{0}}\subsetneq \l_{\v}\,,
\eeq
we just notice that
$$\v_{0}(s_n)\approx 2^{-2^{2^{2^n}}}\,2^{2^n}\,n\,,$$
while
$$\v^{(s)}(s_n)\leq \mu(s_n)\,(1+\log n)\approx  2^{-2^{2^{2^n}}}\,2^{2^n}\,\log n\,.$$
From here we conclude
$$\lim_{n\rightarrow\infty} \frac{\v(s_n)}{\v_{0}(s_n)}=0\,,$$
thus showing \eqref{strict} and  ending the proof of \eqref{llarge}.
$\newline$

We pass now to proving that there exists $\v$ obeying \eqref{limit} such that
\beq\label{secondpart}
\l_{\v}\:\:\textrm{is not a}\:\:\mathcal{C}_L-\textrm{space}\:.
\eeq
For this, appealing again to Main Theorem 1, we notice that it will be enough to show that for a proper choice of $\v$, one has
 \beq\label{enough}
\underline{\lim}_{k\rightarrow\infty} \frac{1}{\log k}\,\sum_{j=\frac{k}{2}+1}^{k} 2^{\log k\,2^{2^j}}\,\v(|F_j|)\,=0\:.
\eeq
For each $n\geq 100$, take $k_n=2^{2^n}$ in the counterexample provided by Main Theorem 1 and define $l_n$ to be the line passing through the points $A_{\frac{k_n}{2}+1}=(|F_{\frac{k_n}{2}+1}|,\,\v_{0}(|F_{\frac{k_n}{2}+1}|))$ and $A_{k_n}=(|F_{k_n}|,\,\v_{0}(|F_{k_n}|))$, respectively. Define now $\v$ as follows:
\beq \label{syst}
\v(t):=\left\{
                        \begin{array}{ll}
l_n(t)\:\:\:\:\:\:\textrm{for}\:\:t\in[|F_{\frac{k_n}{2}+1}|,\,|F_{k_n}|],\:\:\:n\in \N,\:\:n\geq 100\,,\\\\
\v_0(t)\:\:\:\:\:\:\textrm{otherwise}                      \end{array} \right.
\:.\eeq
One can now easily check that $\v$ verifies \eqref{limit}.

Moreover, for $k=k_n=2^{2^n}$, one has that
$$\sum_{j=\frac{k}{2}+1}^{k} 2^{\log k\,2^{2^j}}\,\v(|F_j|)=\sum_{j=\frac{k}{2}+1}^{k} 2^{\log k\,2^{2^j}}\,l_n(|F_j|)$$
$$\approx  2^{\log k\,2^{2^j}}\,l_n(|F_j|)\,|_{j=\frac{k}{2}+1}\approx\frac{\log k}{k}\,,$$
and thus \eqref{enough} holds.

\subsection{Proof of 2 i)}
$\newline$

Our goal here is to prove that if $X$ is a r.i. Banach space with its fundamental function $\v_{X}=\v$ obeying \eqref{growth} then $X$ is a $\mathcal{C}_L-$space.

Now, based on \eqref{continclusion}, it will be enough to show that \eqref{growth} implies $M_{\v_{*}}\subset\l_{\v_{0}}$. But this last relation follows easily from the following:
\begin{itemize}
\item $f\in M_{\v_{*}}$ implies that there exists $C>0$ such that
\beq\label{drecrear}
f^{**}(t)<\frac{C}{\v(t)}\:\:\:\:\:\:\forall\:t\in(0,1]\:.
\eeq
\item $\|f\|_{\l_{\v_{0}}}=\int_{0}^{1} f^{*}(t)\,d(\v_{0})(t)\approx\int_{0}^1 f^{**}(t)\,[-t\,\v_{0}''(t)]\,dt\:.$
\end{itemize}

\subsection{Proof of 2 ii) and iii)} With the notations made at point 3) of our Main Theorem 2, assume that the following holds:
\beq\label{weakb}
\|T f_k\|_{1,\infty}\gtrsim \|f_k\|_V:=\inf_{\sigma\in S_k} \sum_{j=1}^k r_j\,|F_j|\,2^j\,\log (\sigma(j)+1)\:,
\eeq
where here $S_k$ designates the class of all permutations of the set $\{1\,,\ldots,\,k\}$ (A proof of this statement will be outlined when proving 3) below).

Next, imposing the condition that $\{r_j\}_{1\leq j\leq k}$ is an increasing sequence of positive numbers and using the definition of the Marcinkiewicz space $M_{\v_{*}}$, we have that
\beq\label{marcinkiew}
\|f_k\|_{M_{\v_{*}}}\approx \sup_{1\leq n\leq k}\, \frac{\v(|F_n|)}{|F_n|}\,\sum_{j=n}^k r_j\,|F_j|\:.
\eeq
Thus if we assume that (B) holds, taking $X=M_{\v^*}$, we deduce that there exists $C'>0$ absolute constant such that
\beq\label{ineq}
\inf_{\sigma\in S_k} \sum_{j=1}^k r_j\,|F_j|\,2^j\,\log (\sigma(j)+1)\leq C'\,\sup_{1\leq n\leq k}\, \frac{\v(|F_n|)}{|F_n|}\,\sum_{j=n}^k r_j\,|F_j|\:,
\eeq
holds for any $k\in \N$ large enough.

Choosing now $r_j:=\frac{1}{|F_j|\,2^j\,j}$ relation \eqref{ineq} becomes
\beq\label{ineq1}
\sup_{n\leq k}\frac{\v(|F_n|)}{2^n\,|F_n|\,n}\geq C'' (\log k)^2\,,
\eeq
which further implies
\beq\label{ineq2}
\overline{\lim}_{{s\rightarrow0}\atop{s>0}}\frac{\v(s)}{\v_0(s)}=\infty\:,
\eeq
thus proving 2 ii).

Notice here that for a choice of the form $r_j:=\frac{1}{|F_j|\,2^j\,j\,(\log j)^2}$ one can improve \eqref{ineq2} to
\beq\label{ineq3}
\overline{\lim}_{{s\rightarrow0}\atop{s>0}}\frac{\v(s)}{\v_0(s)\,\log\log\log \frac{1}{s}\,\log\log\log\log \frac{1}{s}}=\infty\:.
\eeq

Passing now to the proof of 2 iii), consider wlog that
\beq\label{assump}
s\,\mapsto\,\frac{\v_0(s)}{\v(s)}\:\:\textrm{increasing on}\:(0,1)\,,
\eeq
i.e. $\ep=1$ (otherwise trivial modifications).

In this setting it is enough to show that (B) implies (A).

Making the choice $r_j=\frac{1}{\v(|F_j|)}$ and using \eqref{assump} we have that
\beq\label{marcinkiew1}
\|f_k\|_{M_{\v_{*}}}\approx \sup_{1\leq n\leq k}\, \frac{\v(|F_n|)}{|F_n|}\,\sum_{j=n}^k \frac{|F_j|}{\v(|F_j|)}\lesssim 1\:,
\eeq
and
\beq\label{ineq1}
\|T f_k\|_{1,\infty}\gtrsim \sum_{j=1}^{k/2}\frac{\v_0(|F_j|)}{\v(|F_j|)}\:.
\eeq

Take now $|F_j|$ such that
    \beq\label{choicefj}
\frac{|F_j|}{\v(|F_j|)}=\int_{y_{j+1}}^{y_j}\frac{s}{\v(s)}\,(-\frac{\v''_0(s)}{\v'_0(s)})\,ds\:.
    \eeq
Then, from \eqref{marcinkiew1},  \eqref{ineq1}, \eqref{choicefj} and the requirement that $M_{\v^*}$ is a $C_L-$space we deduce that
\beq\label{conclude}
\lim_{k\rightarrow\infty}\int_{y_{k/2}}^{1}\frac{-s\,\v''_0(s)}{\v(s)}\,ds\lesssim 1\:,
\eeq
finishing the proof of 2 iii).

\subsection{Some remarks on the proof of Theorem 3.} Notice that our claim follows if we are able to show \eqref{weakb}.

For this though, all that we need is to notice that we can follow the reasonings in the proof of Main Theorem 1 and construct $F_j\subseteq\TT$ measurable such that:
\begin{itemize}
\item the \textit{size} of $F_j$ obeys $|F_j|\in [y_{j+1},\, y_j]$.
\item the \textit{locations} of $\{F_j\}_{j}$ are arising from a \textbf{CME}-structure following the steps at Sections 5 and 6.
\end{itemize}

In this new setting, defining $f_k$ as in \eqref{functk} and proceeding as in Section 8, one obtains
\beq\label{conclusion}
\|T f_k\|_{1,\infty}\gtrsim \inf_{{\sum_{j=1}^k y_j\leq \frac{1}{2}}\atop{y_j\geq 0}} \sum_{j=1}^k r_j\,|F_j|\,\log\log\frac{1}{|F_j|}\,\log\frac{1}{y_j}\approx \|f_k\|_V\:.
\eeq
Thus, in particular \eqref{weakb} holds.

It is worth mentioning that one actually has
\beq\label{functv}
 C'\,\|f_k\|_{\W}\leq \|f_k\|_V\leq C''\,\|f_k\|_{\W}\,,
\eeq
for some absolute $C''\geq C'>0$.
$\newline$

\section{A discussion regarding the (Lacey-Thiele) discretized Carleson model and the (discretized) Walsh model}

This section explains in great detail the comments expressed in Observation \ref{W}. We will thus present a detailed antithesis between the a.e. pointwise convergence properties of the (lacunary) Carleson operator $C_{lac}^{\{n_j\}_j}$ and those corresponding to the (lacunary) Walsh-Carleson operator $C_{W}^{\{n_j\}_j}$ and (more briefly) to the (lacunary) Lacey-Thiele discretized Carleson model.

In fact, in most of our analysis we will insist on analyzing the Walsh-Carleson operator as the corresponding analysis of the Lacey-Thiele discretized Carleson model will then become immediately transparent.

We first list (a very detailed explanation will follow) the key aspects that will make a difference in the behavior of $C_{W}^{\{n_j\}_j}$:
\beq\label{Walshp}
\eeq
\begin{itemize}
\item (I)  the \textbf{algebraic} properties of the Walsh wave-packets;

\item (II) the \textbf{discrete/dyadic} character of the Walsh-Carleson operator.
\end{itemize}

In order to better explain the above items, let us recall\footnote{For this we will follow closely \cite{T}.} some of the definitions/properties of the Walsh system/Walsh-Carleson operator.

We will only discuss the periodic setting since this the one of interest for us. Thus we define the Walsh phase plane
$\textbf{W}=[0,1)\times [0,\infty)$. We fix the canonical dyadic grid (having origin at $0$ and scales in power of two) on both
$[0,1)$ and $[0,\infty)$. As before in our paper, a dyadic interval will be of the form $[2^{-j}\,n,\,2^{-j}\,(n+1))$ with $j\in\N$ (or $j\in\Z$ for positive real axis) and $n\in\N$. Keeping the same notations from Introduction, we refer to $P$ as a tile if $P=I_{P}\times \o_{P}\subset\textbf{W}$ with $I_{P},\,\o_{P}$ dyadic intervals such that $|I_P|\,|\o_{P}|=1$ and we let $\P$ be the collection of all such tiles. Unlike the Fourier setting, here we will also need to work with the so-called \textit{bitiles}, that is dyadic rectangles $R=I_{R}\times \o_{R}\subset\textbf{W}$ of area two (that is $I_{R},\,\o_{R}$ dyadic intervals such that $|I_R|\,|\o_{R}|=2$). We refer to the collection of all bitiles as $\mathcal{R}$.

Next for $R\in\mathcal{R}$ with $I_{R}=[x_0,\,x_1)$ and $\o_{R}=[\xi_0,\,\xi_1)$ we define
\beq\label{sons}
\eeq
\begin{enumerate}
\item the upper son of $R$ as $R_u:=[x_0,\,x_1)\times[\frac{\xi_0+\xi_1}{2},\,\xi_1)\in\P$;

\item the lower son of $R$ as $R_l:=[x_0,\,x_1)\times[\xi_0,\,\frac{\xi_0+\xi_1}{2})\in\P$;

\item the left son of $R$ as $l^R:=[x_0,\,\frac{x_0+x_1}{2})\times[\xi_0,\,\xi_1)\in\P$;

\item the right son of $R$ as $r^R:=[\frac{x_0+x_1}{2},\,x_1)\times[\xi_0,\,\xi_1)\in\P$.
\end{enumerate}

These being said, let us first recall the definition of the Walsh system.

Fix $n\in\N$ and let
\beq\label{2base}
n=\sum_{i=0}^{\infty}\ep_{i}\,2^{i}\:\:\:\textrm{with}\:\:\:\ep_{i}\in\{0,\,1\}\,,
\eeq
be it's dyadic decomposition. Then we define the Walsh system $\{w_n\}_{n\in\N}$ as:
\beq\label{defw}
\eeq
\begin{itemize}
\item if $x\in\R\setminus[0,1)$ then $w_n(x)=0$ for any $n\in\N$;

\item if $x\in[0,1)$ and $n=0$ then we set $w_n(x)=1$;

\item if $x\in[0,1)$ and $n\geq 1$ obeys \eqref{2base} then we let\footnote{Notice that in \eqref{2base} only finitely many $\ep_{i}$'s are nonzero.}
$$w_{n}(x):=\prod_{i=0}^{\infty}\left( \textrm{sgn} (\sin 2^{i+1}\, \pi\,x)\right)^{\ep_i}\;.$$
\end{itemize}
Next, given $P\in\P$  with $P=[2^{-j}\,l,\,2^{-j}\,(l+1))\times[2^j\,n,\,2^j\,(n+1))$ and $j,\,l,\,n\in\N$ we define the associated Walsh wave-packet as
\beq\label{wavepacket}
w_{P}(x)=w_{n,l,j}(x)=2^{\frac{j}{2}}\,w_{n}(2^{j}\,x-l)\;.
\eeq
Given $R\in \mathcal{R}$ and using now definitions \eqref{sons}, \eqref{defw} and \eqref{wavepacket} we have the following key \textit{algebraic} properties of the Walsh wave-packets:
\beq\label{alg1}
w_{R_u}=\frac{1}{\sqrt{2}}\,\left(w_{l^R}-w_{u^R}\right)\;,
\eeq
and
\beq\label{alg2}
w_{R_l}=\frac{1}{\sqrt{2}}\,\left(w_{l^R}+w_{u^R}\right)\;.
\eeq

Now, as a consequence of relations \eqref{alg1} and \eqref{alg2}, we have that the following key identity holds\footnote{For a proof of this statement see \cite{T}.}:
\beq\label{Walshid}
W_{n}f(x):=\sum_{k=0}^{n} <f,\,w_k>\,w_k(x)=\sum_{R\in\mathcal{R}} <f,w_{R_l}> \,w_{R_l}(x)\,\chi_{\o_{R_u}}(n)\;,
\eeq
where here $f\in L^1(\TT)$ and $n\in\N$.

\begin{o0}\label{WW} We deduce from \eqref{Walshid} that the (lacunary) Walsh-Carleson operator obeys
\beq\label{WalshC}
C_{W}^{\{n_j\}_j}f(x)=\sup_{j\in\N}|\sum_{R\in\mathcal{R}} <f,w_{R_l}> \,w_{R_l}(x)\,\chi_{\o_{R_u}}(n_j)|\;.
\eeq
This is precisely the meaning of item (II) that the Walsh-Carleson operator is of \textit{discrete} nature, namely, it can be written as a superposition of projection operators associated with a \textit{single} dyadic grid.
\end{o0}

We end this discussion on the basic properties of the Walsh system by mentioning another relevant \textit{algebraic} feature of it: as a consequence of definition \eqref{defw}, we have
\beq\label{alg3}
\sum_{n=0}^{2^L-1} w_{n}(x)=\prod_{i=0}^{L-1} (r_{0}(x)+r_{2^i}(x))\;,
\eeq
where here $L\in\N$.

This implies that for any $f\in L^1(\TT)$ one has
\beq\label{alg4}
W_{2^L-1}(f)(x)=<f(\cdot),\,\prod_{i=0}^{L-1} \left(r_{0}(x)\,r_{0}(\cdot)+r_{2^i}(x)\,r_{2^i}(\cdot)\right)>\;.
\eeq
One can reshape \eqref{alg4} for other subsequences of partial Fourier-Walsh sums or even differences of partial sums. For example choosing $L>>M$ with $L,\,M\in\N$ one has
\beq\label{alg5}
\begin{array}{ll}
&W_{2^L-1}(f)(x)-W_{2^L-2^M-1}(f)(x)=\\\\
&<f(\cdot),\,\left(\prod_{i=M}^{L-1}r_{2^i}(x)\,r_{2^i}(\cdot) \right) \left(\prod_{i=0}^{M-1} \left(r_{0}(x)\,r_{0}(\cdot)+r_{2^i}(x)\,r_{2^i}(\cdot)\right)\right)>\;.
\end{array}
\eeq

Let us now switch our attention towards the Carleson operator $C$ (or it's lacunary version $C_{lac}^{\{n_j\}_j}$). For making the analogy with the Walsh case more transparent we will use a different, cleaner decomposition than the one in \cite{lt3}.\footnote{This simple, elegant decomposition of the real line Carleson operator appears in junior paper \cite{DJ} written under the supervision of E. Stein and sporadic guidance of the author.}

Following the Walsh case development, we first define the real Fourier phase plane as \textbf{$F_r$}$=\R\times\R$ and the periodic Fourier phase plane\footnote{For technical reasons, in order to remove the boundary terms we enlarge the canonical interval $[0,1)$ to $[-1,2)$.} as $\textbf{$F_p$}=[-1,2)\times \R$. Unlike the Walsh case, in order to connect $C$ with its \textit{discretized} model, we will need to use a continuum of grids. Let $\ll,\,\mu$ be parameters ranging in $[0,1)$.
On the two real axis we will use ``dyadic" grids defined by the following structure:
\beq\label{inte}
\eeq
\begin{itemize}
\item in space  - $\mathcal{I}_{-m}^{-\ll,y}$ that is the set of intervals $2^{-m-\ll}[n+y,\,n+y+1)$ with $m\in \Z$, $n\in\Z$, $y\in[0,1]$.

\item in frequency - $\mathcal{J}_{m}^{\ll,\mu}$ that is the set of intervals $2^{m+\ll}[n+\mu,\,n+\mu+1)$ with $m\in Z$ and $n\in\Z$.
\end{itemize}

Define now

- the $(y,\ll,\mu)$-collection of tiles (real case)
\beq\label{Collect}
 \P^{y,\ll,\mu}:=\bigcup_{m\in\Z} \P^{y,\ll,\mu}_{m}:=\bigcup_{m\in\Z}\mathcal{I}_{-m}^{-\ll,y}\times\mathcal{J}_{m}^{\ll,\mu}\:.
\eeq

- the $(y,\ll,\mu,+)$-collection of tiles (periodic case)
\beq\label{Collectp}
 \P^{y,\ll,\mu,+}:=\bigcup_{m\in\N} \P^{y,\ll,\mu,+}_{m}:=\bigcup_{{m\in\N}}\mathcal{I}_{-m}^{-\ll,y,+}\times\mathcal{J}_{m}^{\ll,\mu}\:,
\eeq
where here $\mathcal{I}_{-m}^{-\ll,y,+}=\{I\in\mathcal{I}_{-m}^{-\ll,y}\,|\,I\subset [-1,2)\}$.

We pass now to defining the Fourier wave-packets adapted to $\P^{y,\ll,\mu}$ (with the obvious changes for the periodic case).

Let $\phi\in \mathcal{S}(\R)$ such that $\textrm{supp}\:\hat{\phi}\subseteq [-0.1,\,0.1]$, $\hat{\phi}\geq 0$ and
$\hat{\phi}\equiv 1$ on $[-0.07,\,0.07]$.

Next, we introduce the classes of symmetries entering in the structure of the Carleson operator:
\begin{itemize}
\item translations - $T_{z}f(x)=f(x-z)$ with $x,\,z\in\R$;

\item modulations - $M_{\xi}f(x)=e^{2\pi i x\xi}\,f(x)$ with $\xi\in\R$;

\item dilations - $D_{\ll}^{p}f(x)=\ll^{-\frac{1}{p}}\,f(\ll^{-1}x)$ with $\ll>0$ and $p\in(0,\infty]$.
\end{itemize}

Let now $P$ be a generic tile, that is belonging to $\bar{\P}=\bigcup_{y,\ll,\mu\in[0,1)}\P^{y,\ll,\mu}$. If $P=I_{P}\times\o_{P}$ (recall that $P$ has area one) then we set $c(\o_P)$ the center of $I_{P}$ let $\o_{P_{l}}=(-\infty,\,c(\o_{P})]\cap\o_{P}$ and $\o_{P_{u}}=\o_{P}\setminus\o_{P_{l}}$ and as before define
\begin{enumerate}
\item the upper son of $P$: $P_u=I_{P}\times \o_{P_{u}}$;

\item the lower son of $P$: $P_l=I_{P}\times \o_{P_{l}}$.
\end{enumerate}
We are now ready to define the wave-packet associated with $P_{l}$ by
\beq\label{Fwv}
 \phi_{P_l}(x):=M_{c(\o_{P_{l}})}\, T_{c(I_P)}\,D^2_{|I_P|}\,\phi(x)\:,
\eeq
or equivalently $ \phi_{P_l}(x)=e^{2\pi i c(\o_{P_{l}}) x}\,|I_{P}|^{-\frac{1}{2}}\,\phi(\frac{x-c(I_P)}{|I_P|})$.

We now define
\begin{itemize}
\item the $(y,\ll,\mu,\xi)$-\textbf{discretized Carleson real model} as
\beq\label{CMr}
 \tilde{\c}^{(y,\ll,\mu)}_{\xi}f(x):=\sum_{P\in\P^{y,\ll,\mu}} <f,\,\phi_{P_l}>\,\phi_{P_l}(x)\,\chi_{\o_{P_{u}}}(\xi)\:,
\eeq
where here $x\in\R$, $\xi\in\R$ and $f\in L^1(\R)$;

\item the $(y,\ll,\mu,\xi)$-\textbf{discretized Carleson periodic model} as
\beq\label{CMp}
 \tilde{C}^{(y,\ll,\mu)}_{\xi}f(x):=\sum_{P\in\P^{y,\ll,\mu,+}} <f,\,\phi_{P_l}>\,\phi_{P_l}(x)\,\chi_{\o_{P_{u}}}(\xi)\:,
\eeq
where here $x\in[0,1)$, $\xi\in\R$ and $f\in L^1(\R)$ with $\textrm{supp}\,f\subseteq[0,1)$.
\end{itemize}

With these being said, we have the following key result:
$\newline$

\noindent\textbf{Proposition.} In what follows, $c\in\R$ is an absolute constant that is allowed to change from line to line. The following are true:
\begin{itemize}

\item for any $\mu\in[0,\frac{1}{4}]$ we have for the real case
\beq\label{CM1}
\chi_{(-\infty,0]}(\xi)=c\,\int_0^1
 \sum_{m\in \Z} 2^{m+\ll}\sum_{P\in 2^{-m-\ll}[0,\,1)\times \mathcal{J}_{m}^{\ll,\mu}} |\hat{\phi}_{P_l}(\xi)|^2\,\chi_{\o_{P_{u}}}(0)\,d\ll\:,
\eeq
and for the periodic case
\beq\label{CM2}
\tilde{\chi}_{(-\infty,0]}(\xi)=c\,\int_0^1
 \sum_{m\in \N} 2^{m+\ll}\sum_{P\in 2^{-m-\ll}[0,\,1)\times \mathcal{J}_{m}^{\ll,\mu}} |\hat{\phi}_{P_l}(\xi)|^2\,\chi_{\o_{P_{u}}}(0)\,d\ll\:,
\eeq
where here $\tilde{\chi}_{(-\infty,0]}$ stands for a smooth version of $\chi_{(-\infty,0]}$ that is $\tilde{\chi}_{(-\infty,0]}\in C^{\infty}(\R)$ with $\textrm{supp}\,\tilde{\chi}_{(-\infty,0]}\subseteq (-\infty,0]$ and
$\tilde{\chi}_{(-\infty,0]}(\xi)=1$ for $\xi\leq -1$.

\item If $N\in\R$ is a fixed number then for the real case
\beq\label{CM3}
\chi_{(-\infty,N]}(\xi)=c\,\int_0^1\int_0^1
 \sum_{m\in \Z} 2^{m+\ll}\sum_{P\in 2^{-m-\ll}[0,\,1)\times \mathcal{J}_{m}^{\ll,\mu}} |\hat{\phi}_{P_l}(\xi)|^2\,\chi_{\o_{P_{u}}}(N)\,d\ll\,d\mu\:,
\eeq
while for the periodic case
\beq\label{CM4}
\tilde{\chi}_{(-\infty,N]}(\xi)=c\,\int_0^1\int_0^1
 \sum_{m\in \N} 2^{m+\ll}\sum_{P\in 2^{-m-\ll}[0,\,1)\times \mathcal{J}_{m}^{\ll,\mu}} |\hat{\phi}_{P_l}(\xi)|^2\,\chi_{\o_{P_{u}}}(N)\,d\ll\,d\mu\:,
\eeq
where as before $\tilde{\chi}_{(-\infty,N]}$ stands for a smooth version of $\chi_{(-\infty,N]}$.

\item for $f\in \mathcal{S}(\R)$, let us define the real axis Carleson operator
\beq\label{fulCarl}
\c f(x):=\sup_{N\in\Z}|\c_N f(x)|=\sup_{N\in\Z}\large|p.v.\,\int_{\R} e^{2\pi i N (x-y)}\,\frac{1}{x-y}\,f(y)\,dy\large|\,.
\eeq
Further, for $N\in \Z$, set
\beq\label{tilN}
\tilde{\c}_{N} f(x):=\int_{-\infty}^N \hat{f}(\xi)\,e^{2\pi i x \xi}\,d\xi\,,
\eeq
and notice that there exist $c_1,\,c_2\in\R$ absolute constants such that
\beq\label{translN}
\c_N f= c_1\,\tilde{\c}_{N} f\,+\,c_2\,f\:.
\eeq
As a consequence, one can reduce $\c$ to the study of
\beq\label{fulCarl1}
\tilde{\c}f(x):=\sup_{N\in\Z}|\tilde{\c}_{N} f(x)|:=\sup_{N\in\Z}\large|\int_{-\infty}^N \hat{f}(\xi)\,e^{2\pi i x \xi}\,d\xi\large|\,.
\eeq
Now, for $f,\,g\in\mathcal{S}(\R)$, one has
\beq\label{fulCarl2}
<\tilde{\c}_N f,\,g>=c\,\int_{0}^1\,\int_{0}^1\,\int_{0}^1
\sum_{P\in\P^{y,\ll,\mu}}  <f,\,\phi_{P_l}>\,<\phi_{P_l},\,g>\chi_{\o_{P_{u}}}(N)\,d\ll\,d\mu\,dy\,.
\eeq
Thus, linearizing the supremum, we deduce
\beq\label{fulCarl3}
\tilde{\c} f(x)=c\,\int_{0}^1\,\int_{0}^1\,\int_{0}^1 \sum_{P\in\P^{y,\ll,\mu}}  <f,\,\phi_{P_l}>\,\phi_{P_l}(x)\,\chi_{\o_{P_{u}}}(N(x))\,d\ll\,d\mu\,dy\,,
\eeq
or equivalently
\beq\label{fulCarl3e}
\tilde{\c} f(x)=c\,\int_{0}^1\,\int_{0}^1\,\int_{0}^1\tilde{\c}^{(y,\ll,\mu)}_{N(x)}f(x)\,d\ll\,d\mu\,dy\,.
\eeq
\item In what follows, we consider $f\in C_{0}^{\infty}(\R)$ with $\textrm{supp}\,f\subseteq[0,1)$ and $x\in[0,1)$. Let us recall the definition of the periodic Carleson operator:
\beq\label{fulCarlp}
C f(x):=\sup_{N\in\Z}|C_N f(x)|=\sup_{N\in\Z}\,\left|p.v.\,\int_{\TT}e^{i\,2\pi\,N\,(x-y)}\,\cot(\pi\,(x-y))\,f(y)\,dy\right|\,.
\eeq
Following the real case reasonings, for $N\in \Z$, one defines
\beq\label{fulCarl2p}
\tilde{C}_N f(x):=\int_{0}^1\,\int_{0}^1\,\int_{0}^1
\sum_{P\in\P^{y,\ll,\mu,+}}  <f,\,\phi_{P_l}>\,\phi_{P_l}(x)\,\chi_{\o_{P_{u}}}(N)\,d\ll\,d\mu\,dy\,.
\eeq
Then relation \eqref{translN} is replaced by
\beq\label{relp}
C_N f=c\,\tilde{C}_N+\A_{N}f\,,
\eeq
with
\beq\label{relp1}
\sup_{N}|\A_{N}f|\leq Mf\,,
\eeq
where here $M$ stands for the Hardy-Littlewood maximal operator.

Thus, setting
\beq\label{fulCarl3p}
\tilde{C} f(x)=\int_{0}^1\,\int_{0}^1\,\int_{0}^1
\sum_{P\in\P^{y,\ll,\mu,+}}  <f,\,\phi_{P_l}>\,\phi_{P_l}(x)\,\chi_{\o_{P_{u}}}(N(x))\,d\ll\,d\mu\,dy\,.
\eeq
or equivalently
\beq\label{fulCarl3ep}
\tilde{C} f(x)=\int_{0}^1\,\int_{0}^1\,\int_{0}^1 \tilde{C}^{(y,\ll,\mu)}_{N(x)}f(x)\,d\ll\,d\mu\,dy\,.
\eeq
based on \eqref{relp} and \eqref{relp1} one has that there exists $c\in\R$ such that
\beq\label{relp2}
|Cf-c\,\tilde{C}f|\leq Mf\,.
\eeq
\end{itemize}

\begin{o0}\label{WF} From \eqref{fulCarl3ep} and \eqref{relp2} we notice that\footnote{In what follows, for notational simplicity, we drop the subscript \textit{lac} from the definition of $C^{\{n_j\}_j}_{lac}$.}
\beq\label{keyCp}
\begin{array}{lc}
&\left|C^{\{n_j\}_j}f(x)-c\,\sup_{j}|\int_{0}^1\,\int_{0}^1\,\int_{0}^1
\sum_{P\in\P^{y,\ll,\mu,+}}  <f,\,\phi_{P_l}>\,\phi_{P_l}(x)\,\chi_{\o_{P_{u}}}(n_j)\,d\ll\,d\mu\,dy|\right|\\\\
&\leq Mf(x)\,,
\end{array}
\eeq
whenever $f\in C_{0}^{\infty}(\R)$ with $\textrm{supp}\,f\subseteq[0,1)$ and $x\in[0,1)$.

Compare now \eqref{keyCp} with \eqref{WalshC}. Deduce that \textbf{unlike} the Walsh-Carleson operator, the Carleson operator
is obtained as a \textbf{continuum average} of discrete models of type \eqref{CMp}. This will play a key role in explaining the differences between the a.e. pointwise behavior of the two operators.
\end{o0}

This ends the prerequisites about the basic definitions, concepts and properties regarding $C_{W}^{\{n_j\}_j}$ and $C^{\{n_j\}_j}$.

We pass now to a detailed motivation of the statements made in Observation \ref{W}.

From the proof of our Main Theorem 1, it should be by now obvious that there is no fundamental distinction between the a.e. pointwise behavior of $C^{\{n_j\}_j}$ and that of $C_{AW}^{\{n_j\}_j}$. Indeed, we immediately have the following

\begin{c0}\label{Walsh}[\textsf{An \textbf{``averaged"} Walsh-Carleson model}]\footnote{In the initial version of this paper, the corollary addressing the Walsh-Carleson model appeared label as Corollary 3 and was stated incorrectly. We thank Michael Lacey and the referee for pointing this to us. The subtle difficulties arising from the discretization of the Walsh-Carelson operator are now discussed in great details in the present section with the corresponding implications being summarized in Observation \ref{W}.}

Recall definition \eqref{carlac} of the (lacunary Fourier) Carleson operator. Given $N\in\N$ one can replace the Fourier mode $e^{i\,2\pi\,N\,x}$ with the corresponding Walsh mode $w_N(x)$ and so for a given sequence $\{n_j\}_{j\in\N}\subset\N$ one can define an averaged Walsh-Carleson model by
\beq\label{carlw}
C_{AW}^{\{n_j\}_j}f(x):=\sup_{j\in\N}\left|\int_{\TT} w_{n_j}(x)\,w_{n_j}(-y)\,\cot(\pi\,(x-y))\,f(y)\,dy\right|\:.
\eeq
(Here $\{w_{n_j}\}_j$ are regarded as periodic functions on $\R$.)

Then, the conclusions of Main Theorem 1, Main Theorem 2, Corollary \ref{Orlicz} and Corollary \ref{conjK} remain true for
$C_{AW}^{\{n_j\}_j}$.
\end{c0}

Indeed, one notices that Corollary  \ref{Walsh} requires only trivial modifications. In fact, the entire \textbf{CME} structure of $\F$ and the construction of the sets $\{F_j\}_j$ are left untouched. The only required modifications are expressed in the actual form of the (lacunary) averaged Walsh-Carleson operator $C_{AW}^{\{n_j\}_j}$ and part of the intermediate estimates provided in the last two sections.

To complete our antithesis and thus fully address Observation \ref{W}, it remains to discuss
\beq\label{ant}
\eeq
\begin{itemize}
\item why Theorem \ref{ww} holds for $C^{\{n_j\}_j}$;
\item why Theorem \ref{ww} \textsf{does not hold} for $C_{W}^{\{n_j\}_j}$ (and for the corresponding Lacey-Thiele discretized Carleson model).
\end{itemize}

We start our analysis with an easy but the same time important observation on the behavior of $C^{\{n_j\}_j}$ versus its discretized model(s):

Theorem \ref{ww} states that for any lacunary $\{n_j\}_j\subset\N$ we can find a sequence of functions $\{f_k\}_k\subset L^1(\TT)$ such that \eqref{normb} and \eqref{normub} hold and for some absolute constant $C>0$ one has
\beq\label{blowup1}
\|C^{\{n_j\}_j}(f_k)\|_{L^{1,\infty}}\geq C\,\|f_k\|_{L\,\log\log L\,\log\log\log\log L}\,.
\eeq

(As a consequence of the proof of Main Theorem 1 and of Remark 1 above Theorem \ref{ww} is proved.)

Now in view of \eqref{keyCp}, we thus know that \eqref{blowup1} can be rewritten as
\beq\label{blowup2}
\|\sup_{j}\,\left|\int_{0}^1\,\int_{0}^1\,\int_{0}^1 \tilde{C}^{(y,\ll,\mu)}_{n_j} (f_k)\,d\ll\,d\mu\,dy\right|\,\|_{1,\infty}\geq C'\,\|f_k\|_{L\,\log\log L\,\log\log\log\log L}\,.
\eeq
We now present the following

\begin{o0}\label{Fav}
Remark that while \eqref{blowup2} holds it is possible to find triples $(y,\ll,\mu)$ for which the corresponding \textbf{un-averaged} relation does \textbf{not} hold. In fact, we have that there exist triples $(y,\ll,\mu)$ and
there exists $c>0$ such that for any $f\in L^1(\TT)$ one has
\beq\label{blowup3}
\|\sup_{j}\,\left|\tilde{C}^{(y,\ll,\mu)}_{n_j} (f) \right|\,\|_{1,\infty}\leq c\,\|f\|_{L^1(\TT)}\,.
\eeq
\end{o0}

Indeed, in order to justify the above observation, let us pick $\{n_j=2^j\}_{j\in\N}$ and take a closer look to the discrete operators
$$\{\tilde{C}^{(y,0,0)}_{2^j}\}_{j\in\N}\:.$$
Given now definitions \eqref{inte} and \eqref{CMp} we deduce that
\beq\label{discj}
\tilde{C}^{(y,0,0)}_{2^j}f(x)=\sum_{P\in\P_{j+1}^{y,0,0,+}} <f,\,\phi_{P_l}>\,\phi_{P_l}(x)\,\chi_{\o_{P_{u}}}(2^j)\:.
\eeq
Notice now that this is a \textbf{single scale} operator, and hence we immediately deduce that for any $f\in L^1(\TT)$
\beq\label{cont}
\sup_{j}\,\left|\tilde{C}^{(y,0,0)}_{2^j} (f) \right|\lesssim M f(x)\,,
\eeq
and hence one has
\beq\label{cont1}
\|\sup_{j}\,\left|\tilde{C}^{(y,0,0)}_{2^j} (f) \right|\,\|_{1,\infty}\lesssim \|f\|_{L^1(\TT)}\,,
\eeq
explaining thus Observation \ref{Fav}.

Notice that exactly the same argument (with no modifications) can be carried over for the Walsh operator, thus giving
\beq\label{cont2}
\|\sup_{j}\,\left|C_{W}^{\{2^j\}_j} (f) \right|\,\|_{1,\infty}\lesssim \|Mf\|_{1,\infty}\lesssim\|f\|_{L^1(\TT)}\,.
\eeq

\begin{o0}\label{Fav1}
 Thus, we have identified a \textsf{first} reason of why Theorem \ref{ww} does not hold for $C_{W}^{\{n_j\}_j}$:

  relying on the second item displayed in \eqref{Walshp}, we have that the discrete nature of the Walsh operator makes possible the manifestation of a ``boundary effect" in which the specific choice of a sequence $\{n_j\}_j$ interacts with the specific choice of the dyadic grid appearing in the decomposition of the operator, allowing a single scale behavior at the frequencies defined by the sequence. In this way, the key geometric configurations of tiles defined by a \textbf{CME} \textsf{can not} be present in the time-frequency decomposition of $C_{W}^{\{n_j\}_j}$ since this structure requires a superposition of \textsf{multiple scales}\footnote{The number of such scales must tend to infinity as the parameter $k$ in \eqref{blowup2} goes to infinity.} at \textsf{each} frequency defined by - possibly a subsequence of - our lacunary sequence.
\end{o0}

Once reached in this point, one may ask the following natural

$\newline$
\noindent\textbf{Question:} \textit{Is it the true that the only way in which the Walsh analogue of \eqref{blowup1} can fail is when the time-frequency decomposition of $C_{W}^{\{n_j\}_j}$ lacks \textbf{CME} structures?}
$\newline$

\noindent\textbf{Answer:} \textit{No!}
$\newline$

Indeed, one can easily check that shifting the previous dyadic sequence, that is defining $n_j:= 2^j-1$, ($j\in\N$), one has:
\beq\label{ex}
\eeq
\begin{itemize}
\item the operator $C_{W}^{\{2^j-1\}_j}$ maps $L^1$ into $L^{1,\infty}$;

\item the time frequency decomposition of $C_{W}^{\{2^j-1\}_j}$ does admit \textbf{CME} structures.
\end{itemize}

In order to understand the general framework that allows the above answer (as well as the underlying mechanism behind the example provided by \eqref{ex}) we need to elaborate on a subtle element that plays a key role in the proof of our Main Theorem 1 and Theorem \ref{ww}.

Preserving the usual notations and in particular identifying $C^{\{n_j\}_j}$ with the operator $T$ used in the proof of our Main Theorem 1, suppose we have constructed our \textbf{CME} $\F=\bigcup_{j=\frac{k}{2}+1}^{k}\F_{j}$ with $\F_{j}=\bigcup_{l=1}^{\log k}\F_{j}^l$. When next constructing the corresponding sets
$\{F_{j}\}_{j=\frac{k}{2}+1}^{k}$ there are two important properties that we want to fulfill for each set $F_{j}$:
\beq\label{keyF}
\eeq
\begin{itemize}
\item for each $P\in \F_{j}$ the function $\textrm{Re}\,(\chi_{F_{j}}\,T_{P}^{*}(\chi_{[0,1]}))$ is positive;

\item for each $P\in \F_{j}^l$ one has $\frac{|I_P^{*}\cap F_{j}|}{|I_P|}\approx \frac{|I_P\cap F_{j}|}{|I_P|}\approx 2^{l}\,|F_j|\;.$
\end{itemize}
Decompose now $\F_{j}^l=\bigcup_{m} \F_{j}^{l,m}$ into maximal USGTF's. Next, fix such an $\F_{j}^{l,m}$ and choose $a\in\a(\F_{j}^{l,m})$; denote with $\F_{j}^{l,m}(a)$ all the tiles $P\in\F_{j}^{l,m}$ that live at frequency $a$.
Recall also that $\F_{j}^{l,m}[n]$ stands for those tiles $P\in \F_{j}^{l,m}$ having $A(P)=A_0(P)=2^{-n}$ where here $n\in\{2^{j-1}+\log\log k,\,2^j\}$. Set now $\F_{j}^{l,m}[n](a):=\F_{j}^{l,m}[n]\cap\F_{j}^{l,m}(a)$. Deduce then that from our \textbf{CME} construction we have that $\F_{j}^{l,m}[n](a)$ has precisely one element; moreover denoting with $P_{j}^{l,m}(a)$ the unique element in $\F_{j}^{l,m}(a)$ with $I_{P_{j}^{l,m}(a)}\in IBtm(\F_{j}^l)$ one has that
\beq\label{keyyF1}
\forall\:P\in\F_{j}^{l,m}(a)\:\:\:\:\:E(P)=E(P_{j}^{l,m}(a))\;.
\eeq

As a consequence of \eqref{keyF} we have the following important fact:

\beq\label{keyF1}
\begin{array}{ll}
&\textrm{if}\:\:P\in\F_{j}^{l,m}[n](a)\:\: \textrm{then}\\\\
&\textrm{Re}\,(\int\chi_{F_j}(x)\,T_{P}^{*}(\chi_{[0,1]})(x)\,dx)
\approx\int\chi_{F_j}(x)\,|T_{P}^{*}(\chi_{[0,1]})(x)|\,dx\\\\
&\approx |I_{P}^{*}\cap F_{j}|\,\frac{|E(P)|}{|I_P|}\approx_{\eqref{keyyF1}} 2^{l}\,|F_j|\,|E(P_{j}^{l,m}(a))|\:.
\end{array}\eeq
It is now the moment to include the above relation into the following

\begin{o0}\label{F} The proof of our Theorem \ref{ww} involving the operator $T\approx C^{\{n_j\}_j}$ relies crucially on the fact that the for each $P\in \F_{j}^{l,m}(a)$ the quantity $\textrm{Re}\,(\int\chi_{F_j}\,T_{P}^{*}(\chi_{[0,1]}))$ has \textbf{the same signum and approximative size} given by $2^{l}\,|F_j|\,|E(P_{j}^{l,m}(a))|$, a quantity which is \textbf{independent} (at least for our construction of USGTF's) on the scale and on the mass of the tile $P$.

As a consequence one immediately deduces that
\beq\label{keyF2}
|\sum_{P\in \F_{j}^{l,m}(a)} \int\chi_{F_j}\,T_{P}^{*}(\chi_{[0,1]})|\approx 2^j\,2^{l}\,|F_j|\,|E(P_{j}^{l,m}(a))|\:.
\eeq
\end{o0}

The analogue of \eqref{keyF2} for the Walsh case would read as
\beq\label{WkeyF2}
\begin{array}{lc}
&|\sum_{{R\in\mathcal{R}}\atop{R_u\in \F_{j}^{l,m}(a)}} <\chi_{F_j},\,w_{R_l}> \,<\chi_{[0,1]} w_{a}\,w_{R_l}\,\chi_{E(P_{j}^{l,m}(a))}>\,\chi_{\o_{R_u}}(a)|\\\\
&\approx 2^j\,2^{l}\,|F_j|\,|E(P_{j}^{l,m}(a))|\:.
\end{array}
\eeq

However the above relation \textsf{does not hold} for arbitrary values of $a\in\N$.

Indeed, taking as an example the lacunary sequence $\{n_r:= 2^r-1\}_{r\in\N}$ and choosing $a=a_0=2^r-1$ (for some large $r\in\N$), we
appeal to \eqref{alg1}-\eqref{Walshid} and \eqref{alg5} and deduce that
\beq\label{WkeyF22}
\begin{array}{lc}
&|\sum_{{R\in\mathcal{R}}\atop{R_u\in \F_{j}^{l,m}(a_0)}} <\chi_{F_j},\,w_{R_l}> \,<\chi_{[0,1]} w_{a_0}\,w_{R_l},\,\chi_{E(P_{j}^{l,m}(a_0))}>\,\chi_{\o_{R_u}}(a_0)|\\\\
&=|<W_{a_0}(\chi_{F_j})-W_{a_0-2^{2^{j-1}+\log\log k}}(\chi_{F_j}),\,\chi_{[0,1]} w_{a_0}\,\chi_{E(P_{j}^{l,m}(a_0))}>|\\\\
&\lesssim \sup_{P\in \F_{j}^{l,m}}\frac{\int_{I_{P}}\chi_{F_j}}{|I_{P}|}\,|E(P_{j}^{l,m}(a_0))|\lesssim 2^{l}\,|F_j|\,|E(P_{j}^{l,m}(a_0))|\:.
\end{array}
\eeq
This last fact violates \eqref{keyF2} for large enough $j\in\N$.

By simple modifications of \eqref{alg5}, one can see that \eqref{WkeyF22} continues to hold for many other choices of (frequencies within) lacunary sequences $\{n_r\}_r$ that have some suitable dyadic structure. Some examples of such lacunary sequences are given by $\{n_r=2^r+m\}_{r\in\N}$ (here $m\in\N$ fixed) or more generally any lacunary sequence  $\{n_r\}_{r\in\N}$ having the property that there exists $p\in\N$ such that for any $r\in\N$ the dyadic expansion of $n_r$ has at most $p$ nonzero terms\footnote{Of course, in this latter case the constant appearing in the second to last inequality of \eqref{WkeyF22} depends on $p$.}. Notice that all these examples of $\{n_r\}_r$ have the property that
\beq\label{Wa}
C_{W}^{\{n_r\}_r}\:\:\textrm{maps}\:\:L^1\:\:\textrm{into}\:\:L^{1,\infty}\,.
\eeq
We should mention here that the fact that \eqref{Wa} holds for the above examples of lacunary sequences has been known for some time (see e.g. \cite{ko3}). Indeed, based on an interesting observation of Konyagin (\cite{ko3}), one can find many other (classes of) sequences $\{n_j\}_j$ for which $C_{W}^{\{n_j\}_j}$ maps $L^1$ to  $L^{1,\infty}$ boundedly. Exactly because of this peculiar behavior, in his ICM address, Konyagin posed the following
$\newline$

\noindent\textbf{Open problem (OP I)}. (\cite{ko1}) \textit{Find a necessary and sufficient condition on a sequence $\{n_j\}_j\subseteq\N$ for which the associated Walsh-Carleson operator obeys
\beq\label{pat}
 C_{W}^{\{n_j\}_j}\:\:\:\:\:\textit{maps}\:\:\:L^1\:\:\:\textrm{to}\:\:\:L^{1,\infty}\,.
\eeq}

In view of the second item in \eqref{ant} and of the discussion following it, it is natural to ask if \eqref{pat} can only hold for rather ``exceptional" (lacunary) sequences as a manifestation of a ``boundary effect" due to the special \textsf{algebraic} and \textsf{dyadic/discrete} structures of the Walsh system. In this context, we raise the following

$\newline$
\noindent\textbf{Open Problem (OP II).} \textit{Decide if the Walsh analogue of Main Theorem 1 holds, that is, if it is true that in \eqref{blowup} $C_{lac}^{\{n_j\}_j}$ can be replaced by $C_{W}^{\{n_j\}_j}$.}
$\newline$

However, due to dichotomy \eqref{WkeyF2}-\eqref{WkeyF22}, we start now to better understand what type of obstacles we encounter when dealing with  Open Problems (OP) I and II. For example, a reasonable strategy in approaching OP II relies on the following

\begin{o0}\label{WW}
In order for the Main Theorem 1 to hold in the Walsh case (or equivalently, for our Open Problem II to be affirmatively decided) one needs to search for lacunary sequences $\{n_r\}_r$ for which
\begin{itemize}
\item there is a \textbf{CME} structure compatible with the time-frequency discretization of $C_{W}^{\{n_r\}_r}$;

\item once we identify a \textbf{CME}, relation \eqref{WkeyF2} must be satisfied.
\end{itemize}
\end{o0}
Moving our attention towards the Fourier case, we continue with the following

\begin{o0}\label{FF}
One should notice the following important antithesis
\begin{itemize}
\item it is possible to use a \textbf{single} (time-frequency) dyadic grid for decomposing the Carleson operator. Indeed, (proceeding as in Section 3) one can write for a generic operator $C_{lac}^{\{n_r\}_r}$ the following equality
 \beq\label{FE}
C_{lac}^{\{n_r\}_r}=\sum_{P\in \P} C_P\,.
\eeq
The key fact is that in this case $C_P$ are \textbf{not} projection operators of rank one, but they are convolution operators
with the form imposed by \eqref{TP}. Due to this specific form, once we fix a certain frequency $N(x)=a$, all the
adjoint operators $C_P^{*}$ with the property that $a\in\o_{P}$ will oscillate at the \textbf{same} frequency.

\item in contrast with the previous situation, one can use a decomposition of the Carleson operator involving \textbf{infinitely many} time-frequency ``dyadic" grids. In this latter case, in the sense described in
    \eqref{fulCarl3p} - \eqref{keyCp} one has
\beq\label{fullc1}
C f(x)\approx\int_{0}^1\,\int_{0}^1\,\int_{0}^1 \tilde{C}^{(y,\ll,\mu)}_{N(x)}f(x)\,d\ll\,d\mu\,dy\,
\eeq
 with each
 \beq\label{fullc2}
 \tilde{C}^{(y,\ll,\mu)}_{\xi}=\sum_{P\in\P^{y,\ll,\mu,+}} C_{P,\xi} \:,
 \eeq
where the operators $C_{P,\xi}$ are rank one projection operators defined by
\beq\label{fullc3}
C_{P,\xi}f = <f,\,\phi_{P_l}>\,\phi_{P_l}(x)\,\chi_{\o_{P_{u}}}(\xi) \:.
 \eeq
Notice that unlike the previous case, if we now fix the frequency $N(x)=a$, the
adjoint operators $\{C_{P,a}^{*}\}_{P\in \P^{y,\ll,\mu,+}}$ will oscillate at pairwise \textbf{distinct} frequencies.
\end{itemize}
\end{o0}

Now the entire proof of our Main Theorem 1, including the construction of a \textbf{CME}, is realised using a \textsf{single} dyadic time-frequency grid corresponding to the family of tiles $\P\approx\P^{0,0,0,+}$ and involving decomposition \eqref{FE} as described in the first item above. Taking now such a \textbf{CME} $\F$ and using the standard notations from before let us assume that $\max\{|I_P|\,|\,P\in \F_{j}^{l,m}(a)\}=2^{-m_0}$ and $\min\{|I_P|\,|\,P\in \F_{j}^{l,m}(a)\}=2^{-m_1}$. Then, following the same approach and in the same spirit as in \eqref{fulCarl3p} - \eqref{keyCp} we have that
\beq\label{keyFR}
\begin{array}{lc}
&\sum_{P\in \F_{j}^{l,m}(a)} C_{P}(\chi_{F_j})(x)\approx
\int_{0}^1\,\int_{0}^1\,\int_{0}^1\sum_{\bar{m}=m_0}^{m_1}\sum_{P\in\P_{\bar{m}}^{y,\ll,\mu,+}}\\\\
& <\chi_{F_j},\,\phi_{P_l}>\,\phi_{P_l}(x)\,\chi_{E(P_{j}^{l,m}(a))}(x)\,e^{-2\,\pi i\,a\,x}
\,\chi_{\o_{P_{u}}}(a)\,d\ll\,d\mu\,dy\:.
\end{array}
\eeq
Notice that while the RHS of \eqref{WkeyF2} corresponds in the Fourier setting to
\beq\label{discj1}
\sum_{P\in\P^{0,0,0,+}\atop{P\in \F_{j}^{l,m}(a)}} <\chi_{F_j},\,\phi_{P_l}>\,<\phi_{P_l}(\cdot)\,e^{-2\,\pi i\,a\cdot},\,\chi_{E(P_{j}^{l,m}(a))}(\cdot)\,\chi_{[0,1]}(\cdot)>\,\chi_{\o_{P_{u}}}(a)\:,
\eeq
after the \textsf{averaging process} we get using \eqref{keyFR} and \eqref{keyF2} that
\beq\label{keyFRR}
\begin{array}{lc}
&\sum_{P\in \F_{j}^{l,m}(a)} \int C_{P}(\chi_{F_j})\,\chi_{[0,1]}\approx\int_{0}^1\,\int_{0}^1\,\int_{0}^1\sum_{\bar{m}=m_0}^{m_1} \sum_{P\in\P_{\bar{m}}^{y,\ll,\mu,+}} \\\\
&<\chi_{F_j},\,\phi_{P_l}>\,<\phi_{P_l}(\cdot)\,e^{-2\,\pi i\,a\cdot},\,\chi_{E(P_{j}^{l,m}(a))}(\cdot)\,\chi_{[0,1]}(\cdot)>\,\chi_{\o_{P_{u}}}(a)\,d\ll\,d\mu\,dy\\\\
&\approx 2^j\,2^{l}\,|F_j|\,|E(P_{j}^{l,m}(a))|\:.
\end{array}
\eeq

Thus we can now conclude with the following

\begin{o0}\label{Fav2}
 We now understand the \textsf{second} reason of why Theorem \ref{ww} does not hold for $C_{W}^{\{n_j\}_j}$:

 while necessary, the existence\footnote{Relative to the time-frequency decomposition of $C_{W}^{\{n_j\}_j}$.} of \textbf{CME} structures is not sufficient for the Walsh analogue of \eqref{blowup1}. This is due to the \textsf{combined} effect of two facts:
\begin{itemize}
\item  the \textbf{CME} structure has \textsf{different frequency implications} in the Fourier (``continuous") case and in the Walsh (``discrete") case. Indeed, as hinted in Observation \ref{FF}, in the first case, all the operators $\{C_P^{*}\}_{P\in \F_{j}^{l,m}(a)}$ oscillate at the \textsf{same} frequency $a$. In the second case, the operators $\{W_R^{*}\}_{R_u\in \F_{j}^{l,m}(a)}$ defined by $W_{R}f(x):=<f,w_{R_l}>w_{R_l}(x)\,\chi_{\o_{R_u}}|_{N(x)=a}$ will oscillate at pairwise \textsf{distinct} frequencies. Notice that the latter property is not due to the particular (algebraic) structure of the Walsh wavepackets but is due to the \textsf{geometry} of the tile discretization: if $\tilde{\mathcal{R}}$ is any family of bi-tiles with the property that there exists $a\in\R$ such that $R\in\tilde{\mathcal{R}}$ implies $a\in \o_{R_u}$ then any two tiles within the family $\{R_l\}_{R\in\tilde{\mathcal{R}}}$ are pairwise disjoint.
\item based on the algebraic properties of the Walsh wave packets (such as \eqref{alg2} - \eqref{alg5}) one can find counterexamples to the Walsh analogue of \eqref{keyF2}. However, as described in Observation \ref{F}, relation \eqref{keyF2} is of key importance in the proof of Theorem \ref{ww}.
\end{itemize}
Finally, we notice that taking averages over discrete models makes the pointwise behavior of the lacunary Carleson operator $C^{\{n_j\}_j}$ independent on the structure of the lacunary sequence $\{n_j\}_j$.
\end{o0}
$\newline$

\section{Final remarks.}

\noindent 1) We start our section with a discussion on Theorem \ref{ww}. In what follows we will briefly present the relevant modifications in the proof of our Main Theorem 1 that are required in order to deal with a general lacunary sequence. The key observation is that the mechanism involved in constructing the \textbf{CME} $\F$ is \textsf{independent} on the specific choice of our lacunary sequence and the only thing that depends on $\{n_j\}_j$, and hence needs special care, is the construction of the sets $\{F_j\}_{j=\frac{k}{2}+1}^{k}$. \footnote{This should be integrated in the lines of thoughts supporting Observation \ref{Cme}.} Indeed, it is natural to expect that the structure of the frequencies plays a role in the corresponding structure of the sets. In the case of a perfectly dyadic sequence, once we constructed
the multi-tower $\F=\bigcup_{j=\frac{k}{2}+1}^{k}\F_j$ and decomposed each $\F_j=\bigcup_{l=1}^{\log k}\F_j^l$ and then each tower $\F_j^l$ into maximal USGTF's $\bigcup_{m}\F_j^{l,m}$, taking for simplicity $j=k$, we could arrange for \eqref{monot} to hold for each $I\in\I Btm (\F_k^{\log k})$ and thus in turn we were able to construct $F_k$ to obey \eqref{REQ}.
This further implied that
\beq\label{ce}
\begin{array}{ll}
&\textit{for any}\: I\in\I Btm (\F_k^{\log k})\:\textit{the set}\:I\cap F_k\:\\
&\textit{has a fractal structure (Cantor set) of the same size}\;.
\end{array}
\eeq

 For an arbitrary lacunary sequence $\{n_j\}_j$ we no longer aim for preserving the \textsf{exact} form of \eqref{ce}. Instead, we seek for a good \textsf{approximation} of \eqref{ce}, meaning that we want the sets $I\cap F_k$ to have fractal structure with ``almost" the same (relative) size as $I$ ranges through $\I Btm (\F_k^{\log k})$. The precise meaning of the above heuristic is explained in what follows:

Let $\{n_j\}_j$ be our favorite choice for a lacunary sequence.

Fix as before $k=2^{2^K}$ for some large $K\in \N$. We want to construct the sets $\{F_j\}_{j=\frac{k}{2}+1}^k$ such that the function
$$f_k:=\frac{1}{k}\sum_{j=\frac{k}{2}+1}^{k} \frac{1}{|F_j|\,\log\log \frac{1}{|F_j|}}\,\chi_{F_j}$$
obeys \eqref{blowup}.

Now since $C_{lac}^{\{n_j\}_j}$ is a maximal operator and we are interested in a lower bound for its $L^{1,\infty}$ norm, we can always restrict our attention to any subsequence of the initial $\{n_j\}_j$. Taking now $\bar{K}:=2^{2^{2^k}}$ and by possibly passing to a subsequence, we can wlog assume that $n_1>1$ and that the following holds:
\beq\label{nj}
\frac{n_{j+1}}{n_j}> 2^{\bar{K}}\:\:\:\:\:\;\:\:\forall\:j\in\N\:.
\eeq
We now adapt \eqref{imN} to our new context by requiring
\beq\label{imN1}
  \textsf{Image}(N)\subseteq\{n_{\bar{K}+m}\}_{m\in \{0,\,\log k\,2^{2^k-1}\}}\:.
\eeq

With this done, we follow line by line the construction of our \textbf{CME} $\F:=\bigcup_{l=\frac{k}{2}+1}^{k}\F_{l}$ presented in Section 7 with the only trivial change regarding the frequency locations of our tiles in each $USGTF$, that is, whenever we see a frequency $\a(P)=2^{2^{2^{2^k}}}+100m$ we replace it with the corresponding analogue $\a(P)=n_{\bar{K}+m}$.

Having constructed our \textbf{CME} we now adapt the construction of our sets $\{F_j\}_{j=\frac{k}{2}+1}^{k}$ in Section 8.1. to our new setting. For simplicity and space constrains we only focus on the case $j=k$. As in Section 8.1., we fix $I\in\I Btm (\F_k^{\log k})$ and define $\S[I,\F_{k}^{1}]$ and $\U_{\F_k^{1}}[I](a)$ exactly as in  \eqref{zero}-\eqref{u0}.

 Remark that  \eqref{monot} ceases to remain true; however, based on the fact that for any $P\in\F$ (in particular $A_0(P)\leq \frac{1}{2}$) one has $|I_P|> 2^{-2^{2k}}$ and hence $\frac{1}{n_j}<<|I|$ for any $j\geq \bar{K}$, we deduce that
 \begin{itemize}
 \item for all $a\in \a(\F_{k})$
 \beq\label{u11}
(\frac{1}{2}-\frac{1}{2^{\frac{\bar{K}}{2}}})\,|I|\leq |\U_{\F_k^1}[I](a)|\leq (\frac{1}{2}+\frac{1}{2^{\frac{\bar{K}}{2}}})\,|I|\:.
\eeq
 \item for all $a,\,b\in \a(\F_{k}^{1})\:\:\:\textrm{with}\:\;\:a<b$
\beq\label{monot1}
(\frac{1}{2}-\frac{1}{2^{\frac{\bar{K}}{2}}})\,|\U_{\F_k^1}[I](a)|\leq|\U_{\F_k^1}[I](a)\cap \U_{\F_k^1}[I](b)|
\leq (\frac{1}{2}+\frac{1}{2^{\frac{\bar{K}}{2}}})\,|\U_{\F_k^1}[I](a)|\:.
\eeq
\end{itemize}
 Notice now that \eqref{monot1} becomes a very good approximation of \eqref{monot}.

With these being said, and following Section 8.1., we impose that the set $F_k$ obeys in a first instance \eqref{reqFk1} and then the more general \eqref{reqFk2}. In fact, we can now slightly simplify the initial dyadic scenario and directly define $F_k$ as the set obeying:
\begin{itemize}
\item $F_k=\bigcup_{I\in \I Btm (\F_k^{\log k})} I\cap F_k$;

\item $I\cap F_k:=\bigcap_{l=1}^{\log k}\:\bigcap_{a\in \a(\F_{k}^{l})} \U_{\F_k^l}[I](a)$.
\end{itemize}

With this done, and making essential use of \eqref{u11} and \eqref{monot1} we get the analogue of \eqref{imp}, that is
\beq\label{impk}
(\frac{1}{2}-\frac{1}{2^{\frac{\bar{K}}{2}}})^{\#\a(\F_k)}\leq\frac{|I\cap F_k|}{|I|}\leq (\frac{1}{2}+\frac{1}{2^{\frac{\bar{K}}{2}}})^{\#\a(\F_k)}\,,
\eeq
 which from the choice of $\bar{K}$ further implies that
\beq\label{impk1}
 \frac{1}{e}\,\frac{1}{2^{\log k\,2^{2^{k}-1}}}\leq \frac{|I\cap F_k|}{|I|}\leq e\,\frac{1}{2^{\log k\,2^{2^{k}-1}}}\:.
\eeq

In particular, since $|\tilde{\I} Btm (\F_k^{\log k})|\approx 2^{-\log k}=\frac{1}{k}$, we get that
\beq\label{impk2}
 |F_k|\approx\frac{1}{k}\,2^{-\log k\,2^{2^{k}-1}}\:.
\eeq

Notice now that from the above construction one gets another key property that our set $F_k$ must satisfy, that is
\beq\label{prp}
 \forall\:l\in\{1,\ldots,\log k\}\:\:\textrm{and}\:\:\forall\:P\in\F_{k}^{l}\:\:\Rightarrow\:\:\frac{|I_P\cup F_k|}{|I_P|}\approx 2^{l}\,|F_k|\:.
\eeq

The above process can be repeated with the obvious changes for a general $j\in\{\frac{k}{2}+1,\,k\}$ and in this situation \eqref{impk1} - \eqref{prp} become
\beq\label{impk12}
 \frac{1}{e}\,\frac{1}{2^{\log k\,2^{2^{j}-1}}}\leq \frac{|I\cap F_j|}{|I|}\leq e\,\frac{1}{2^{\log k\,2^{2^{j}-1}}}\:.
\eeq
\beq\label{impk22}
 |F_j|\approx\frac{1}{k}\,2^{-\log k\,2^{2^{j}-1}}\:.
\eeq
and respectively
\beq\label{prp2}
 \forall\:l\in\{1,\ldots,\log k\}\:\:\textrm{and}\:\:\forall\:P\in\F_{j}^{l}\:\:\Rightarrow\:\:\frac{|I_P\cup F_j|}{|I_P|}\approx 2^{l}\,|F_j|\:.
\eeq

Thus at the end of this process one is able to construct the desired sequence of sets $\{F_j\}_{j=\frac{k}{2}+1}^k$.

Finally, based again on the properties \eqref{nj}, \eqref{u11} and \eqref{monot1} one can easily check that the fundamental Claim in Section 8.3., more precisely \eqref{KEY}, holds.

With this done, the reasonings in Sections 9 and 10 can be repeated concluding the proof of our Theorem \ref{ww}.
$\newline$

\noindent 2) Recalling the description of our \textbf{CME}
$\F=\bigcup_{j=\frac{k}{2}+1}^{k}\F_j$  (see Section 7), this remark seeks to explain \textit{why} we chosen the height of each tower $\F_j$ to be of the order $\log k$.  As one can notice, writing $\F_j=\bigcup_{l=1}^{m}\F_{j}^{l}$, a first impulse would be to aim for a height $m$ as big as possible, (\textit{e.g.} $m=(\log k)^2$), hoping that this would increase the lower bound of the $L^{1,\infty}$ norm of $C_{lac} f_k$. However, in the view of the previous remark we immediately notice that this is not the case. Indeed, based on relations \eqref{JN}, \eqref{use} and  \eqref{obs} and on the fact that
\beq\label{base}
\{\n_j\geq m\}\subset Basis (\F_j^m)\subset \{\n_j>\frac{m}{2}\}\,,
\eeq
we notice that there exists $c>0$ absolute constant such that if $m>C\,\log k$ with $C>0$ absolute constant large enough then
\beq\label{base1}
\sum_{j=\frac{k}{2}+1}^{k} |Basis (\F_j^m)|\leq k\,e^{-c m} = o(k)\,,
\eeq
and thus the set $\bigcup_{j=\frac{k}{2}+1}^{k} Basis (\F_j^m)$ becomes an ``exceptional" (removable) set that consequently has a negligible impact on the  size of $\|C_{lac} f_k\|_{1,\infty}$. From this, one further notices the strong connection between the tile structure maximizing the size of $\|C_{lac} f_k\|_{1,\infty}$ and the one corresponding to the maximal size of $\|\n^{[k]}\|_{1,\infty}$.
$\newline$

\noindent 3) Corollaries  \ref{Orlicz}, \ref{conjK}, \ref{Halo} and \ref{Limextr} are straightforward applications of the first part of Main Theorem 2.
$\newline$

\noindent 4) For Corollary \ref{revista}, preserving the notations in Section 7, one simply takes $\P^{\a}:=\F_{k}^1$ and $f=\chi_{F_k}$. We leave further details to the interested reader.
$\newline$

\noindent 5) Our entire paper relies on a completely new idea revealing the deep relationship between the behaviour of the grand maximal counting function and the pointwise convergence of Fourier series near $L^1$. This is the first result in the math literature regarding this topic that provides a counterexample using multiscale analysis of the wavepackets. All the previous counterexamples, were focusing on identifying and working directly with input functions and on doing involved computations addressing the local pointwise estimates of (sub)sequences of partial Fourier sums applied to such input functions. In our case, we change completely the point of view, by first designing a geometric construction of tiles that are encoding both the nature of the lacunary sequence (reflected in their frequencies) and the extremizer property relative to the $L^{1,\infty}$ norm of the grand maximal counting function. The input functions are now naturally obtained as a direct byproduct of this construction, relying on the alignment of the oscillations requirement.

Thus the entire realization of this paper is based on a new perspective. It is not about a technical $\log \log\log\log L$ factor addition but is about a \textsf{conceptual} advancement that identifies a \textsf{structural} mechanism in approaching the problem of the pointwise convergence of Fourier series near $L^1$.
$\newline$

\noindent\textbf{Acknowledgement.} I would like to thank Charlie Fefferman and Terence Tao for reading parts of the manuscript and providing useful comments. Also, I would like to express my gratitude to Alex Stokolos for referring me to paper \cite{PH} and Fernando Soria for several bibliographical suggestions (including \cite{ch} and \cite{P}) and a helpful discussion on the Halo conjecture. Finally, I want to thank the anonymous referees whose comments and remarks greatly improved the presentation of our paper.

\section{Appendix}

In this section we review some of the basic facts concerning the theory of rearrangement invariant (quasi-)Banach spaces.
We follow closely the description made in \cite{MMR} which further relies on \cite{BeSa}.

Denote with $L^{0}(\TT)$ the topological linear space of all periodic Lebesgue-measurable functions equipped with the topology of convergence in measure. Given now $f\in L^{0}(\TT)$, we define its \textit{distribution function} as
\beq\label{distribf}
m_{f}(\ll):=m(\{x\in\TT\,|\,|f(x)|>\ll\})\,,
\eeq
where here $m$ stands for the Lebesgue measure on $\TT$.

The \textit{decreasing rearrangement} of $f$ is now defined as
\beq\label{decrrea}
f^{*}(t):=\inf\{\ll\geq 0\,|\,m_f(\ll)\leq t\}\,,\:\:\;\:\:t\geq 0\;.
\eeq

All the quasi-Banach spaces $X$ mentioned in our paper are considered as subspaces of $L^{0}(\TT)$.

We say that $X$ is a \textit{(quasi-)Banach lattice} if the following properties are satisfied:

- $\exists\:h\in X$ with $h>0\:\:a.e.$;

- if $|f|\leq |g|$ a.e. with $g\in X$ and $f\in L^0(\TT)$ then $\|f\|_{X}\leq \|g\|_X$.

A (quasi-)Banach lattice $(X,\,\|\cdot\|_X)$ is called a \textit{rearrangement invariant (quasi-)Banach space} iff
given $f\in X$ and $g\in L^0(\TT)$ with $m_f=m_g$ one has $g\in X$ and $\|f\|_X=\|g\|_X$.

If $X$ is a r.i. (quasi-)Banach space and $\chi_{A}$ stands for the characteristic functions of a measurable set $A\subseteq\TT$, the function
\beq\label{fundamfunction}
\v_{X}(t)=\|\chi_A\|_{X}\;\:\:\:\:\:\textrm{with}\:\:\:\:m(A)=t\:\:\:\textrm{and}\:\:\:t\in[0,1]\;,
\eeq
is called the \textit{fundamental function} of $X$.

In what follows we introduce two fundamental classes of r.i. Banach spaces: the Marcinkiewicz and the Lorentz spaces.

We say that $\v:\,\TT\,\rightarrow\,\R_{+}$ is a \textit{quasi-concave function} if the following hold:

- $\v(0)=0$ and $\v(t)>0$ for all $t\in (0,1)$;

- the functions $\v(t)$ and $\v_{*}(t)=\frac{t}{\v(t)}$ are non-decreasing functions on $\TT$.

It is worth noticing that, for our purposes here, we can always replace a quasi-concave function $\v$ with its least concave majorant $\tilde{\v}$  since one has
$\tilde{\v}(t)\leq 2\v(t)\leq 2\tilde{\v}(t)$ for all $t\in\TT$.

The \textit{Marcinkiewicz space} $M_{\v}$ is the r.i. Banach space of all $f\in L^{0}(\TT)$ such that
\beq\label{Marc}
\|f\|_{M_{\v}}:=\sup_{t\in(0,1]} \frac{1}{\v(t)}\,\int_{0}^{t} f^{*}(s)\,ds <\infty\;.
\eeq

The setting for defining a Lorentz space always requires first defining a function $\v:\,\TT\,\rightarrow\,[0,\infty)$ with
the following properties:
\beq\label{proplorentz}
\eeq
\begin{itemize}
\item $\v(0)=0$;
\item $\v$ is nondecreasing;
\item $\v$ is concave.
\end{itemize}
With $\v$ as above, we define the \textit{Lorentz space} $\l_{\v}$ as the r.i. Banach space of all the functions $f\in L^{0}(\TT)$ such that
\beq\label{Lorentz}
\|f\|_{\l_{\v}}:=\int_{0}^{1} f^{*}(s)\,d(\v(s)) <\infty\;.
\eeq
Finally, we close with the following important observation: if $X$ is a r.i. Banach space with the fundamental function $\v$, then
\begin{itemize}
\item $\v$ is quasi-concave;

\item the following continuous inclusion holds:
\beq\label{continclusion}
\l_{\tilde{\v}}\,\hookrightarrow\, X\,\hookrightarrow\,M_{\v_{*}}\:.
\eeq
\end{itemize}

\end{document}